\documentclass[12pt]{amsart}
\usepackage{amsmath}
\usepackage{amsfonts}
\usepackage{amssymb}
\usepackage{amscd}
\usepackage{graphicx}
\usepackage[abbrev,alphabetic]{amsrefs}
\usepackage{hyperref}
\RequirePackage[dvipsnames,usenames]{color}
\usepackage{stmaryrd}
\usepackage{mathtools}
\usepackage{booktabs}
\usepackage{multirow}
\usepackage{stmaryrd} 
\newtagform{tiny}{\tiny(}{)}
\usepackage{tikz}
\usepackage{tikz-3dplot}

\makeatletter
\@addtoreset{section}{part}
\makeatother

\usepackage{mathtools}
\usepackage{hyperref}
\usepackage[margin=1.25in]{geometry}

\usepackage{amsthm}
\usepackage{comment}
\usepackage[all,cmtip]{xy}
\usepackage{tikz-cd}
\usetikzlibrary{cd}
\usepackage{mathtools}
\usepackage[all]{xy}

\makeatletter
\def\@tocline#1#2#3#4#5#6#7{\relax
  \ifnum #1>\c@tocdepth 
  \else
    \par \addpenalty\@secpenalty\addvspace{#2}%
    \begingroup \hyphenpenalty\@M
    \@ifempty{#4}{%
      \@tempdima\csname r@tocindent\number#1\endcsname\relax
    }{%
      \@tempdima#4\relax
    }%
    \parindent\z@ \leftskip#3\relax \advance\leftskip\@tempdima\relax
    \rightskip\@pnumwidth plus4em \parfillskip-\@pnumwidth
    #5\leavevmode\hskip-\@tempdima
      \ifcase #1
       \or\or \hskip 1em \or \hskip 2em \else \hskip 3em \fi%
      #6\nobreak\relax
    \hfill\hbox to\@pnumwidth{\@tocpagenum{#7}}\par
    \nobreak
    \endgroup
  \fi}
\makeatother

\newsavebox{\pullback}
\sbox\pullback{%
\begin{tikzpicture}%
\draw (0,0) -- (1ex,0ex);%
\draw (1ex,0ex) -- (1ex,1ex);%
\end{tikzpicture}}

\newsavebox{\pullbackdl}
\sbox\pullbackdl{%
\begin{tikzpicture}%
\draw (-1ex,0ex) -- (0ex,0ex);%
\draw (0ex,-1ex) -- (0ex,0ex);%
\end{tikzpicture}}

\newsavebox{\pushoutdr}
\sbox\pushoutdr{%
\begin{tikzpicture}%
\draw (-1ex,-1ex) -- (-1ex,0ex);%
\draw (-1ex,0ex) -- (0ex,0ex);%
\end{tikzpicture}}

\newcommand{\cyan}{\color{cyan}}

\newcommand{\cred}{\color{red}}

\renewcommand{\mod}{\ \textrm{mod}\ }
 
\newcommand{\A}{\mathbb{A}}
\renewcommand{\P}{\mathbb{P}}
\newcommand{\Z}{\mathbb{Z}}
\newcommand{\Q}{\mathbb{Q}}

\newcommand{\R}{\mathbb{R}}
\newcommand{\F}{\mathbb{F}}

\newcommand{\G}{\mathbb{G}}
\newcommand{\Cris}[0]{{\operatorname{cris}}}

\newcommand{\cF}{\mathcal{F}}

\newcommand{\MO}{\mathcal{O}}
\newcommand{\sO}{\mathcal{O}}

\newcommand{\m}{\mathfrak{m}}

\newcommand{\Frac}{\mathrm{Frac}}
\newcommand{\Proj}{\mathrm{Proj}}
\renewcommand{\Im}{\mathrm{Im}}
\newcommand{\wt}{\widetilde}
\newcommand{\ol}{\overline}
\newcommand{\et}[0]{{\operatorname{\acute{e}t}}}

\DeclareMathOperator{\Br}{Br}
\DeclareMathOperator{\Gr}{Gr}
\DeclareMathOperator{\CGr}{CGr}
\DeclareMathOperator{\Gal}{Gal}
\DeclareMathOperator{\ch}{ch}

\DeclareMathOperator{\Sing}{Sing}
\DeclareMathOperator{\rank}{rank}
\DeclareMathOperator{\cris}{cris}

\DeclareMathOperator{\Supp}{Supp}
\DeclareMathOperator{\Spec}{Spec}

\DeclareMathOperator{\Ext}{Ext}

\DeclareMathOperator{\GL}{GL}
\DeclareMathOperator{\Pic}{Pic}

\DeclareMathOperator{\Ex}{Ex}
\DeclareMathOperator{\Cone}{Cone}

\DeclareMathOperator{\Ker}{Ker}

\DeclareMathOperator{\Bl}{Bl}
\DeclareMathOperator{\dlog}{dlog}
\DeclareMathOperator{\la}{\langle}
\DeclareMathOperator{\ra}{\rangle}

\theoremstyle{plain}
\newtheorem{theorem}{Theorem}[section]
\newtheorem{thm}[theorem]{Theorem}

\newtheorem{prop}[theorem]{Proposition}

\newtheorem{lem}[theorem]{Lemma}

\newtheorem{cor}[theorem]{Corollary}

\newtheorem*{claim*}{Claim}
\newtheorem{step}{Step}

\theoremstyle{definition}

\newtheorem{dfn}[theorem]{Definition}

\newtheorem{nota}[theorem]{Notation}

\newtheorem*{setup*}{Setup}
\newtheorem{theo}{Theorem}

\newtheorem{nasi}[theorem]{}

\theoremstyle{remark}
\newtheorem{rem}[theorem]{Remark}



\numberwithin{equation}{theorem}

\title[Liftability and vanishing for Fano threefolds]{Liftability and vanishing theorems for Fano threefolds in positive characteristic I}

\author{Tatsuro Kawakami}
\address{Graduate School of Mathematical Sciences, 
The University of Tokyo, 
3-8-1 Komaba, Meguro-ku, Tokyo 153-8914, JAPAN} 
\email{kawakami@ms.u-tokyo.ac.jp}
\author{Hiromu Tanaka} 
\address{Department of Mathematics, 
Graduate School of Science, 
Kyoto University, 
Kyoto 606-8502, JAPAN} 
\email{tanaka.hiromu.7z@kyoto-u.ac.jp}

\begin{document}

\begin{abstract}
In our series of papers, we prove that smooth Fano threefolds in positive characteristic lift to the ring of Witt vectors. 
Moreover, we show that they 
satisfy Akizuki-Nakano vanishing, 
$E_1$-degeneration of the Hodge to de Rham spectral sequence, 
and torsion-freeness of Crystalline cohomologies. 
In this paper, we establish these results for the case when $|-K_X|$ is very ample and the Picard group is generated by $\omega_X$. 
\end{abstract}

\maketitle

\setcounter{tocdepth}{2}

\tableofcontents

\section{Introduction}

In the context of the minimal model program, Fano varieties play a significant role in the classification of algebraic varieties. 
The classification of Fano varieties in characteristic zero has a long history. 
In the early nineteenth century, 
Gino Fano established 
a partial classification result for smooth Fano threefolds, 
in order to attack the rationality problem of cubic threefolds. 
In 1980s, the classification of smooth Fano threefolds was carried out 
by Mori--Mukai (\cite{MM81}, \cite{MM83}), based on earlier works by Iskovskih and Shokurov. 
Recently, this classification result has been extended to the case of positive characteristic \cites{
Meg98, FanoI, FanoII, FanoIII, FanoIV}.

From now on, we work over an algebraically closed field $k$ of characteristic $p>0$.
Both for geometric and arithmetic purposes, it is important to ask when a smooth projective variety over $k$ lifts 
to the ring $W(k)$ of Witt vectors, in other words, when $X$ comes from a smooth projective variety in characteristic zero. 
As a main theorem of our series of papers, we prove that every smooth Fano threefold lifts to $W(k)$.

\begin{theo}[Liftability to $W(k)$]\label{Introthm:W(k)-lift}
Let $X$ be a smooth Fano threefold over an algebraically closed field $k$ of positive characteristic.
Then $X$ lifts to $W(k)$.
\end{theo}

The liftability provides us useful consequence such as Akizuki-Nakano vanishing and $E_1$-degeneration of Hodge to de Rham spectral sequence when $p\geq \dim X = 3$ \cite{DI}. 
We shall establish these results in arbitrary characteristic.

\begin{theo}[Akizuki-Nakano vanishing]\label{Introthm:ANV}
Let $X$ be a smooth Fano threefold over an algebraically closed field $k$ of characteristic $p>0$. 
Then Akizuki-Nakano vanishing holds on $X$, that is, 
if $A$ is an ample  Cartier divisor $A$ on $X$, then we have
\[
H^j(X,\Omega^i_X \otimes \MO_X(-A))=0
\]
for all integers $i,j\geq 0$ satisfying  $i+j<3$.
\end{theo}

\begin{theo}
\label{Introthm:E_1-degeneration}
Let $X$ be a smooth Fano threefold over an algebraically closed field $k$ of characteristic $p>0$. 
Then the following hold. 
\begin{enumerate}
    \item The Hodge to de Rham spectral sequence 
\[
E_1^{i,j} 
= H^j(X, \Omega_X^i) 
\Rightarrow 
H^{i+j}(X, \Omega_X^{\bullet}) =E^{i+j}
\]
degenerate at $E_1$.
    \item Crystalline cohomology $H^i_{\cris}(X/W(k))$ is torsion-free for every $i\geq 0$.
\end{enumerate}
\end{theo}

\begin{rem}\label{r ANV Intro}
We here summarise some previous works on Akizuki-Nakano vanishing for smooth Fano threefolds (Theorem \ref{Introthm:ANV}). 
\begin{enumerate}
\item  $H^1(X, \sO_X(-A))=0$ 
and 
$H^0(X, \Omega_X^1 \otimes \MO_X(-A))=0$  
have been already proven 
(see \cite{SB97}, \cite[Theorem 1.1]{Kaw2}, \cite[Theorem 2.4]{FanoI}). 
\item 
$H^2(X, \sO_X(-A))=0$ has been proven under the additional assumption that $\rho(X)=r_X=1$
\cite[Corollary 4.5]{FanoI}.
\end{enumerate}
\end{rem}

Note that Theorem \ref{Introthm:E_1-degeneration} is a formal consequence 
of $W(k)$-liftablity (Theorem \ref{Introthm:W(k)-lift}) and 
invariance of Hodge numbers (Theorem \ref{Introthm:hoge number}).   

\begin{theo}[Hodge numbers]\label{Introthm:hoge number}
Let $X$ be a smooth Fano threefold over an algebraically closed field $k$ of positive characteristic. 
Take a lift $f\colon \mathcal{X}\to W(k)$ of $X$ to $W(k)$, whose existence is ensured by Theorem \ref{Introthm:W(k)-lift}. 
Let $X_{\overline{K}}$ be the geometric generic fibre over $W(k)$. 
Then all the Hodge numbers $h^{i,j}(X)\coloneqq \dim_{k} H^j(X,\Omega^i_X)$ of $X$ coincide with those of $X_{\overline{K}}$, that is, 
\[
h^j(X,\Omega^i_X)=h^j(X_{\overline{K}},\Omega^i_{X_{\overline{K}}})
\]
hold for all $i,j\geq 0$.
\end{theo}

{Here, we remark that Petrov \cite{Petrov} recently proved that a smooth projective quasi-$F$-split variety satisfies Akizuki-Nakano vanishing (hence lift to $W(k)$ if $X$ is Fano) and the $E_1$-degeneration of Hodge to de Rham spectral sequence.
Unfortunately, 
there exist smooth Fano threefolds which are not quasi-$F$-split (see \cite[Example 7.2]{KTY} and \cite[Example 8.7]{KTLift2}).}

\subsection{Overview of proofs}\label{ss intro proofs}

In this paper, we prove Theorem \ref{Introthm:W(k)-lift}, 
\ref{Introthm:ANV}, \ref{Introthm:E_1-degeneration}, and \ref{Introthm:hoge number} 
for the case when $|-K_X|$ is very ample and $\Pic X$ 
is generated by $\MO_X(K_X)$. 
The remaining case will be treated in our next paper \cite{KTLift2}.

In what follows, we overview some of ideas of the proofs in this paper. 
Since $|-K_X|$ is very ample, we have the closed embedding $X \subset \P^{g+1}$ 
induced by $|-K_X|$, where $g$ is defined by $(-K_X)^3 = 2g -2$. 
Recall that $g$ is an integer satisfying $3 \leq g \leq 12$ and $g \neq 11$ \cite[Theorem 1.1]{FanoII}. 
As the case $g \leq 5$ is easy, we here treat the case when $g \geq 6$. 

As mentioned above, Theorem \ref{Introthm:E_1-degeneration} follows from 
Theorem \ref{Introthm:W(k)-lift} and Theorem \ref{Introthm:hoge number}. 
Hence it is enough to show Theorems \ref{Introthm:W(k)-lift}, \ref{Introthm:ANV}, and \ref{Introthm:hoge number}. 
Set $\Omega_X^i(nK_X) := \Omega_X^i \otimes \MO_X(nK_X)$. 
By $\Pic X = \Z K_X$, it suffices to compute 
\[
h^j(X, \Omega_X^i(nK_X)) 
\]
for suitable triples $(i, j, n)$, e.g., 
$H^1(X, \Omega_X^1(K_X))=0$ and $H^0(X, \Omega_X^i)=0$ for $i>0$. 
We treat the following two cases separately. 
\begin{enumerate}
\item[(I)] $g \geq 7$. 
\item[(II)] $g=6$. 
\end{enumerate}

(I) Assume $g \geq 7$, i.e., $X \subset \P^{g+1}$ is an anti-canonically embedded Fano threefold with $g \geq 7$ and $\Pic X = \Z K_X$. 
In this case, by using two-ray game, 
we prove that $X$ is rational or birational to a smooth cubic threefold. 
We can check $h^1(X, \Omega_X^1)=1$ and $H^0(X, \Omega^i_X)=0$ 
for $i>0$ (Theorem \ref{t-rat3-Pic-Dol}, Lemma \ref{l-sm-RET}). 

Then it is enough to show that $H^1(X, \Omega^1_X(nK_X))=0$ for $n>0$. 
The essential part is the case when $n=1$. 
By standard argument, the problem is reduced to the injectivity of $\beta : H^1(X, \Omega_X^1) \to H^1(S, \Omega_S^1)$. 
Consider the following commutative diagram consisting of natural $\F_p$-linear maps: 
\[
\begin{CD}
\Pic\,X \otimes_{\Z} \F_p @>\alpha >> \Pic\,S \otimes_{\Z} \F_p\\
@VV\dlog_X V @VV\dlog_S V\\
H^1(X, \Omega_X^1) @>\beta >> H^1(S, \Omega_S^1). 
\end{CD}
\]
By $\Pic X = \Z K_X$ and $k \simeq \Pic\,X \otimes_{\Z} k \xrightarrow{\simeq} H^1(X, \Omega_X^1)$ (Proposition \ref{p h^{1, 1}}), 
it is enough to show that both $\alpha$ and $\dlog_S$ are injective. 

\medskip

\underline{Injectivity of of $\dlog_S$}:
The injectivity of  $\dlog_S$ is ensured by (i) and (ii) below.  
\begin{enumerate}
\item[(i)] 
For a smooth K3 surface $S$, 
$S$ is not superspecial if and only if 
$\dlog_S : \Pic S \otimes_{\Z} \F_p \to H^1(S, \Omega_S^1)$ is injective. 
\item[(ii)] 
If $S$ is a general member, then $S$ is not superspecial. 
\end{enumerate}
The fact (i) is a result by van der Geer--Katsura \cite{vdGK00}. 
We now overview how to show (ii). 
Recall that there exists a unique superspecial smooth K3 surface $S$ over $k$. 
Suppose that every smooth member of $|-K_X|$ is isomorphic to $S$. 
It suffices to derive a contradiction. 
By two-ray game starting with a blowup at a point on $X$, 
we can check that the anticanonical embedding $X \subset \P^{g+1}$ is a Lefschetz embedding (Theorem \ref{t-Lef-emb}). 
Then, for two general members $S, S' \in |-K_X|$ and their intersection $B := S \cap S'$, 
the blowup $\sigma :Y \to X$ along $B$ admits a Lefschetz fibration $\pi : Y \to \P^1$, 
and hence every fibre of $\pi$ has at worst RDP (i.e., canonical) singularities. 
Then the pair $(Y, \sigma^*(-K_X))$  is a family of polarised RDP K3 surfaces over $\P^1$,  
which corresponds to a morphism $\theta : \P^1 \to  \mathcal M_{2g-2}$ 
to the moduli stack $\mathcal M_{2g-2}$ parametrising all the polarised RDP K3 surfaces $(T, A)$ satisfying $A^2 = 2g-2$. 
Since all but finite fibres of $\pi$ are isomorphic to $S$, 
we can check that the image $\theta(\P^1)$  to $\mathcal M_{2g-2}$ 
is a point. 
Since $\mathcal M_{2g-2}$ is a Deligne-Mumford stack and $\P^1$ has no non-trivial \'etale covers, $\pi : Y \to \P^1$ must be trivial, i.e., 
$Y \simeq S \times \P^1$ and $\pi$ coincides with the second projection.  
However, this leads to a contradiction: $2 = \rho(X)+1 =\rho(Y) \geq \rho(S) =22$. 

\medskip

\underline{Injectivity of of $\alpha$}: 
By $\Pic X \otimes_{\Z} \F_p = \F_p K_X$, 
the injectivity of 
the restriction map 
$\alpha : \Pic X \otimes_{\Z} \F_p \to \Pic S \otimes_{\Z} \F_p$ is equivalent to the condition that 
$\MO_X(-K_X)|_S$ is not divisible by $p$ in $\Pic S$. 
Suppose not, i.e., we have $\MO_X(-K_X)|_S \simeq \MO_S (pD)$ for some divisor $D$ on $S$. 
We then  get  
\[
2g -2 = (-K_X)^3 = ((-K_X)|_S)^2 = (pD)^2 =p^2D^2 \in 2p^2 \Z, 
\]
where $D^2 \in 2\Z$ follows from the Riemann-Roch theorem: 
$\chi(S, D) = \chi(S, \MO_S) + \frac{1}{2}D^2$. 
By $3 \leq g \leq 12$, we get two solutions $(p, g) = (2, 9)$ and $(p, g) = (3, 10)$. 
In order to eliminate these cases, 
we  apply a two-ray game starting with a blowup 
$\sigma : Y \to X$ along a conic $\Gamma$ on $X$. 
We have a flopping contraction $\psi :Y \to Z$. 
For its flop $\psi^+ : Y^+ \to Z$, we get another contraction $\tau: Y^+ \to W$. 
\begin{equation}\label{e flop1 Intro}
\begin{tikzcd}
Y \arrow[d, "\sigma"'] \arrow[rd, "\psi"]& & Y^+ \arrow[ld, "\psi^+"'] \arrow[d, "\tau"]\\
X & Z & W.
\end{tikzcd}
\end{equation}

We now overview the case $(p, g)=(3, 10)$. 
Recall that $|-K_Y|$ is base point free and $-K_Y$ is big. 
Set $\overline{Z} := \varphi_{|-K_Y|}(Y)$ and let $\theta : Z \to \overline{Z}$ be the induced morphism. 
In this case, a general member $T_Y$ of $|-K_Y|$ is smooth, 
because  $\theta : Z \to \overline{Z}$ is either an isomorphim or a double cover. 
A key part is to ensure that 
the induced morphism $T_Y \to T := \sigma_*T_Y$ is an isomorphism. 
We then get $T \simeq T_Y \simeq T_Z \simeq T_{Y^+}$ 
for the proper transforms $T_Z$ and $T_{Y^+}$ of $T_Y$ on $Z$ and $Y^+$, respectively. 
Moreover, our assumption $\MO_X(-K_X)|_S \simeq \MO_S(p D)$ implies the same conclusion for $T$, i.e., 
we can find a Cartier divisor $D_T$ on $T$ satisfying $\MO_X(-K_X)|_T \simeq \MO_T(p  D_T)$. 
Then we get a contradiction by computing suitable intersection numbers on $T_{Y^+}$. 

The other case  $(p, g)=(2, 9)$ is similar to but more complicated than the case $(p, g)=(3, 10)$. 
This is because $\theta : Z \to \overline{Z}$ can be inseparable. 
If $\deg \theta=1$, then we may apply almost the same argument as in the case $(p, g)=(3, 10)$. 
Otherwise, $\deg \theta =2$. 
In this case, we have $W=\P^1$, which implies that $\theta : Z \to \overline{Z}$ is inseparable. 
We then get the diagram corresponding to (\ref{e flop1 Intro}) which consists of contractions of projective normal varieties: 
\begin{equation}\label{e flop2 Intro}
\begin{tikzcd}
\ol{Y} \arrow[d, "\sigma'"'] \arrow[rd, "\psi'"]& & \ol{Y}^+ \arrow[ld, "\psi'^+"'] \arrow[d, "\tau'"]\\
\ol{X} & \ol{Z} & \ol{W}=\P^1. 
\end{tikzcd}
\end{equation}
In particular, 
the diagrams (\ref{e flop1 Intro})
and 
(\ref{e flop2 Intro}) are homeomorphic, 
and  $\ol{X}, \ol{Y}, \ol{Z}, \ol{Y}^+$ are birational. 
A key observation is the fact that $\overline Z$ is a toric variety. 
This is confirmed by proving that $\overline Z$ is a variety of minimal degree. 
Since being toric is stable under small birational maps and contraction images, 
all of the varieties in the diagram (\ref{e flop2 Intro}) are toric. 
By explicit toric computation, we see that $\overline{Y}^+$ is a smooth toric Fano threefold. 
However, this is a contradiction, because a smooth Fano threefold cannot admit 
a small birational morphism $\psi'^+ : \overline{Y}^+ \to \overline Z$ by classification of extremal rays \cite{Kol91}. 

\medskip

(II) Assume $g=6$, i.e., $X \subset \P^{7}$ is an anti-canonically embedded Fano threefold of genus $6$ with $\Pic X = \Z K_X$. 
The main part is to establish (a)-(d) below. 
\begin{enumerate}
\item[(a)] One of the following holds.  
\begin{enumerate}
\item[(a1)] 
$X$ is isomorphic to a complete intersection 
    $\Gr(2, 5) \cap H_1 \cap H_2  \cap Q$, 
    where each $H_i$ is a hyperplane of $\P^9$ and $Q$ is a quadric hypersurface of $\P^9$.
\item[(a2)] 
There exists a double cover $f: X \to V$, 
where $V$ is a smooth Fano threefold with $r_{V}=2$ and $(-K_{V}/2)^3 = 5$. 
\end{enumerate}
\item[(b)] $f: X \to V$ is separable for the case (a2). 
\item[(c)] $h^1(X, \Omega^1_X)=1$. 
\item[(d)] $H^3(X, \Omega^1_X)=0$. 
\end{enumerate}
Roughly speaking, the proof of (a) is carried out by applying the same argument as in the paper by Debarre-Kuznetsov \cite{DK18}, although we need to be careful with 
some minor parts. 
As for (b), it is enough to show that 
\[
c_3(\Omega_V^1 \otimes \MO_V(2D_V)) \neq 0
\]
for an ample Cartier divisor $D_V$ on $V$ satisfying $-K_V \sim 2D_V$. 
This can be checked by an explicit computation, because 
$V$ is a linear section of $\Gr(2, 5)$ and 
intersection theory on $\Gr(2, 5)$ is governed by Schubert calculus. 

The assertions  (c) and (d) hold for  the case (a1) by the fact that 
these cohomologies can be computed by the ones of the ambient space $\Gr(2, 5)$. 
Hence we may assume that (a2) holds. 
Let us overview the proof of (c). 
Since (the proof of) (a) implies that $X$ lifts to $W(k)$, 
we get $h^1(X, \Omega_X) \geq 1$ by upper semi-continuity of cohomologies. 
Hence it suffices to show $h^1(X, \Omega_X) \leq 1$.  
The exact sequence $0 \to \MO_V \to f_*\MO_X \to \MO_V(-D_V) \to 0$ 
induces the following one: 
\[
0 \to \Omega_V^1 \to f_*f^*\Omega_V^1 \to \Omega_V^1(-D_V) \to 0. 
\]
By the Akizuki-Nakano vanishing on $V$, 
we get $h^1(X, f^*\Omega_V^1) = h^1(V, f_*f^*\Omega_V^1) 
= h^1(V, \Omega_V^1)=1$. 
Then 
the problem is reduced to the vanishing $H^1(X, \Omega^1_{X/V})=0$ by 
\[
0 \to f^*\Omega_V^1 \to \Omega_X^1 \to \Omega^1_{X/V} \to 0. 
\]
Recall that $X$ can be embedded into 
the $\P^1$-bundle $P := \mathbb P_V(\MO_V \oplus \MO_V(-D_V))$ over $V$. 
Consider the exact sequence 
\[
0 \to I/I^2 \to \Omega_{P/V}^1 \to \Omega^1_{X/V} \to 0, 
\]
where $I$ denotes  the ideal sheaf on $P$ defining $X$. 
Since $I/I^2$ and $\Omega_{P/V}^1|_X$ are invertible sheaves on $X$, 
the equality $\Pic X = \Z K_X$ and Kodaira vanishing on $X$ imply $H^i(X, I/I^2) =H^i(X,\Omega_{P/V}^1|_X)=0$ for $i \in \{1, 2\}$. 
Therefore, we get $H^1(X, \Omega^1_{X/V})=0$.

For the proof of (d), the problem is reduced to $H^3(V, \Omega^1_V(-D_V))=0$ by applying a similar argument to the one of (c). 
By Serre duality, this is equivalent to $H^0(V, \Omega_V^2(D_V))=0$. 
In order to prove this, we apply 
two-ray game for $V$ starting with a blowup $W\to V$ along a conic on $V$. 
We then get a birational morphisms
\[
\begin{tikzcd}
 & W \arrow[rd] \arrow[ld] & \\
\P^3 & & V, 
\end{tikzcd}
\]
where $W \to \P^3$ is a blowup along a smooth rational curve $B$ on $\P^3$ of degree $4$. 
By comparing $\Omega^2_V$ and $\Omega^2_{\P^3}$ via $\Omega_W^2$, 
the problem is reduced to 
\begin{equation}\label{e P3 IB intro}
H^0(\P^3, \Omega_{\P^3}^2 \otimes \MO_{\P^3}(3) \otimes I_B)=0. 
\end{equation}
The  Euler sequence 
\[
0 \to \Omega^1_{\P^3} 
\to 
\bigoplus_{i=0}^3 
\MO_{\P^3}(-1)e_i \to \MO_{\P^3} \to 0 
\]
and its $\wedge^2$ induce an exact sequence 
\[
0 \to \Omega^2_{\P^3} \to \bigoplus_{0\leq i < j \leq 3} \MO_{\P^3}(-2) e_i \wedge e_j \to  \bigoplus_{0\leq i \leq 3}\MO_{\P^3}(-1) e_i. 
\]
We then get (\ref{e P3 IB intro}) by using this explicit description of $\Omega^2_{\P^3}$.

\subsection{Overview of contents}

In Section \ref{s prelim}, 
we summarise notation and some known results. 
Let $X \subset \P^{g+1}$ be an anti-canonically embedded Fano threefold. 
As explained in Subsection \ref{ss intro proofs}, 
our goal is to compute $h^j(X, \Omega_X^i(nK_X))$ for suitable triples $(i, j, n)$. 
The case $g \geq 7$  is treated in Section \ref{s Lef} and Section \ref{s lift}. 
In  Section \ref{s Lef}, we prove that $X \subset \P^{g+1}$ is a Lefschetz embedding. 
Section \ref{s lift} is mainly devoted to the injectivity of 
$\alpha : \Pic X\otimes_{\Z} \F_p \to \Pic S \otimes_{\Z} \F_p$ 
(cf. (I) in Subsection \ref{ss intro proofs}). 
In Section \ref{s genus 6}, we treat the case $g=6$ 
(cf. (II) in Subsection \ref{ss intro proofs}).  
In Section \ref{s main proofs}, we give proofs of our main theorems. 

Section \ref{s Veronese} and Section \ref{s-rat3-Dol} are appendices of this paper. 
In Section \ref{s Veronese}, we extend some classical results on Veronese surfaces to arbitrary characteristics. 
If ${\rm char}\,k \neq 2$, then 
the classical method 
in characteristic zero  works without any changes, 
whilst we need to modify some details in characteristic two. 
The main result in Section \ref{s Veronese} is used to show the existence of conics for the case $g=7$ (Subsection \ref{ss exist conic}). 
The main result of  Section \ref{s-rat3-Dol} is to prove 
that $\dlog : \Pic X \otimes_{\Z} k \to H^1(X, \Omega_X^1)$ is an isomorphism if $X$ is a smooth projective threefold which is rational or birational to a smooth cubic threefold (Theorem \ref{t-rat3-Pic-Dol}). 


\medskip
\noindent {\bf Acknowledgements.}
The authors thank 
Emiliano Ambrosi, Jefferson Baudin, Katsuhisa Furukawa, Keiji Oguiso, Teppei Takamatsu, and Burt Totaro for answering questions and constructive suggestions. 
Kawakami was supported by JSPS KAKENHI Grant number JP22KJ1771 and JP24K16897.
Tanaka was supported by JSPS KAKENHI Grant number JP22H01112 and JP23K03028. 

\section{Preliminaries}\label{s prelim}



\subsection{Notation}
\label{ss:notation}

In this subsection, we summarise notation and basic definitions used in this article. 
\begin{enumerate}
\item Throughout the paper, $p$ denotes a prime number and we set $\F_p \coloneqq \Z/p\Z$. 
Unless otherwise specified, we work over an algebraically closed field $k$ of characteristic $p>0$. 
We denote by $F \colon X \to X$ the absolute Frobenius morphism on an $\F_p$-scheme $X$.
    \item We say that $X$ is a {\em variety} (over a field $\kappa$) if 
    $X$ is an integral scheme 
    that is separated and of finite type over $\kappa$. 
    We say that $X$ is a {\em curve} (resp. {\em surface}, resp. {\em threefold}) 
    if $X$ is a variety of dimension one (resp. two, resp. three). 
\item For a variety $X$, 
we define the {\em function field} $K(X)$ of $X$ 
as the stalk $\MO_{X, \xi}$ at the generic point $\xi$ of $X$. 
\item 
Given a coherent sheaf $\cF$ and a Cartier divisor $D$ on a variety $X$, 
we set $\cF(D) \coloneqq \cF \otimes \MO_X(D)$ unless otherwise specified. 
\item 
Given a closed subscheme $X$ of $\P^n$, we set $\MO_X(a) := \MO_{\P^n}(a)|_X$ for $a \in \Z$ unless otherwise specified. 
\item Given two closed subschemes $Y$ and $Z$ on a scheme $X$, we denote by $Y \cap Z$ the scheme-theoretic intersection, i.e., 
$Y\cap Z \coloneqq Y \times_X Z$. 
\item 
Recall that $\dlog : \MO_X^{\times} \to \Omega_X^1, f \mapsto df/f$ is a group homomorphism, 
and hence so is 
\[
\dlog :  \Pic\,X =\check{H}^1(X, \MO_X^{\times})  \to \check{H}^1(X, \Omega_X^1) \simeq H^1(X, \Omega_X^1). 
\]
By abuse of notation, 
each of the induced maps 
$\Pic\,X  \otimes_{\Z} \F_p \to H^1(X, \Omega_X^1)$ and 
$\Pic\,X  \otimes_{\Z} k \to H^1(X, \Omega_X^1)$ 
is also denoted by $\dlog$. 
\item 
Let $X$ be  a smooth projective variety. 
We say that $X$ {\em satisfies Akizuki-Nakano vanishing} 
if $H^j(X, \Omega_X^i(-A))=0$ 
for an ample divisor $A$ and integers $i$ and $j$ satisfying $i+j <3$. 
The spectral sequence 
\[
E_1^{i, j} = H^j(X, \Omega_X^i) \Rightarrow H^{i+j}(X, \Omega_X^{\bullet})
\]
is called the {\em Hodge to de Rham spectral sequence} of $X$. 
\item\label{ss nota lift}
Given a projective variety $X$ over $k$, 
we say that $\mathcal X$ is a {\em $W(k)$-lift} of $X$ (or a {\em lift of $X$ over} $W(k)$)  
if $\mathcal X$ is a flat projective scheme over $W(k)$ 
satisfying $\mathcal X \times_{\Spec W(k)} \Spec k \simeq X$. 
We say that $X$ {\em lifts to} $W(k)$ if there exists a $W(k)$-lift of $X$.
\end{enumerate}

\subsubsection{Fano threefolds}

We say that $X$ is a {\em Fano threefold} if 
$X$ is a smooth projective threefold such that $-K_X$ is ample. 
We say that $X \subset \P^{g+1}$ is an {\em anti-canonically embedded Fano threefold} 
if $X$ is a Fano threefold, $X$ is a closed subscheme of $\P^{g+1}$, 
$\MO_X(-K_X) \simeq \MO_{\P^{g+1}}(1)|_X$, 
and the induced $k$-linear map 
\[H^0(\P^{g+1}, \MO_{\P^{g+1}}(1)) 
\to H^0(X, \MO_{\P^{g+1}}(1)|_X)\] is an isomorphism.

\subsection{Two-ray game}

\begin{prop}\label{p FanoY cont}
Let $X \subset \P^{g+1}$ be an anti-canonically embedded Fano threefold. 
For a point (resp. a conic, resp. a line) $\Gamma$ on $X$, 
let $\sigma : Y \to X$ be the blowup along $\Gamma$. 
Assume that $g \geq 6$ (resp. $g \geq 5$, resp. $g \geq 4$). 
Moreover, when $\Gamma$ is a point, we assume that 
$\Gamma$ is not contained in any line on $X$. 
Then 
\begin{enumerate}
\item $|-K_Y|$ is base point free and $-K_Y$ is big. 
\end{enumerate}
For the induced morphism $\ol{\psi} : Y \to \overline{Z}$ 
onto the image $\overline{Z}$ by $\varphi_{|-K_Y|}$, 
let 
\[
\ol{\psi} : Y \xrightarrow{\psi} Z \xrightarrow{\theta} \ol{Z}
\]
be the Stein factorisation of $\ol{\psi}$. 
Assume that $\dim \Ex(\psi)=1$. 
Then the following hold. 
\begin{enumerate}
\setcounter{enumi}{1}
\item $\deg \theta =1$ or $\deg \theta =2$. 
\item If $\deg \theta =1$, then $\theta$ is an isomorphism. 
\item If $\deg \theta = 2$, then $\Delta(\overline{Z}, A_{\overline Z})=0$, {where $\Delta(\ol{Z}, A_{\ol{Z}}) := \dim \ol{Z} +A_{\ol{Z}}^3 -h^0(\ol Z,  A_{\ol{Z}})$.}
\end{enumerate}
\end{prop}

\begin{proof}
The assertion (1) follows from the same argument as in \cite{Pro}. 
By construction, we can find  an ample Cartier divisor $A_{\overline Z}$ on $\ol{Z}$ 
satisfying $\ol{\psi}^*A_{\ol{Z}} \sim -K_Y$ and 
$\ol{\psi} : H^0(\ol{Z}, A_{\ol{Z}}) \xrightarrow{\simeq} H^0(Y, -K_Y)$ \cite[Lemma 3.13]{FanoI}. 
The assertions (2) and (4) hold by 
\[
0 \leq \Delta(\ol{Z}, A_{\ol{Z}}) = \dim \ol{Z} +A_{\ol{Z}}^3 -h^0(\ol Z,  A_{\ol{Z}})
\]
\[
=3 + \frac{(-K_Y)^3}{ \deg \ol{\psi}} -h^0(Y, -K_Y) 
\leq  (-K_Y)^3(\frac{1}{\deg \ol{\psi}} - \frac{1}{2}). 
\]
The assertion (3) follows from the same argument as in \cite[Theorem 6.2]{FanoI}. 
\end{proof}

\begin{lem}[cf.~\cite{IP99}*{Remark 4.1.10}]\label{l FanoY cont insep 2-to-1}
We use the same notation as in Proposition \ref{p FanoY cont}. 
Let $Y^+$ be the flop of the flopping contraction $\psi : Y \to Z$. 
Assume that $\deg \theta =2$ and 
$Y$ is not isomorphic to $Y^+$ as $k$-schemes. 
Then the following hold. 
\begin{enumerate}
\item $p=2$. 
\item $\theta : Z \to \overline{Z}$ is a finite inseparable morphism of degree two. 
\item $\overline Z$ is not $\Q$-factorial. 
\end{enumerate}
\end{lem}

\begin{proof}
We first finish the proof by assuming (2). 
Clearly, (2) implies (1), and hence 
it suffices to prove (3). 
If $\overline Z$ is $\Q$-factorial, then  so is $Z$ \cite[Lemma 2.5]{Tan18b}, which contradicts 
the fact that $\psi : Y \to Z$ is a flopping contraction. 

It suffices to show (2). 
Suppose that $\theta : Z \to \overline Z$ is separable, i.e., 
$K(Z)/K(\overline{Z})$ is a Galois extension. 
Then 
\[
\Z  \simeq \Pic Z \subset {\rm Cl}(Z) = {\rm Cl}(Y) 
=\Pic Y \simeq \Z^2. 
\]
Since the Galois group $G := \Gal(K(Z)/K(\overline{Z}))$ acts on ${\rm Cl}(Z)$ 
and the linear equivalence class of $-K_Z$ is $G$-invariant, we have that 
\[
\Pic Z \otimes_{\Z} \Q =  {\rm Cl}(Z)^G\otimes_{\Z} \Q. 
\]
Pick an ample Cartier $H_Y$ on $Y$. For its pushforward $H_Z := \psi_*H_Y$ on $Z$, 
we have that $m(H_Z +\iota^*H_Z)$ is a $G$-invariant Cartier divisor for some 
integer $m>0$, 
where $\iota: Z \to Z$ denotes the Galois involution. 
Therefore, we get the following contraction: 
\[
Y = \Proj\,(\bigoplus_{d \geq 0} \MO_Z(mH_Z)) \overset{\iota}{\simeq} 
\Proj\,(\bigoplus_{d \geq 0} \MO_Z(-m\iota^*H_Z)) = Y^+, 
\]
as required. 
\end{proof}

\subsection{K3 surfaces}

Throughout this subsection, 
we work over an algebraically closed field of characteristic $p>0$. 
We say that $S$ is a {\em smooth (resp.~RDP) K3 surface} 
if $S$ is a smooth (resp.~canonical) projective surface such that $K_S \sim 0$ and $H^1(S, \MO_S)=0$. 
Note that the miniamal resolution of an RDP K3 surface is a smooth K3 surface.

Given a smooth K3 surface $S$, it is known that the following are equivalent 
(cf.~\cite[Section 2]{Ito18}). 
\begin{enumerate}
    \item $\rho(S) =22$. 
    \item The Artin-Mazur height $h(S)$ of $S$ is equal to $\infty$.  
\end{enumerate}
We say that $S$ is {\em supersingular} if one (each) of (1) and (2) holds.









\begin{dfn}\label{d sspecial}
We say that a K3 surface $X$ is {\em superspecial} if 
$X$ is supersingular and the Artin invariant of $X$ is equal to $1$ 
(cf.~\cite[Definition 2.3 and Proposition 2.4]{Ito18}). 
\end{dfn}

\begin{rem}\label{r sspecial unique}
Over an algebraically closed field $k$ of characteristic $p>0$, 
there exists a unique superspecial K3 surface up to isomorphisms 
\cite[Theorem 1.1]{DK03} ($p=2$), \cite[Corollary 7.14]{Ogu79} ($p>2$). 
\end{rem}

\begin{prop}\label{p ss pic inje}
Let $S$ be a smooth K3 surface. Then the following are equivalent. 
\begin{enumerate}
\item $S$ is not superspecial. 
\item $\dlog : \Pic\,S \otimes_{\Z} \Z/p\Z \to H^1(S, \Omega_S^1)$ is injective. 
\end{enumerate}
\end{prop}


\begin{proof}
The assertion follows from 
\cite[the proof of Lemma 11.2]{vdGK00}. 
\qedhere

\end{proof}

\begin{prop}\label{p K3 ss locus}
Let $\pi : V \to W$ be  a smooth projective morphism 
of quasi-projective schemes over an algebraically closed field $k$ of characteristic $p>0$. 
Assume that {a fibre} $V_w := \pi^{-1}(w)$ is a smooth K3 surface for every closed point $w \in W$. 
Then the following hold. 
\begin{enumerate}
\item $\{ w \in W \,|\, V_{\ol{w}} \text{ is supersingular}\}$ is a closed subset of $W$. 
\item $\{ w \in W \,|\, V_{\ol{w}} \text{ is superspecial}\}$ is a closed subset of $W$. 
\end{enumerate}
Here $\overline{w}$ denotes the geometric point of $w$, i.e., 
for the algberaic closure $\overline{k(w)}$ of the residue field $k(w)$ of $w$, 
we set $V_{\ol{w}} := V \times_W \Spec\,\overline{k(w)}$. 
\end{prop}

\begin{proof}
The assertions (1) and (2) follow from \cite[Theorem 6.2]{vdGK00} and 
\cite[Theorem 11.10]{vdGK00}, respectively. 
\end{proof}

\begin{prop}\label{p no isotriv K3}
Let $\pi : Y \to \P^1$ be a flat projective morphism, 
where $Y$ is a smooth projective threefold. 
Assume that the following hold. 
\begin{enumerate}
\item Every fibre is an RDP K3 surface. 
\item There exist infinitely many closed points of $\P^1$ whose fibres are smooth superspecial K3 surfaces. 
\end{enumerate}
Then $Y \simeq S \times \P^1$ and $\pi$ coincides with the second projection, 
where $S$ denotes the superspecial smooth K3 surface. 
\end{prop}

\begin{proof}
Set $B := \P^1$. 
By Proposition \ref{p K3 ss locus}, every smooth fibre of $\pi$ is isomorphic to $S$.

We first treat the case when $k$ is uncountable. 
Fix an ample Cartier divisor $H$ on $Y$. 
By the invariance of intersection numbers in flat families, 
the number $2d:=(H|_{\pi^{-1}(b)})^2 \in 2\Z_{>0}$ is independent of $b \in B$. 
Let $\theta : B \to \mathcal M_{2d}$ be the induced morphism to the moduli stack that parametrises 
all the polarised RDP K3 surfaces $\{ (S, A)\}$ with $A^2 = 2d$, whose existence is guaranteed by \cite[Proposition 2.1 and its proof]{Mau14}. 
We then get $\dim \theta(B) = 0$, because there are only countably many polarisation on a fixed K3 surface $S$. 
Since $B=\P^1$ is irreducible, $\mathcal P := \theta(B)$ is a closed point, which we equip with the reduced closed substack: 
\[
B=\P^1 \to \mathcal P \to \mathcal M_{2d}. 
\]
Since $\mathcal M_{2d}$ is a Deligne-Mumford stack, 
so is $\mathcal P$. 
Then there is an \'etale surjective morphism $u:Q := \Spec\,k \to \mathcal P$. 
Set $B' := B \times_{\mathcal P} Q$, which is an algebraic space \cite[Proposition 8.1.10]{Ols16}. 
Since $B' \to B$ is an \'etale morphism, 
we get $\dim B' = \dim B = 1$, and hence $B'$ is a scheme. 
Note that the base change $\pi' : Y \times_B B' \to B'$ of $\pi : Y \to B=\P^1$ is trivial, 
i.e., $Y \times_B B' \simeq S \times B'$ and $\pi'$ coincides with the second projection. 
As $B=\P^1$ has no non-trivial \'etale covers, we get $B \times_{\mathcal P} Q \simeq B \amalg \cdots \amalg B$. 
Then also $\pi$ is trivial. 
This completes the proof for the case when $k$ is uncountable.

Let us go back to the general case. 
We reduce the problem to the case when $k$ is uncountable. 
Take a field extension $k \subset \wt{k}$, 
where $\wt{k}$ is an uncountable algebraically closed field. 
In order to apply the argument in the previous paragraph, 
it suffices to show that the stack-theoretic image $\mathcal P := \theta(B)$ is zero-dimensional. 
This follows from the fact that 
$\mathcal P \times_{\Spec k} \Spec \wt{k} \to \mathcal P$ is faithfully flat and 
$\mathcal P \times_{\Spec k} \Spec \wt{k}$ is a zero-dimenisonal Deligne-Mumford  stack. 
\qedhere



\end{proof}

\subsection{Complete intersections}

In this subsection, we discuss Akizuki-Nakano vanishing and Hodge numbers for complete intersections inside good ambient spaces.

\begin{lem}\label{lem:ANV for CI}
Take integers $d >0$ and $n>0$.
    Let $P$ be a locally complete intersection (lci, for short) projective variety of $\dim\,P=d+n$. 
    Let $H_1, ..., H_n$ be ample effective Cartier divisors on $P$ 
    such that $X := H_1 \cap \cdots \cap H_n$ is complete intersection (i.e., $\dim X =d$) and $X$ is a smooth variety. 
Set $Z_{k}\coloneqq \bigcap_{i=k+1}^{n} H_i$ for $0 \leq k \leq n-1$. 
Assume that the equality 
    \begin{equation}\label{e1 ANV for CI}
    H^j(P,\Omega^{i}_P(-H))=0
    \end{equation}
 holds    if $i<d$, $i+j<d+n$, and $H$ is an ample Cartier divisor on $P$. 
    Then the following hold. 
    \begin{enumerate}
        \item $H^j(Z_{k},\Omega^{i}_{Z_{k}}(-H)))=0$ 
        if $i < d$, $i+j<d+k$, and 
        $H$ is an ample Cartier divisor on $P$. 
        \item $H^{j}(Z_{k}, \Omega^{i}_{Z_{k}}) \simeq H^j(Z_{k+1}, \Omega^{i}_{Z_{k+1}})$ 
        if $i < d$ and $i+j < d+k$. 
    \end{enumerate}
    In particular, the following hold by setting $k =0$. 
    \begin{enumerate}
        \item[(1)'] $H^j(X,\Omega^{i}_{X}(-H))=0$ if $i+j<d$ and 
        $H$ is an ample Cartier divisor on $P$. 
        \item[(2)'] $H^{j}(X, \Omega^{i}_{X})\simeq H^j(P, \Omega^{i}_{P})$ 
        if $i+j<d$.
    \end{enumerate}
\end{lem}
\begin{proof}
Let $\mathcal A$ be the set of all the ample divisors on $P$. 
For  $Z_n := P$, the following hold. 
    \begin{enumerate}
        \item[(a)] $Z_n=P$ and $Z_0=X$.
        \item[(b)] $Z_{k}=Z_{k+1}\cap H_{k+1}$ for all $0\leq k\leq n-1$,
        \item[(c)] $\dim Z_{k}=d+k$ and  $\dim (\Sing Z_{k}) \leq k-1$, 
        where $\Sing Z_k$ denotes the singular locus of $Z_k$. 
        In particular, $\dim Z_{k} - \dim (\Sing Z_{k}) \geq d+1$. 
    \end{enumerate}

{
We now prove Claim below.     
\begin{claim*}
The following hold for all $k$ and $i<d$. 
\begin{enumerate}
\renewcommand{\labelenumi}{(\roman{enumi})}
\item $Z_k$ is an lci normal variety. 
\item 
The sequence 
\begin{equation}\label{e2 ANV for CI}
    0\to \Omega^{i}_{Z_{k+1}}(-H_{k+1}) \to \Omega^{i}_{Z_{k+1}} \to \Omega^{i}_{Z_{k+1}}|_{Z_{k}} \to 0
\end{equation}
is exact. 
\item Each of $\Omega_{Z_k}^i$ and $\Omega_{Z_{k+1}}^i|_{Z_k}$  
is $S_2$. 
\item The sequence 
\begin{equation}\label{e3 ANV for CI}
    0\to \Omega^{i-1}_{Z_{k}}(-H_{k+1}) \to \Omega^{i}_{Z_{k+1}}|_{Z_{k}} \to \Omega^{i}_{Z_{k}} \to 0 
\end{equation}
is exact. 
\end{enumerate}
\end{claim*}

\begin{proof}[Proof of Claim]
Let us show (i).  By (c) and $d \geq 1$, we get 
    $\dim Z_{k} - \dim (\Sing Z_{k})\geq d+1 \geq 2$. 
    As $Z_k$ is Cohen-Macaulay, 
    we see that $Z_k$ is normal by Serre's criterion. 
    In particular, $Z_n = P$ is a normal variety. 
    Recall that $Z_{n-1}$ is an ample effective Cartier divisor on $Z_n=P$. 
    Then     $Z_{n-1}$ is connected, and hence $Z_{n-1}$ is a normal variety. 
    Repeating this procedure, we see that each $Z_k$ is an lci normal  variety. 
    Thus (i) holds. 

We now prove  ($\star$) below. 
\begin{enumerate}
\item[($\star$)] 
Fix $i <d$. 
Then $\Omega_{Z_1}^i$ is $S_2$ and $\Omega_{Z_k}^i$ is $S_3$ if $k \neq 1$. 
\end{enumerate}
If $k=0$, then ($\star$) follows from the fact that $Z_0 = X$ is smooth. 
The other case (i.e., $k \geq 1$) follows from \cite[Lemma 3.11]{Sato-Takagi}, 
which is applicable by (c). 
This completes the proof of ($\star$). 
Then the assertion (ii) holds, because 
$\Omega^i_{Z_{k+1}}(-H_{k+1})$ is torsion-free by ($\star$). 

Let us show (iii). 
The case when $k >0$ follows from ($\star$) and the exact sequence (\ref{e2 ANV for CI}). 
If $k=0$, then $Z_{k+1} = Z_1$ is smooth around $Z_k = Z_0 = X$, 
which implies that $\Omega^i_{Z_{k+1}}|_{Z_k}$ is locally free and hence $S_2$. 
This completes the proof of (iii). 



Let us show (iv). 
Set $Z_{k, {\rm sm}} := Z_k \setminus \Sing Z_k$, which is the smooth locus of $Z_k$. 
Since $Z_k$ is an effective Cartier divisor on $Z_{k+1}$, 
$Z_{k+1}$ is smooth around $Z_{k, {\rm sm}}$. 
Then the restricted sequence (\ref{e3 ANV for CI})$|_{Z_{k, {\rm sm}}}$ is exact, because it is obtained by applying the $i$-th wedge product $\wedge^i$ to the conormal sequence 
$0 \to \MO_P(-H_{k+1})|_{Z_{k, {\rm sm}}} \to \Omega^1_{Z_{k+1}}|_{Z_{k, {\rm sm}}}  \to \Omega^1_{Z_{k, {\rm sm}}} \to 0$. 
Taking the pushforward $j_*$ by the open immersion $j :Z_{k, {\rm sm}} \hookrightarrow Z_k$, 
we obtain the sequence (\ref{e3 ANV for CI}) 
(recall that the sheaves 
in (\ref{e3 ANV for CI}) are $S_2$ by (iii)). 
As $j_*$ is left exact, the sequence (\ref{e3 ANV for CI}) is exact. 
Then it is enough to show that 
$\Omega^{i}_{Z_{k+1}}|_{Z_{k}} \to \Omega^{i}_{Z_{k}}$ is surjective, 
which follows from the fact that this homomorphism coincides with 
the one obtained by applying $\wedge^i$ to the surjection 
$\Omega^{1}_{Z_{k+1}}|_{Z_{k}} \to \Omega^{1}_{Z_{k}}$. 
Thus (iv) holds. This completes the proof of Claim. 
\end{proof}
}





    We prove (1) by descending induction on $k$. 
    The base case  $k=n$ of this induction follows from $Z_k = Z_n = P$ 
    and our assumption (\ref{e1 ANV for CI}). 
    Fix $ 0 \leq k < n$. 
    Assume that 
\begin{equation}\label{e4 ANV for CI}
    H^j(Z_{k+1},\Omega^{i}_{Z_{k+1}}(-H))=0\quad
    \text{ if\quad $i< d$, $i+j<d+(k+1)$, and $H \in \mathcal A$.}
\end{equation}
It is enough to prove $H^j(Z_{k},\Omega^{i}_{Z_{k}}(-H))=0$ 
when $i< d$, $i+ j<d+k$, and $H \in \mathcal A$. 
We prove this by induction on $i$. 
The base case $i=0$ 
follows from  
(\ref{e4 ANV for CI}) and 
the exact sequence
    \[
    0\to \sO_{Z_{k+1}}(-H-H_{k+1})\to \sO_{Z_{k+1}}(-H) \to \sO_{Z_{k}}(-H) \to 0. 
    \]
Fix $0<i < d$, $j<d+k-i$, and $H \in \mathcal A$. 
By Claim(iv), we have the following exact sequence: 
    \[
    0\to \Omega^{i-1}_{Z_{k}}(-H-H_{k+1}) \to \Omega^{i}_{Z_{k+1}}(-H)|_{Z_{k}} \to \Omega^{i}_{Z_{k}}(-H) \to 0.
    \]
    By the induction hypothesis 
    (i.e., $H^j(Z_{k},\Omega^{i-1}_{Z_{k}}(-H))=0$ 
    for all $j<d+k-(i-1)$ and $H \in \mathcal A$), we obtain
    \[
    H^{j}(Z_{k}, \Omega^{i}_{Z_{k+1}}(-H)|_{Z_{k}}) \xrightarrow{\simeq}
    H^{j}(Z_{k}, \Omega^{i}_{Z_{k}}(-H)). 
    \]
Taking the tensor product (\ref{e2 ANV for CI}) $\otimes\,\, \MO_P(-H)$, 
we get  an exact sequence
    \[
    0\to \Omega^{i}_{Z_{k+1}}(-H-H_{k+1}) \to \Omega^{i}_{Z_{k+1}}(-H) \to \Omega^{i}_{Z_{k+1}}(-H)|_{Z_{k}} \to 0. 
    \]
    This, together with (\ref{e4 ANV for CI}), 
    implies $H^{j}(Z_{k}, \Omega^{i}_{Z_{k+1}}(-H)|_{Z_{k}})=0$. 
    Thus (1) holds.

    \medskip

    Let us show (2). 
    Fix integers $i, j, k$ satifying $i<d$ and $i+j < d+k$. 
    By (1) and the exact sequence (\ref{e3 ANV for CI}),  
    we get $H^{j}(Z_{k}, \Omega^{i}_{Z_{k+1}}|_{Z_{k}}) \simeq H^{j}(Z_{k}, \Omega^{i}_{Z_{k}})$. 
    By (1) and the exact sequence  (\ref{e2 ANV for CI}),  
    we obtain $H^j(Z_{k}, \Omega^{i}_{Z_{k+1}}) \simeq H^j(Z_{k}, \Omega^{i}_{Z_{k+1}}|_{Z_{k}})$. 
    These isomorphisms imply $H^{j}(Z_{k}, \Omega^{i}_{Z_{k}}) \simeq H^j(Z_{k+1}, \Omega^{i}_{Z_{k+1}})$. 
    Thus (2) holds. 
\qedhere 


\end{proof}

\begin{thm}\label{thm:CI}
Take integers $d >0$ and $n>0$.
    Let $P$ be an lci projective variety of $\dim\,P=d+n$. 
    Let $H_1, ..., H_n$ be an ample effective Cartier divisors on $P$ 
    such that $X := H_1 \cap \cdots \cap H_n$ is complete intersection (i.e., $\dim X =d$) and $X$ is a smooth variety. 
Assume that the following hold. 
    \begin{enumerate}
        \item $H^j(P,\Omega^{i}_P(-H))=0$ if $i<d$, $i+j<d+n$,  and $H$ is an ample Cartier divisor on $P$.
        \item 
        { There exists a lift $(\mathcal P; \mathcal H_1, ..., \mathcal H_n)$ 
        of $(P; H_1, ..., H_n)$ over $W(k)$, i.e.,        
$\mathcal P$ is a flat projective scheme over $W(k)$ and 
$\mathcal H_1, ..., \mathcal H_n$ are effective Cartier divisors on $\mathcal P$ 
such that there exists a $k$-isomorphism 
$\theta : \mathcal P \times_{\Spec W(k)} \Spec k \xrightarrow{\simeq}  P$ 
satisfying $\theta(\mathcal H_1 \times_{\Spec W(k)} \Spec k) = H_1, ..., \theta(\mathcal H_n \times_{\Spec W(k)} \Spec k) = H_n$.}
        \item $h^{i,j}(P)=h^{i,j}(P_{\overline{K}})$ for all $i,j\geq 0$, 
        where $P_{\overline K}$ denotes the geometric generic fibre of $\mathcal P \to \Spec W(k)$. 
    \end{enumerate}
    {Set $\mathcal{X}\coloneqq \mathcal{H}_1 \cap \cdots \cap \mathcal H_n$, which is a $W(k)$-lift of $X$.
    Let $X_{\overline{K}}$ be the generic fibre over $W(k)$. }
    Then 
    $h^{i,j}(X)=h^{i,j}(X_{\overline{K}})$     for all $i, j \geq 0$. 
\end{thm}

\begin{proof}
By (1), Lemma \ref{lem:ANV for CI} is applicable. 
If $i+j<d$, then the assertion follows from (3) and 
Lemma \ref{lem:ANV for CI}(2)'. 
By Serre duality: $h^{i, j}(X) = h^{d-i, d-j}(X)$, the assertion holds for the case when $i+j > d$. 

Hence we may assume   $i+j=d$. 
Fix $i$. 
Since Euler characteristic preserves under flat families, 
we obtain $\chi(\Omega^i_X)=\chi(\Omega^i_{X_{\overline{K}}})$. 
This equality implies 
$h^{i,d-i}(X)=h^{i,d-i}(X_{\overline{K}})$, 
because  we have already shown 
the equality $h^{i, j}(X) = h^{i,j}(X_{\overline{K}})$ 
for the other cases $j \neq d-i$. 
\qedhere

\end{proof}

\subsection{Existence of conics ($g \geq 7$)}\label{ss exist conic}

The purpose of this subsection is to prove that $X$  has a conic 
if $X \subset \P^{g+1}$ is an anti-canonically embedded Fano threefold with 
$g \geq 7$ and 
$\Pic\,X = \Z K_X$ (Theorem \ref{t conic g=7}).

\begin{nota}\label{n genus 7}
Let $X \subset \P^{8}$ be an anti-canonically embedded Fano threefold 
of genus $7$ such that 
$\Pic\,X = \Z K_X$. 
Fix a closed point $P \in X$ which no line on $X$ passes through. 
Take the blowup $\sigma: Y \to X$  at $P$. 
Recall that $|-K_Y|$ is base point free and $-K_Y$ is big 
(Proposition \ref{p FanoY cont}). 
Let 
\[
\varphi_{|-K_Y|} : Y \xrightarrow{\overline{\psi}} \overline Z \hookrightarrow \P^{h^0(Y, -K_Y)-1}
\]
be the induced morphism, where 
$\varphi_{|-K_Y|}$ is the morphism $\varphi_{|-K_Y|} : Y \to \P^{h^0(Y, -K_Y)-1}$ induced by $|-K_Y|$ and $\overline Z := \varphi_{|-K_Y|}(Y)$. 
Let $\psi : Y \to Z$ be the Stein factorisation of $\overline{\psi}$. 
\[
\begin{tikzcd}
& Y \arrow[ld, "\sigma"']\arrow[rd, "\psi"] \arrow[rdd, bend right, "{\overline{\psi}}"']\\
X && Z \arrow[d, "\theta"]\\
 & & \overline{Z}
\end{tikzcd}
\]
\end{nota}

\begin{prop}\label{p genus 7 prelim}
We use Notation \ref{n genus 7}. 
Assume $\dim \Ex(\psi) = 1$. 
Then the following hold. 
\begin{enumerate}
\item $h^0(E, -K_Y|_E)=6$. 
\item $H^0(Y, -K_Y)  \to H^0(E, -K_Y|_E)$ is injective. 
\item $h^0(Y, -K_Y) \in \{ 5, 6\}$. 
\end{enumerate}
\end{prop}

\begin{proof}
The assertion (1) follows from \cite[Proposition 6.1(3)]{FanoII}. 
By \cite[Proposition 6.6]{FanoII}, (2) holds. 
Let us show (3). 
It follows from (2) that 
\[
h^0(Y, -K_Y) \leq h^0(E, -K_Y|_E) = h^0(\P^2, \MO_{\P^2}(2)) = 6. 
\]
Hence the assertion (3) holds by 
\[
h^0(Y, -K_Y) \geq g -2 = 5, 
\]
where the inequality is guaranteed by \cite[Proposition 6.1(3)]{FanoII}.     
\end{proof}

\begin{lem}\label{l easy conic g=7}
We use Notation \ref{n genus 7}. 
Assume that there exists no conic on $X$. 
Then the following hold. 
\begin{enumerate}
\item $\dim \Ex(\psi) = 1$. 
\item $h^0(Y, -K_Y) = 5$. 
\item $|-K_Z|$ is very ample and $\theta : Z \to \overline Z$ is an isomorphism. 
\end{enumerate}
\end{lem}

\begin{proof}
Let us show (1). 
If $\dim \Ex(\psi) \neq 1$, then  
we get $\dim \Ex(\psi) =2$ \cite[Corollary 6.5]{FanoII}. 
In this case,  there exists a conic on $X$ by \cite[Proposition 6.2]{FanoII}, which is a contradiction. 
Thus (1) holds. 

Let us show (2). 
Suppose $h^0(Y, -K_Y) \neq  5$. 
In this case, the restriction map 
\[
\rho : H^0(Y, -K_Y ) \to H^0(E, -K_Y|_E)
\]
is surjective (Proposition \ref{p genus 7 prelim}). 
Then $X$ has a conic  by applying the same argument as in 
\cite[Theorem 6.7]{FanoII}, 
because the assumption $g \geq 8$ in \cite[Theorem 6.7]{FanoII} is used only to assure the surjectivity of $\rho$. 
Thus (2) holds.

Let us show (3). 
By construction, we can find  an ample Cartier divisor $A_{\overline Z}$ on $\ol{Z}$ 
satisfying $\ol{\psi}^*A_{\ol{Z}} \sim -K_Y$ and 
$\ol{\psi} : H^0(\ol{Z}, A_{\ol{Z}}) \xrightarrow{\simeq} H^0(Y, -K_Y)$. 
It is enough to prove that $\deg \theta \neq 2$ (Proposition \ref{p FanoY cont}). 

Suppose $\deg \theta =2$. 
Let us derive a contradiction. 
We have $A_{\overline Z}^3 = (-K_Y)^3/2 = (2g-10)/2 = 2$, 
which implies that $\overline Z$ is a quadric hypersurface. 
Recall that $h^0(Y, -K_Y) = 5$, $h^0(E, -K_Y|_E)=6$, and $\rho : H^0(Y, -K_Y) \to H^0(E, -K_Y|_E)$ is injective. 
Then the induced morphism $\overline{\psi}|_E : E \to \overline{\psi}(E)$ 
coincides with a projection of the Veronese surface $S \subset \P^5$ from a closed point $R \in \P^5 \setminus S$. 
Since the quadric hypersurface $\overline{Z} \subset \P^4$ contains $\overline{\psi}(E)$, 
$\overline Z$ must be singular (Theorem \ref{t Veronese no quad}). 
{
Then we can find distinct prime divisors $D$ and $D'$ on $\ol{Z}$ satisfying $D \neq  \overline{\psi}(E)$, 
$D' \neq  \overline{\psi}(E)$, and
$\MO_{\P^4}(1)|_{\ol{Z}} \simeq \MO_{\overline Z}(D+D')$ (Lemma \ref{l sing quad split}). 
Since $\overline{\psi} : Y \to \overline Z$ is surjective, 
there exist prime divisors $D_Y$ and $D'_Y$ such that 
$\overline{\psi}(D_Y) = D$ and 
$\ol{\psi}(D'_Y) =D'$. 
As $D+D'$ is an effective Cartier divisor on $\ol{Z}$, 
we get $\overline{\psi}^*(D+D') = D_Y + D'_Y +F$ for some effective divisor $F$ on $Y$. 
We have 
\[
\MO_Y(-K_Y) \simeq \varphi^*_{|-K_Y|}\MO_{\P^4}(1) 
= \overline{\psi}^*(\MO_{\P^4}(1)|_{\overline Z}) 
\simeq \overline{\psi}^*\MO_{\ol{Z}}(D+D') 
\simeq \MO_Y(D_Y+D'_Y + F). 
\]
By $D \neq  \overline{\psi}(E)$ and
$D' \neq  \overline{\psi}(E)$, 
we get $D_Y \neq E$ and $D'_Y \neq E$. 
Then both $D_X := \sigma(D_Y)$ and  $D'_X := \sigma(D'_Y)$  are prime divisors on $X$. 
However, this is absurd, because we get 
\[
-K_X \sim \sigma_*(-K_Y) \sim \sigma_*(D_Y+D'_Y+F) =D_X+D'_X+\sigma_*(F), 
\]
which contradicts $\Pic X = \Z K_X$.} 
\end{proof}

{ 
\begin{lem}\label{l sing quad split}
Let $V \subset \P^4$ be a quadric hypersurface. 
Assume that $V$ is an integral scheme which is not smooth. 
Fix a prime divisor $E_V$ on $V$. 
Then there exist distinct prime divisors $D$ and $D'$ on $V$ such that $D \neq E_V$, $D' \neq E_V$, and $\MO_{\P^4}(1)|_V \simeq \MO_V(D+D')$. 
\end{lem}

\begin{proof}
For $\P^4 = \Proj\,k[x_0, x_1, x_2, x_3, x_4]$, 
we may assume that one of the following holds. 
\begin{itemize}
\item $V = \{ x_0x_1 + x^2_2=0\} \subset \P^4$. 
\item $V = \{ x_0x_1 + x_2x_3=0\} \subset \P^4$. 
\end{itemize}
Take a hyperplane $H:= \{ x_2=0\} \subset \P^4$. 
For $D := \{ x_0 = x_2 =0 \} \subset \P^4$ and $D' := \{ x_1 = x_2 = 0\}$, we get 
\[
V \cap H = \{ x_0x_1 = x_2 =0\} = D \cup D', 
\]
which induces 
\[
\MO_{\P^4}(1)|_V \simeq \MO_{\P^4}(H)|_V \simeq \MO_V(V \cap H)  = \MO_V(D+D'). 
\]
This completes the proof for the case when 
$E_V \neq \{ x_0 = x_2 =0 \}$ and $E_V \neq \{ x_1 = x_2 = 0\}$. 

In what follows, we treat the case when 
$E_V = \{ x_0 = x_2 =0 \}$ or $E_V = \{ x_1 = x_2 = 0\}$. 
By symmetry, we may assume that $E_V = \{ x_0 = x_2 =0\}$. 
In this case, we may apply the same argument as above by setting 
$H := \{ x_1 -x_2 = 0\}$ (note that we have $E_V \not\subset H$ and $V \cap H$ is not irreducible).  
\end{proof}}

\begin{thm}\label{t conic g=7}
Let $X \subset \P^{g+1}$ be an anti-canonically embedded Fano threefold 
with $g \geq 7$ and 
$\Pic\,X = \Z K_X$. 
Then there exists a conic on $X$. 
\end{thm}

\begin{proof}
If $g \geq 8$, then the assertion holds by \cite[Theorem 6.7]{FanoII}. 
In what follows, we assume $g=7$. 
By applying the same argument as in \cite[Theorem 6.7]{FanoII}, 
it is enough to show that $\psi|_E : E \to E_Z := \psi(E)$ is an isomorphism. 
Suppose that  $\psi|_E$ is not an isomorphism. 
In this case, $\dim {\rm Sing} E_Z = 1$ (Proposition \ref{p Veronese cases}).
However, this is absurd, because $E \to E_Z$ is an isomorphism outside finitely many points of $E_Z$ (Lemma \ref{l easy conic g=7}(1)(3)). 
\end{proof}

\section{Existence of non-superspecial members ($g \geq 7$)}\label{s Lef}

\begin{dfn}[cf. \cite{SGA7II}*{Page 213}]\label{d Lef pt}
Let $X \subset \P^N$ be a smooth projective threefold. 
We say that $P \in X$ is a {\em Lefschetz point} if 
there exists a hyperplane $H$ of $\P^N$ such that 
\begin{enumerate}
\item the scheme-theoretic intersection $X \cap H$ is singular at $P$ and 
\item for the defining element $f \in \widehat{\MO}_{X, P}$ of the pullback of $H$ on $\Spec \widehat{\MO}_{X, P}$, 
the degree two part of $f$  is nonzero and defines a smooth conic on $\P^2_k$. 
\end{enumerate}
\end{dfn}

\begin{dfn}\label{d Lef pencil}
Let $X \subset \P^N$ be a smooth projective threefold. 
We say that a linear system $\Lambda \subset |\MO_{\P^N}(1)|$ is a {\em Lefschetz pencil} (of $X \subset \P^N$) 
if 
\begin{enumerate}
\item $X \cap H$ is smooth for a general member $H$ of $\Lambda$, 
\item 
$\Lambda$ is a pencil (i.e., the corresponding $k$-vector subspace $V_{\Lambda} \subset H^0(\P^N, \MO_{\P^N}(1))$ is of dimension two), and 
\item 
for every member $H$ of $\Lambda$ such that $X \cap H$ is singular, 
$X \cap H$ 
has a unique singular point $P$ and  
satisfies 
the property (2) of Definition \ref{d Lef pt} 
at $P$ (in particular, $H$ is a projective normal surface which has at worst canonical singularities
by Remark \ref{r Lef canonical}). 
\end{enumerate}
Here a {\em member} of ${\Lambda}$ is a 
hyperplane $H$ of $\P^N$ whose corresponding element of $H^0(\P^N, \MO_{\P^N}(1))$ is contained in $V_{\Lambda}$. 
\end{dfn}

\begin{rem}\label{r Lef canonical}
Let $X \subset \P^N$ be a smooth projective threefold. 
For a Lefschetz pencil $\Lambda \subset |\MO_{\P^N}(1)|$ 
of $X \subset \P^N$, 
take a member $H \in \Lambda$ such that $S := X \cap H$ is singular. 
By Definition \ref{d Lef pencil}(3), 
$S$ has a unique singular point.

Let us prove that $P$ is at worst canonical singularity of $S$. 
Since the singular locus of $S$ is zero-dimensional, 
$S$ is normal. 
By definition, the singularity $P \in S$ is defined by 
\[
f = q +(\text{higher}) \in k[[x, y, z]], 
\]
where $q$ is a homogeneous polynomial of degree $2$ such that 
$\Proj\,k[x, y, z]/(q)$ is a smooth conic. 
After a suitable coordinate change, we may assume that $q = xy+z^2$. 
For $\widehat{S} := \Spec \widehat{\MO}_{S, P}$, 
we see that $\widehat{S} \cap \{ z =0\} \subset \Spec\,k[[x, y]]$ is a nodal curve, and hence semi-log canonical. 
For $H := \{ z=0\}$, $(\widehat{S}, H)$ is log canonical by inversion of adjunction \cite[Theorem 5.1]{Tan18m}. 
Since $H$ is a nonzero effective Cartier divisor, $(\widehat{S}, 0)$ is canonical, and hence so is $S$. 
\end{rem}

\begin{lem}\label{l Lef dim4 pt}
Assume that $k$ is of characteristic two. 
Let $X \subset \P^N_k$ be a smooth projective threefold. Take 
a closed point $P \in X$. 
Assume that $\dim ({\rm Im}\,\theta_P) \geq 4$ for the induced $k$-linear map 
\[
\theta_P : H^0(X, (\MO_{\P^N}(1)|_X)  \otimes \m_P^2) 
\to 
H^0(X, (\MO_{\P^N}(1)|_X)  \otimes (\m_P^2/\m_P^3)). 
\]
Then $P$ is a Lefschetz point of $X \subset \P^N$. 
\end{lem}

\begin{proof}
In what follows, we make the following identification:  
$H^0(X, (\MO_{\P^N}(1)|_X)  \otimes (\m_P^2/\m_P^3))= \m_P^2/\m_P^3$. 
Fix a regular system of parameters $\m_P \MO_{X, P} = (x, y, z)$. 
Note that 
\[
\m_P^2/\m_P^3 = kx^2 + ky^2 + kz^2 + kxy + kyz + kzx \supset 
{\rm Im}\,\theta_P =:V.  
\]
After possibly permuting $x, y, z$, 
we can find $f \in V$ such that  
\[
f =ax^2 + by^2 + cz^2 + \alpha yz + \beta zx + xy. 
\]
We have 
\[
xy + \alpha yz + \beta zx = (x + \alpha z) (y+ \beta z) -\alpha \beta z^2. 
\]
By applying the coordinate change 
$(x+\alpha, z, y+\beta z, z) \mapsto (x, y, z)$, 
 we may assume $\alpha  = \beta =0$, i.e., 
\[
f = ax^2 + by^2 + cz^2 +xy. 
\]
If $c \neq 0$, then we are done. 
Hence we may assume that $f = xy+ax^2+by^2 \in V$. 

(I) If $V = k x^2 \oplus ky^2 \oplus kz^2 \oplus kf$, 
then  we get $z^2 +xy \in V$, as required. 
In what follows, we assume $V \not\subset k x^2 \oplus ky^2 \oplus kz^2 \oplus kf$.

(II) Take $g \in V  \setminus (k x^2 \oplus ky^2 \oplus kz^2 \oplus kf)$. 
After permuting $x$ and $y$, we may assume that 
$g = yz + \cdots$. 
If $g$ also has the term $zx$, i.e., $g=yz + \beta zx +\cdots$, 
then $yz + \beta zx = z(y+\beta x)$, and apply $y +\beta x \mapsto y$ (note that 
the form $xy +ax^2 + by^2$ is stable under this coordinate change). 
Hence we may assume that $g = yz + a' x^2 + b'y^2 + c'z^2$. 
If $a' \neq 0$, then we are done. 
We may assume that $g = yz + b' y^2 + c'z^2 \in V$.

(III) Assume that $h:=zx + a''x^2 + b'' y^2 + c'' z^2 \in V$. 
Take an element $\varphi \in V \setminus (k f + kg + kh)$. 
Erasing the terms $xy, yz, zx$ by using $f, g, h$, 
we may assume that $\varphi = Ax^2 + By^2 + Cz^2$. 
If $A \neq 0$, then $f+ \varphi$ defines a smooth conic. 
By symmetry, we are done. 

(IV) We may assume that every element of $V$ does not have the term $zx$. 
Recall that $f = x y + ax^2 + by^2, g= yz + a'x^2 +b'y^2$. 
Pick $ h \in V \setminus (kf + kg)$. 
Then, after erasing $xy$ and $yz$ by using $f$ and $g$, 
we may assume that $h = Ax^2 + By^2 + Cz^2$. 
If $C \neq 0$, then $f+h$ defines a smooth conic. Similarly, 
if $A \neq 0$, then $g+h$ defines a smooth conic. 
Hence we may assume that $A= C=0$, i.e., $h = B y^2$. 
Take $\varphi \in V \setminus (kf + kg + kh)$. 
Then we may assume that  $\varphi =A' x^2 + B'y^2 + C' z^2$ and $B' =0$. 
Then $A' \neq 0$ or $C' \neq 0$. 
In any case, we can find an element defining a smooth conic. 
\end{proof}

\begin{lem}\label{l Lef pt-to-pencil}
Let $X \subset \P^N$ be a smooth projective threefold. 
Assume that there exists a Lefschetz point $P \in X$. 
Then there exists a Lefschetz pencil of $X \subset \P^N$. 
\end{lem}


\begin{proof}
The assertion follows from 
 \cite[Corollaire 3.2.1 in page 219 and Lemme 4.1.2 in page 235]{SGA7II}. 
\end{proof}


\begin{thm}\label{t-Lef-emb}
Let $X \subset \P^{g+1}$ be an anti-canonically embedded smooth Fano threefold with $\Pic\,X = \Z K_X$ and $g \geq 7$. 
Then there exists a Lefschetz pencil of $X \subset \P^{g+1}$. 
\end{thm}

\begin{proof}
If $p \neq 2$, then the assertion follows from 
\cite[Corollaire 3.5.0 in page 227]{SGA7II}. In what follows, we assume $p =2$. 
Take a closed point $P$ such that no line on $X$ 
passes through $P$, whose existence is guaranteed by 
\cite[Proposition 5.4]{FanoII}. 
It is enough to show that $P \in X$ is a Lefschetz point of $X \subset \P^{g+1}$ (Lemma \ref{l Lef pt-to-pencil}). 
Let $\sigma :Y \to X$ be the blowup at $P$. 
Since $|-K_Y|$ is  base point free and $-K_Y$ is big 
(Proposition \ref{p FanoY cont}), 
there is a birational morphism $\psi : Y \to Z$ 
such that $\psi_*\MO_Y = \MO_Z$ and $-K_Y \sim \psi^*A_Z$ 
for some ample Cartier divisor $A_Z$ on $Z$. 
It is enough to show that the image of the restriction map 
\[
\rho : H^0(Y, \MO_Y(-K_Y)) \to H^0(E, \MO_Y(-K_Y)|_E)
\]
is of dimension $\geq 4$ (Lemma \ref{l Lef dim4 pt}). 
By $h^0(Y, \MO_Y(-K_Y)) \geq g-2$ \cite[Proposition 6.1(3)]{FanoII}, 
the problem is reduced to proving 
\[
h^0(Y, \MO_Y(-K_Y-E)) \leq g-6. 
\]

We now treat the case when $\dim \Ex(\psi)=1$ or $\psi$ is an isomorphism. 
If $g \geq 8$, then 
$\rho$ is surjective \cite[Proposition 6.6]{FanoII}, 
and hence $\dim \Im \rho =6 \geq 4$. 
If 
$g=7$, then 
$\rho$ is injective \cite[Proposition 6.6]{FanoII}, 
and hence $h^0(Y, \MO_Y(-K_Y-E)) = 0 \leq g-6$. 
This completes the proof for the case when  $\dim \Ex(\psi)=1$ or $\psi$ is an isomorphism.

In what follows, we assume that $\dim \Ex(\psi) = 2$. 
Set $D := \Ex(\psi)$. 
In this case, we can write 
\[
D \sim -\alpha K_Y -\beta E
\]
where $(g, \alpha, \beta) \in \Z^3$ satisfies one of the following possibilities (a)-(c) \cite[Proposition 6.2]{FanoII}. 
\begin{enumerate}
\renewcommand{\labelenumi}{(\alph{enumi})}
\item $(g, \alpha, \beta)= (7, 1, 1)$. 
\item $(g, \alpha, \beta)= (7, 2, 2)$. 
\item $(g, \alpha, \beta)= (9, 1, 2)$. 
\end{enumerate}
Since $D$ is a $\psi$-exceptional prime divisor, the following property 
$(\star)$ holds. 
\begin{enumerate}
\item[($\star$)] 
If $mD \sim F$ for some integer $m>0$ and effective divisor $F$, 
then $mD = F$. 
\end{enumerate}


\underline{(a), (b)}  
Assume (a), i.e., $(g, \alpha, \beta)= (7, 1, 1)$. 
In this case, $h^0(Y, -K_Y-E) = h^0(Y, D) =1 = g-6$, as required. 
Similarly, the case (b) is settled by $h^0(Y, -K_Y -E) =0$.

\underline{(c)} Assume $(g, \alpha, \beta)= (9, 1, 2)$. 
Recall that $B := \psi(D)$ is a curve \cite[Proposition 4.2]{FanoII}.  
Let $\widetilde{D}$ be the normalisation of $D$ and take the Stein factorisation of the composite morphism $\wt{D} \to D \to B$: 
\[
\begin{tikzcd}
\wt{D} \arrow[d] \arrow[r] &D \arrow[d] \\
\wt{B} \arrow[r] & B. 
\end{tikzcd}
\]
It holds that 
\[
-K_Y -E \sim D+E. 
\]
Then we get $H^0(\wt{D}, \MO_Y(D+E)|_{\wt{D}}) =0$, 
because $-(D+E) \sim K_Y +E$ is $\psi$-ample \cite[Proposition 4.2]{FanoII}. 
In particular, $H^0(D, \MO_Y(D+E)|_D)=0$. 
This, together with  an exact sequence   
\[
0 \to \MO_Y(E) \to \MO_Y(D+E) \to \MO_Y(D+E)|_D \to 0, 
\]
implies $H^0(Y, \MO_Y(E)) \xrightarrow{\simeq} 
H^0(Y, \MO_Y(D+E))$. 
Hence we get 
\[
h^0(Y, -K_Y-E) = h^0(Y, D+E) = h^0(Y, E) = 1 \leq 9-6 = g-6, 
\]
as required. 
\qedhere


\end{proof}

\begin{thm}\label{t ex non-ss}
Let $X \subset \P^{g+1}$ be an anti-canonically embedded Fano threefold 
with $g \geq 7$ and $\Pic\,X = \Z K_X$. 
Then a general member of $|-K_X|$ is a smooth K3 surface which is not  superspecial. 
\end{thm}

\begin{proof}
Suppose that all the smooth members of $|-K_X|$ are superspecial K3 surfaces. 
Fix a superspecial K3 surface $S$, which is well known to be unique up to isomorphisms 
(Remark \ref{r sspecial unique}). 
By Theorem \ref{t-Lef-emb}, 
there exists a Lefschetz pencil $\Lambda \subset |-K_X|$.  
For general members $S, S' \in \Lambda$, 
let $\sigma: Y \to X$ be the blowup along the smooth curve $\Gamma := S \cap S'$, so that we get a morphism $\pi : Y \to \P^1$ 
such that each member of $\Lambda$ corresponds to a fibre of $\pi$. 
Note that general fibres of $\pi$ are smooth K3 surfaces. 
Moreover, every fibre of $\pi$ is an RDP K3 surface, i.e., 
a K3 surface which has at worst canonical singularities 
(Remark \ref{r Lef canonical}). 
Then Proposition \ref{p no isotriv K3} is applicable, which implies $Y \simeq S \times \P^1$. 
This leads to the following contradiction: 
\[
2 = \rho(X) +1 = \rho(Y) = \rho(S \times \P^1) \geq \rho(S) =22. 
\]
\end{proof}

\section{Liftability ($g \geq 7$)}\label{s lift}


\subsection{Lifting criterion}

\begin{prop}\label{p rat or birat to cubic}
Let $X \subset \P^{g+1}$ be an anti-canonically embedded 
smooth Fano threefolds 
with $\Pic\,X = \Z K_X$. 
Then the following hold. 
\begin{enumerate}
\item $X$ is rational if $g \in \{9, 10, 12\}$. 
\item $X$ is rational or birational to a smooth cubic threefold if $g 
\geq 7$. 
\end{enumerate}
\end{prop}


\begin{proof}
We may assume $g \geq 7$. 
Then there exists a conic on $X$ by Theorem \ref{t conic g=7}. 
Depending on $g$, the assertions hold as follows. 
\begin{itemize}
\item Use a blowup along a conic when $g \in \{7, 9\}$ 
\cite[Proposition 7.2, Proposition 7.3, Theorem 7.4]{FanoII}. 
\item Use a blowup at a general point when $g \in \{ 8, 10, 12\}$
\cite[Proposition 6.1, Proposition 6.2, Theorem 6.3]{FanoII}. 
\end{itemize}
Here we used the fact that del Pezzo fibrations $\pi : Y^+ \to \P^1$ are rational 
if $(-K_{Y^+})^2 \cdot \pi^*(b) \geq 5$ for a closed point $b \in \P^1$ \cite[Theorem 1.4]{BT24}. 
\end{proof}

\begin{prop}\label{p h^{1, 1}}
Let $X \subset \P^{g+1}$ be an anti-canonically embedded 
smooth Fano threefolds 
with $g \geq 7$ and $\Pic\,X = \Z K_X$. 
Then the induced $k$-linear map 
\[
\dlog: \Pic\,X \otimes_{\Z} k \to H^1(X, \Omega_X^1)
\]
is an isomorphism. 
\end{prop}

\begin{proof}
The assertion follows from Proposition \ref{p rat or birat to cubic} and 
Theorem \ref{t-rat3-Pic-Dol}. 
\end{proof}

The following proposition is a corrected version of  
\cite[Proposition 12.1]{SB97}. 

\begin{prop}\label{p-SB-lifting}
Let $X \subset \P^{g+1}$ be an anti-caninically embedded Fano threefold 
with $g \geq 7$ and  $\Pic\,X = \Z K_X$. 
Assume that $(\star)$ holds. 
\begin{enumerate}
\item[($\star$)] There exists a smooth member $S \in |-K_X|$ such that $S$ is not a superspecial K3 surface and $\MO_X(-K_X)|_S$ is not divisible by $p$, i.e., there exists no Cartier divisor $D$ on $S$ satisfying 
    $\MO_X(-K_X)|_S \simeq \MO_S(pD)$. 
\end{enumerate}
Then $H^1(X,\Omega^1_X(nK_X))=0$ for every $n>0$. 
In particular, $X$ lifts to $W(k)$. 
\end{prop}
\begin{proof}
The 
\lq in-particular\rq\,\,part follows from $H^2(X, \Theta_X)\simeq H^1(X,\Omega^1_X(K_X))=0$ and \cite[Theorem 8.5.19]{FAG}, {where $\Theta_X$ denotes the tangent bundle of $X$}. 
The construction of 
\[
\dlog: \Pic\,X \otimes_{\Z} \F_p\to H^1(X, \Omega_X^1). 
\]
is functorial, so that we get the following commutative diagram consisting of $\F_p$-linear maps of $\F_p$-vector spaces: 
\[
\begin{CD}
\Pic\,X \otimes_{\Z} \F_p @>\alpha >> \Pic\,S \otimes_{\Z} \F_p\\
@VV\dlog_X V @VV\dlog_S V\\
H^1(X, \Omega_X^1) @>\beta >> H^1(S, \Omega_S^1). 
\end{CD}
\]

We now show that $\beta$ is injective. 
By $k \simeq \Pic\,X \otimes_{\Z} k \xrightarrow{\simeq} H^1(X, \Omega_X^1)$ (Proposition \ref{p h^{1, 1}}), 
$\dlog_X(-K_X)$ is a $k$-linear basis of $H^1(X, \Omega_X^1)$. 
Since $\beta$ is a $k$-linear map, 
it suffices to prove that $\beta \circ \dlog_X(-K_X) \neq 0$. 
By ($\star$), $\MO_X(-K_X)|_S$ is nonzero in $\Pic\,S \otimes_{\Z} \F_p$, i.e., 
$\alpha(-K_X) \neq 0$. 
Since $S$ is not superspecial, 
$\dlog_S : \Pic\,S \otimes_{\Z} \F_p \to H^1(S, \Omega_S^1)$ is injective 
(Proposition \ref{p ss pic inje}), and hence 
$\dlog_S \circ \alpha(K_X) \neq 0$. 
This completes the proof of the injectivity of $\beta$. 

\medskip 
It is enough to prove the following:
\begin{enumerate}
    \item $H^1(X, \Omega^1_X(\log S)(-S)) =0$.
    \item If $H^1(X, \Omega^1_X(\log S)(-nS)) =0$ for some $n>0$, then $H^1(X, \Omega^1_X(-nS)) =0$.
    \item If $H^1(X, \Omega^1_X(-nS)) =0$ for some $n>0$, then 
    $H^1(X, \Omega^1_X(\log S)(-(n+1)S)) =0$.
\end{enumerate}
This is because  
every ample divisor on $X$ is linearly equivalent to $nS$ for some $n>0$.

Let us show (1). 
By the exact sequence 
\[
0 \to \Omega^1_X(\log S)(-S) \to \Omega_X^1 \to \Omega_S^1 \to 0, 
\]
we obtain another exact sequence: 
\[
H^0(S, \Omega_S^1) \to H^1(X, \Omega^1_X(\log S)(-S)) 
\to H^1(X, \Omega_X^1) \xrightarrow{\beta} H^1(S, \Omega_S^1).
\]
Since $H^0(S, \Omega_S^1)=0$ \cite[Theorem 7]{RS76} and $\beta$ is injective, we obtain
\[
H^1(X, \Omega^1_X(\log S)(-S))=0.
\]
Thus (1) holds.

Next, we show (2).
Assume that $H^1(X, \Omega^1_X(\log S)(-nS)) =0$ for some $n>0$. We prove $H^1(X, \Omega^1_X(-nS)) =0$.
By the residue exact sequence, we obtain a short exact sequence
\[
0 \to \Omega^1_X(-nS) \to \Omega_X^1(\log S)(-nS) \to \sO_S(-nS) \to 0, 
\]
and thus $H^1(X, \Omega^1_X(-nS))=0$. 
Thus (2) holds. 

{Finally, we show (3).
Fix $n >0$ and assume that $H^1(X, \Omega^1_X(-nS)) =0$. 
We prove $H^1(X, \Omega^1_X(\log S)(-(n+1)S)) =0$.
Since $S$ is a smooth K3 surface, it lifts to $W_2(k)$, 
which implies $H^0(S, \Omega^1_S(-nS))=0$ by  Akizuki-Nakano vanishing \cite[Coroilaire 2.8]{DI}.
By the restriction exact sequence, we obtain a short exact sequence
\[
0 \to \Omega^1_X(\log S)(-(n+1)S) \to \Omega_X^1(-nS) \to \Omega_S^1(-nS) \to 0. 
\]
Therefore, we get $H^1(X, \Omega^1_X(\log S)(-(n+1)S)) =0$. Thus (3) holds.}
\end{proof}
\begin{rem}
   {In the above situation, we have proved $H^1(X, \Omega^1_X(\log S)(-S))=0$, and hence $(X,S)$ lifts to $W(k)$ \cite[Proposition 8.6]{log-deformation(Kato)} (cf.~\cite[Theorem 2.3]{KN}). }
\end{rem}


\subsection{Case $g \neq 9$}

\begin{prop}\label{p-SB-lifting2}
Let $X \subset \P^{g+1}$ be an anti-caninically embedded Fano threefold with $g \geq 7$ and $\Pic\,X = \Z K_X$. 
Assume that $ (p, g) \neq (2, 9)$.
{Then the following hold: 
\begin{enumerate}
\item 
If $m$ is an integer and 
$S$ is a smooth member of $|-K_X|$ such that $m \geq 2$ and 
$\MO_X(-K_X)|_S$ is divisible by $m$ in $\Pic S$, 
then $(m, g) = (2, 9)$ or $(m, g)=(3, 10)$. 
\item 
If $p \neq 2$ and $g =10$, then 
there exists a smooth member 
$S \in |-K_X|$ such that 
$\MO_X(-K_X)|_S$ is primitive in $\Pic S$ $($i.e., 
there exists no integer $m \geq 2$ such that $\MO_X(-K_X)|_S$ is divisible by $m$ in $\Pic S$$)$. 
\end{enumerate}
}
\end{prop}


\begin{proof}
{Let us show (1). Take $m$ and $S$ as in the statement.}  
Since $S$ is a smooth K3 surface, we have $D^2 \in 2\Z$ by the Riemann--Roch theorem $\chi(S, D) = \chi(S, \MO_S) + \frac{1}{2} D^2$. Then
\[
2g -2 = (-K_X)^3 = (-K_X)^2 \cdot S = (-K_X|_S)^2 = (mD)^2 = 2d m^2 
\]
for some $d \in \Z_{>0}$. 
By $7 \leq g \leq 12$ and $g \neq 11$, 
we obtain the two solutions 
{$(m, g) = (2, 9)$ and $(m, g)=(3, 10)$}. 
{Thus (1) holds.}



{Let us show (2). 
Assume $p\neq 2$ and $g=10$. 
Suppose that 
$\MO_X(-K_X)|_S$ is not primitive for every smooth member $S \in |-K_X|$, which is equivalent, by (1), to the condition that $\MO_X(-K_X)|_S$ is divisible by $3$. Let us derive a contradiction.}  
Fix a conic $\Gamma$ on $X$ and 
let $\sigma : Y\to X$ be the blowup along $\Gamma$. 
By Proposition \ref{p FanoY cont}, 
$|-K_Y|$ is base point free and $-K_Y$ is big. 
Let 
$\overline{\psi} : Y \to \overline Z$ be the induced morphism 
onto the image $\overline Z$ of $\varphi_{|-K_Y|}$ and let 
$\psi : Y \to Z$ be its Stein factorisation: 
\[
\overline{\psi} : Y \xrightarrow{\psi} Z \xrightarrow{\theta} \overline Z. 
\]
It follows from $g=10$ and \cite[Proposition 7.3 and Corollary 7.9]{FanoII}
 that 
$\dim \Ex(\psi) =1$, i.e., 
$\psi : Y \to Z$ is a flopping contraction.  
Let $\psi^+ : Y^+ \to Z$ be the flop of $\psi$. 
We have the extremal ray of $Y^+$ not corresponding to $\psi^+$ and 
let $\tau : Y^+ \to W$ be its contraction. 
By \cite[Theorem 7.4]{FanoII}, 
$\tau$ is of type $C_1$. 
We use the same notation as in \cite[Theorem 7.4]{FanoII}. 
In particular, $D := \tau^*L$ for a line $L$ on $W =\P^2$. 
We have $-K_{Y^+} \sim D+E^+$. 
Let $T_Y$ be a general member of $|-K_Y|$. Set $T := \sigma_*T_Y$.

\setcounter{step}{0}
\begin{step}\label{s1-SB-lifting}
$T_Y$ is a smooth K3 surface on $Y$. 
\end{step}

\begin{proof}[Proof of Step \ref{s1-SB-lifting}]
Note that $Z = \overline{Z}$ or $\theta: Z \to \overline Z$ is of degree two 
(Proposition \ref{p FanoY cont}). 
For the former case, $|-K_Z|$ is very ample, and hence there is nothing to show. 
Assume that $\theta: Z \to \overline Z$ is of degree two. 
In this case, $\overline Z$ is a variety of minimal degree, and hence $\overline Z$ is normal. 
Then $\theta: Z \to \overline Z$ is residually separable, i.e., 
for every point $z \in Z$ and its image $\theta(z)$, 
the induced field extension $\kappa(\theta(z)) \subset \kappa(z)$ is separable. 
Therefore, a Bertini theorem \cite[Corollary 4.3]{Spr98} can be applied, so that a general member $T_Y$ of $|-K_Z|$ is smooth. 
This completes the proof of Step \ref{s1-SB-lifting}. 
\end{proof}

\begin{step}\label{s2-SB-lifting}
$\sigma|_{T_Y} : T_Y \to T$ is an isomorphism. 
\end{step}

\begin{proof}[Proof of Step \ref{s2-SB-lifting}]
We first prove that $T$ is normal. 
We have the following commutative diagram in which each horizontal sequence is exact: 
\[
\begin{tikzcd}
0 \arrow[r] & \sigma_*\MO_Y(-T_Y) 
\arrow[r] & \sigma_*\MO_Y 
\arrow[r] & \sigma_*\MO_{T_Y} \arrow[r] 
& R^1\sigma_*\MO_Y(-T_Y) \\
0 \arrow[r] & \MO_Y(-T) 
\arrow[r]\arrow[u] & \MO_X 
\arrow[r]\arrow[u, equal] & \MO_{T} \arrow[r]\arrow[u] & 0.
\end{tikzcd}
\]
In order to show the normality of $T$, it suffices to prove 
$R^1\sigma_*\MO_Y(-T_Y)=0$, which follows from 
\[
R^1\sigma_*\MO_Y(-T_Y)= 
R^1\sigma_*\MO_Y(K_Y)= 0, 
\]
where the latter equality can be checked easily (e.g., \cite[Theorem 0.5]{Tan15}). 

Let us show that $\sigma|_{T_Y} : T_Y \to T$ is an isomorphism. 
To this end, it suffices to prove that $T_Y \cap E$ contains no fibre of $\sigma|_E : E \to \Gamma$. 
It holds that 
\begin{itemize}
\item $h^0(Y, -K_Y) \geq 9$ \cite[Proposition 7.2]{FanoII},
\item $h^0(E, -K_Y|_E) = 6$ \cite[Proposition 7.2]{FanoII}, and 
\item $h^0(Y, -K_Y -E ) =h^0(Y^+, -K_{Y^+}-E^+) =h^0(\P^2, \MO_{\P^2}(1)) = 3$ \cite[Theorem 7.4]{FanoII}.
\end{itemize}
Then the restriction map 
\[
H^0(Y, -K_Y) \to H^0(E, -K_Y|_E)
\]
is surjective.  
Thus we may assume that also $T_Y|_E$ is a general member of $|-K_Y|_E|$. 
Since $|-K_Y|_E|$ is base point free and $-K_Y|_E$ is big \cite[Proposition 7.2(1)(b)]{FanoII}, 
a general member $T_Y|_E$ is a prime divisor \cite[Proposition 2.11]{FanoI}. 
As $T_Y|_E$ is big, $T_Y|_E$ cannot be a fibre of $\sigma|_E : E \to \Gamma$, and hence $T_Y|_E$ does not contain any fibre of $\sigma|_E$.
This completes the proof of Step \ref{s2-SB-lifting}. 
\end{proof}

Set $T_Z := \psi_*T_Y$ and $T_{Y^+} := (\psi^+)_*^{-1}T_Z$. 
We obtain 
\[
T \xleftarrow{\simeq, \sigma} T_Y \xrightarrow{\simeq, \psi} T_Z \xleftarrow{\simeq, \psi^+} T_{Y^+}. 
\]


\begin{step}\label{s3-SB-lifting}
$(-K_{Y^+}+E^+)|_{T_{Y^+}} \sim 3A_{T_{Y^+}}$ for some Cartier divisor $A_{T_{Y^+}}$ on $T_{Y^+}$. 
\end{step}

\begin{proof}[Proof of Step \ref{s3-SB-lifting}]
{By Step \ref{s2-SB-lifting}, 
we have a linear equivalence $-K_X|_T \sim 3A$ 
for some Cartier divisor $A$ on $T$, which} implies 
\[
3A_{T_Y} := 3(\sigma|_{T_Y})^*A \sim (\sigma^*(-K_X))|_{T_Y} \sim (-K_Y+E)|_{T_Y}. 
\]
By $T_Y \simeq T_Z \simeq T_{Y^+}$, we get  
\[
3A_{T_{Y^+}} := 3 (\psi^+)_*^{-1}\psi_*A_{T_Y} \sim 
(-K_{Y^+}+E^+)|_{T_{Y^+}}. 
\]
This completes the proof of Step \ref{s3-SB-lifting}. 
\end{proof}

By taking the intersection with $(-K_{Y^+})|_{T_{Y^+}}$, we obtain 
\[
3\Z \ni 3A_{T_{Y^+}} \cdot ((-K_{Y^+})|_{T_{Y^+}})
= 
(-K_{Y^+}+E^+)|_{T_{Y^+}}  \cdot ((-K_{Y^+})|_{T_{Y^+}})
\]
\[
= (-K_{Y^+}+E^+)  \cdot (-K_{Y^+}) \cdot T_{Y^+}
=(-K_{Y^+}+E^+)  \cdot (-K_{Y^+})^2 \overset{(\star)}{=} (2g-8) +4=16, 
\]
where $(\star)$ holds by \cite[Lemma 4.5 and Proposition 7.2(1)]{FanoII}. 
This is absurd. 
\end{proof}



\subsection{Case $g=9$}

\begin{nota}\label{n-p=2-g=9}
Let $X \subset \P^{10}$ be an anti-canonically embedded smooth Fano threefold of genus $g=9$ such that $\Pic\,X = \Z K_X$. 
Fix a conic $\Gamma$ on $X$ and 
let $\sigma: Y \to X$ be the blowup along $\Gamma$. 
For the induced morphism $\varphi_{|-K_Y|}: Y \to \P^7$ 
(Proposition \ref{p FanoY cont})
and its image $\ol{Z} := \varphi_{|-K_Y|}(Y)$, 
take the Stein factorisation of the induced morphism $\ol{\psi} : Y \to \overline{Z}$: 
\[
\ol{\psi} : Y \xrightarrow{\psi} Z \xrightarrow{\theta} \overline{Z}. 
\]
Set $A_{\overline Z} := \MO_{\P^7}(1)|_{\overline Z}$. 
Note that $\dim \Ex(\psi)=1$ \cite[Proposition 7.3, Corollary 7.9]{FanoII}.  
Let $\psi^+ : Y^+ \to Z$ be the flop of $\psi : Y \to Z$ 
and let $\tau : Y^+ \to W :=\P^1$ be the contraction of the  $K_{Y^+}$-negative extremal ray \cite[Theorem 7.4]{FanoII}. 

\end{nota}

\begin{lem}\label{l p=2 g=9 h^0(-K_Y)}
We use Notation \ref{n-p=2-g=9}. 
Then $h^0(Y, -K_Y) = 8$ and the restriction map 
\[
H^0(Y, -K_Y)  \to H^0(E, -K_Y|_E)
\]
is surjective.  
\end{lem}

\begin{proof}
It holds that 
\begin{itemize}
\item $h^0(Y, -K_Y) \geq 8$ \cite[Proposition 7.2]{FanoII},
\item $h^0(E, -K_Y|_E) = 6$ \cite[Proposition 7.2]{FanoII}, and 
\item $h^0(Y, -K_Y -E ) =h^0(Y^+, -K_{Y^+}-E^+) =h^0(\P^1, \MO_{\P^1}(1)) = 2$ \cite[Theorem 7.4]{FanoII}.
\end{itemize}
Then the sequence 
\[
0 \to H^0(Y, -K_Y-E) \to H^0(Y, -K_Y) \to H^0(E, -K_Y|_E) \to 0
\]
must be exact, and hence $h^0(Y, -K_Y)=8$. 
\end{proof}

\begin{prop}\label{p p=2 g=9 non-birat}
We use Notation \ref{n-p=2-g=9}. 
Assume $\deg \theta =1$. 
Then 
{there exists a smooth member $S \in |-K_X|$ 
such that}  $\MO_X(-K_X)|_S$ is not divisible by $2$ in $\Pic S$.  


\end{prop}


\begin{proof}
Assume $\deg \theta  =1$, i.e., $\theta$ is an isomorphism (Proposition \ref{p FanoY cont}). 
Suppose that $\MO_X(-K_X)|_S$ is divisible by $2$ in $\Pic S$ for {every} smooth member $S \in |-K_X|$.

By $Z= \overline Z$, 
a general member $T_Y$ of $|-K_Y|$ is a smooth K3 surface on $Y$. 
Set $T := \sigma_*T_Y$, $T_Z := \psi_*T_Y$, and $T_{Y^+} := (\psi^+)_*^{-1}T_Z$. 
By the same argument as in Proposition \ref{p-SB-lifting2},  we see that 
\begin{enumerate}
\item[(I)] $T \xleftarrow{\simeq, \sigma} T_Y \xrightarrow{\simeq, \psi} T_Z \xleftarrow{\simeq, \psi^+} T_{Y^+}$ (cf. Step 2 of Proposition \ref{p-SB-lifting2}), and 
\item[(II)] $(-K_{Y^+}+E^+)|_{T_{Y^+}} \sim 2 A_{T_{Y^+}}$ for some Cartier divisor 
$A_{T_{Y^+}}$ on $T_{Y^+}$ (cf. Step 3 of Proposition \ref{p-SB-lifting2}).
\end{enumerate}
We now finish the proof by assuming $(\star)$. 
\begin{enumerate}
\item[($\star$)] $(\tau|_{T_{Y^+}})_*\MO_{T_{Y^+}} = \MO_W$ 
for the induced morphism 
$\tau|_{T_{Y^+}} : T_{Y^+} \to W = \P^1$. 
\end{enumerate}
It follows from (I) that $T_{Y^+}$ is a smooth K3 surface. 
By (II) and $-K_{Y^+} -E^+ \sim \tau^*\MO_{\P^1}(1)$, 
we get  
\[
(\tau|_{T_{Y^+}})^*\MO_{\P^1}(1) \sim 
\tau^*\MO_{\P^1}(1)|_{T_{Y^+}} 
\sim (-K_{Y^+} -E^+)|_{T_{Y^+}} 
\]
\[
= (-K_{Y^+} +E^+)|_{T_{Y^+}} -2E^+|_{T_{Y^+}} \sim 2D
\]
for some Cartier divisor $D$ on $T_{Y^+}$. 
This contradicts the fact that a genus-one fibration 
$\tau|_{T_{Y^+}} : T_{Y^+} \to W = \P^1$ of a smooth K3 surface 
$T_{Y^+}$ has no multiple fibre (cf. \cite[the proof of Proposition 1.6(iii) in Chapter 11]{Huy16}). 

It is enough to show $(\star)$. 
We have an exact sequence 
\[
0 = \tau_*\MO_{Y^+}(-T_{Y^+}) 
\to \tau_*\MO_{Y^+} \to \tau_*\MO_{T_{Y^+}}
\to R^1\tau_*\MO_{Y^+}(-T_{Y^+}).  
\]
Since the geometric generic fibre $F$ of $\tau : Y^+ \to W$ is
a canonical del Pezzo surface, 
we have that 
\[
R^1\tau_*\MO_{Y^+}(-T_{Y^+})|_F \simeq 
R^1\tau_*\MO_{Y^+}(K_{Y^+})|_F \simeq 
H^1(F, K_F)=0. 
\]
Hence 
$\MO_W \to (\tau|_{T_{Y^+}})_*\MO_{T_{Y^+}}$ is isomorphic at the generic point. 
This implies that $\wt{W} \to W$ is birational for  the Stein factorisation
\[
\tau|_{T_{Y^+}} : T_{Y^+} \to \wt{W} \to W. 
\]
Since $W = \P^1$ is normal, we get $\wt{W} = W$. 
Therefore, $(\star)$ holds. 
\qedhere

\end{proof}




We now recall some terminologies and results on rational normal scrolls. 
For details, we refer to \cite{EH87}. 

\begin{dfn}\label{d rat scroll}
For non-negative integers 
\[
0 \leq a_0 \leq a_1 \leq \cdots \leq a_d, 
\]
it is well known that $|\MO_P(1)|$ is base point free 
and $h^0(P, \MO_P(1)) = d+1+\sum_{i=0}^d a_i$
for 
\[
P := \P_{\P^1}(\MO_{\P^1}(a_0) \oplus \MO_{\P^1}(a_1) \oplus\cdots \oplus 
\MO_{\P^1}(a_d)). 
\]
Set $S_{\P^1}(a_0, a_1, ..., a_d) \subset \P^{d+\sum_{i=0}^d a_i}_k$ to be the image of $\varphi_{|\MO_P(1)|}$, 
which is called a {\em rational normal scroll}. 
\end{dfn}

\begin{rem}\label{r rat scroll}
\begin{enumerate}
\item $|\MO_P(1)|$ is very ample if and only if $a_0>0$ (i.e., all $a_0, ..., a_d$ are positive). 
\item $S_{\P^1}(0, a_1, ..., a_d)\subset \P^{d+\sum_{i=1}^d a_i}_k$ 
is the cone over 
$S_{\P^1}(a_1, ..., a_d)\subset \P^{d-1+\sum_{i=1}^d a_i}_k$ \cite[Page 6]{EH87}. 
\item $\deg S_{\P^1}(a_0, a_1, ..., a_d) = a_0+a_1+ \cdots + a_d$ \cite[Page 7]{EH87}. 
\end{enumerate}
\end{rem}

\begin{lem}
We use Notation \ref{n-p=2-g=9}. 
Assume that $\deg \theta =2$. 
Then the following hold. 
\begin{enumerate}
\item 
{$p=2$ and} $\theta : Z \to \overline{Z}$ is a finite inseparable morphism 
of degree two. 
\item 
$\Delta(\overline{Z}, A_{\overline Z}) =0$. 
Moreover, $\overline{Z}$ is isomorphic to $S_{\P^1}(0, 1, 4)$ or $S_{\P^1}(0, 2, 3)$. 
In particular, $\overline{Z}$ is toric. 
\end{enumerate}
\end{lem}

\begin{proof}
Since $Y$ is not isomorphic to  $Y^+$ as a $k$-scheme \cite[Theorem 7.4]{FanoII}, 
(1) holds and $\overline Z$ is not $\Q$-factorial (Lemma \ref{l FanoY cont insep 2-to-1}). 
Moreover, we have $\Delta(\overline Z, A_{\overline Z})=0$ (Proposition \ref{p FanoY cont}). 

It holds that  
$2 A^3_{\overline Z}  = (-K_Y)^3 = 2g-8 = 10$. 
By $1 = \rho(Z) \geq \rho(\overline Z)$, we see that $\rho(\overline Z)=1$. 
By classification of varieties of minimal degree 
(see \cite[Theorem 1]{EH87} or \cite[Theorem 5.10, Theorem 5.15]{Fuj90}), one of the following holds. 
\begin{enumerate}
\renewcommand{\labelenumi}{(\roman{enumi})}
\item $\overline Z = \P^3$ and $A^3_{\overline Z}=1$. 
\item $\overline Z \subset \P^4$ is a quadric hypersurface and $A^3_{\overline Z} =2$. 
\item $\overline Z$ is the cone over the  Veronese surface and 
$A^3_{\overline Z}=4$. 
\item $\overline Z$ is a cone over a rational normal curve. 
\item $\overline Z$ is a cone over a Hirzebrugh surface. 
\end{enumerate}
In any case except for (v), $\overline Z$ is $\Q$-factorial 
(cf. \cite[Proposition 8.5]{FanoI}).  
Hence (v) holds. 
We then get 
$\overline{Z} = S_{\P^1}(0, a, b) \subset \P^7$, 
which is the cone of 
$S_{\P^1}(a, b) \subset \P^6$ for $b \geq a >0$ (Remark \ref{r rat scroll}). 
It follows from Remark \ref{r rat scroll} that $a+b= A^3_{\overline Z} = 5$. 
Therefore, we have two solutions $(a, b) \in \{(1, 4), (2, 3)\}$. 
\qedhere

\end{proof}

\begin{lem}\label{l toric Q-facn}
Let $X$ and $Y$ be projective normal  varieties. 
Assume that \begin{enumerate}
\item $X$ is toric. 
\item $X$ and $Y$ are isomorphic in codimension one, 
i.e., there exist open subsets $X' \subset X$ and $Y' \subset Y$ 
such that $X' \simeq Y'$, ${\rm codim}_X(X \setminus X') \geq 2$, 
and ${\rm codim}_Y(Y \setminus Y') \geq 2$. 
\end{enumerate}
Then $Y$ is toric. 
\end{lem}

In what follows, we use minimal model program for toric varieties \cite{CLS11}. 
\begin{proof}
Taking a small $\Q$-factorialisation of $X$, we may assume that $X$ is $\Q$-factorial, i.e., the corresponding fan is simplicial. 
We have 
\[
{\rm Cl} (X) \simeq {\rm Cl} (Y). 
\]
Fix an ample Cartier divisor $A_Y$ on $Y$. 
Set $A_X$ to be 
 the Weil divisor on $X$ corresponding to $A_Y$. 
By 
\[
 H^0(X, mA_X) = H^0(Y, mA_Y) 
\]
$A_{X}$
is a big divisor. 
We run an $A_{X}$-MMP: 
\begin{equation}\label{e1 toric Q-facn}
X =:X_0 \dashrightarrow X_1 \dashrightarrow \cdots \dashrightarrow X_{\ell} =:X'.  
\end{equation}
Then each $X_i$ is toric and birational to $X =X_0$. In particular, $X'$ is a projective normal $\Q$-factorial toric variety. 
The proper transform $A_{X'}$ of $A_{X}$ on $X'$ is nef and big, 
because (\ref{e1 toric Q-facn}) is an $A_{X}$-MMP and $A_{X}$ is big. 
As $X'$ is toric, $A_{X'}$ is semi-ample. 
Let $\psi : X' \to X''$ be the contraction morphism induced by $A_{X'}$, 
which is birational. 
In particular, $X''$ is toric. 
After replacing $A_Y$ by $mA_Y$ for some integer $m>0$, 
we may assume that $A_{X'} = \psi^*A_{X''}$ for some ample Cartier divisor $A_{X''}$ on $X''$. 
We get that 
\[
X'' = \Proj\,R(X'', A_{X''}) = \Proj\, R(X', A_{X'}) = 
\Proj\, R(X, A_{X}) = \Proj\,R(Y, A_Y) = Y. 
\]
Therefore, $Y$ is toric. 
\end{proof}

\begin{prop}\label{p-p=2-g=9-main}
We use Notation \ref{n-p=2-g=9}. 
Then $\deg \theta \neq 2$. 
\end{prop}

\begin{proof}
Suppose $\deg \theta =2$. 
Let us derive a contradiction. 
Recall that we have the following diagram consisting of contractions of projective normal varieties: 
\[
\begin{tikzcd}
Y \arrow[d, "\sigma"'] \arrow[rd, "\psi"]& & Y^+ \arrow[ld, "\psi^+"'] \arrow[d, "\tau"]\\
X & Z & W=\P^1.
\end{tikzcd}
\]
Since $\theta : Z \to \overline{Z}$ is a finite purely inserapble morphism, 
we have the corresponding diagram consisting of contractions of projective normal varieties: 
\[
\begin{tikzcd}
\ol{Y} \arrow[d, "\sigma'"'] \arrow[rd, "\psi'"]& & \ol{Y}^+ \arrow[ld, "\psi'^+"'] \arrow[d, "\tau'"]\\
\ol{X} & \ol{Z} & \ol{W}=\P^1, 
\end{tikzcd}
\]
where $\sigma', \psi', \psi'^+$ are birational. 
Specifically, this diagram is obtained as follows. 
We have the following decomposition of the absolute Frobenius morphism $F_Z : Z \to Z$ of $Z$: 
\[
F_Z : Z \xrightarrow{\theta} \overline{Z}\to  Z =: Z_1, 
\]
which corresponds to  field extensions 
\[
K(Z) \supset K(\overline{Z}) \supset K(Z)^2 =K(Z_1). 
\]
For the same morphism $\psi_1 : Y_1 \to Z_1$ as $\psi : Y \to Z$, 
we define $\overline Y$ as the normalisation of $Y_1$ in $K(\ol{Z})$. 
By a similar argument, we get  the above diagram. 

Each of $\overline{Y}$ and $\overline{Y}^+$ has exactly two contractions as above. 
In particular, 
$\overline{Y}^+$ has a contraction to $\P^1$, 
whilst $\overline{Y}$ does not. 
Then  $\overline{Y}^+$ is a smooth projective toric threefold 
by Lemma \ref{l-014} and Lemma \ref{l-023} below. 
We have two del Pezzo fibrations 
$V_1:=Y^+ \to \P^1$ and $V_2:= \overline{Y}^+ \to \P^1$, which induce the following commutative diagrams: 
\[
\begin{tikzcd}
V_1 \arrow[r, "\alpha"] \arrow[d, "\tau_1"] &V_2\arrow[d, "\tau_2"] \\
\P^1_1 :=\P^1 \arrow[r, "\beta"] &\P^1_2 :=\P^1 
\end{tikzcd}
\qquad\qquad
\begin{tikzcd}
V_{1, K_1} \arrow[r, "\wt{\alpha}"] 
\arrow[d, "\wt{\tau}_1"] &
V_{2, K_2} \arrow[d, "\wt{\tau}_2"] \\
\Spec K_1 \arrow[r, "\wt{\beta}"] & \Spec K_2. 
\end{tikzcd}
\]
Here 
$K_i$ denotes the function field of $\P^1_i$ for each $i \in \{1, 2\}$ and 
the right diagram is obtained by taking the generic fibres 
over $\P^1_2 = \P^1$. 
Note that we have $V_{1, K_1} := V_1 \times_{\P^1_1} \Spec K_1 = V_1 \times_{\P^1_2} \Spec K_2$. 

We now prove (i)-(v) below. 
\begin{enumerate}
    \item[(i)] 
    For each $i \in \{1, 2\}$, $V_{i, K_i}$ is a regular del Pezzo 
    surface with $\rho(V_{i, K_i})=1$ and $H^0(V_{i, K_i}, \MO_{V_{i, K_i}}) = K_i$. 
    \item[(ii)] $\wt{\alpha} : V_{1, K_1} \to V_{2, K_2}$ is a finite surjective morphism of degree two. 
    \item[(iii)] 
    There exists a Cartier divisor $L$ on $V_2$ 
    such that we have 
    $K_{V_1} \sim \alpha^*(K_{V_2} + L)$ and an exact sequence 
    \begin{equation}\label{e1-p=2-g=9-main}
    0 \to \MO_{V_2} \to \alpha_*\MO_{V_1} \to \MO_{V_2}(-L) \to 0.
    \end{equation} 
    Moreover, for $M:= L|_{V_{2, K_2}}$, 
    we have 
        $K_{V_{1, K_1}} \sim \wt{\alpha}^*(K_{V_{2, K_2}} + M)$ and an exact sequence 
    \begin{equation}\label{e2-p=2-g=9-main}
    0 \to \MO_{V_{2, K_2}} \to \wt{\alpha}_*\MO_{V_{1, K_1}} \to 
    \MO_{V_{2, K_2}}(-M) \to 0.
    \end{equation} 
    \item[(iv)] $K_{V_{1, K_1}}^2 =6$ and $1 \leq K_{V_{2, K_2}}^2 \leq 9$. 
    Note that each intersection number $K_{V_{i, K_i}}^2 = K_{V_{i, K_i}} \cdot_{K_i} 
    K_{V_{i, K_i}}$ is defined over $K_i$, i.e., 
    the equality 
    $\chi(V_i, \MO_{V_i}(-K_{V_{i, K_i}})) = \chi(V_i, \MO_{V_i}) + K_{V_{i, K_i}}^2$ holds for $\chi(V_{i, K_i}, -) := \sum_{j\geq 0} \dim_{K_i} H^j(V_{i, K_i}, -)$.   
    \item[(v)] Either $K_1 = K_2$ or 
 $K_1 \supset K_2$ is an inseparable extension of degree two. 
\end{enumerate}
Fix $i \in \{1, 2\}$. 
Let us show (i). 
Since $\tau_i : V_i \to \P^1_i$ is a del Pezzo fibration with $\rho(V_i/\P^1_i)=1$, 
$V_{i, K_i}$ is a regular del Pezzo surface with $\rho(V_{i, K_i})=1$. 
The equality $H^0(V_{i, K_i}, \MO_{V_{i, K_i}})=K_i$ 
follows from $(\tau_i)_*\MO_{V_i} = \MO_{\P^1_i}$ by the flat base change theorem. 
Thus (i) holds. 
The assertion (ii) follows from $K(V_i) = K(V_{i, K_i})$. 
As for (iii), we can find a Cartier divisor $L$ on $V_2$ satisfying 
an exact sequence (\ref{e1-p=2-g=9-main}) \cite[Lemma A.1]{Kaw2}, which implies  $K_{V_1} \sim \alpha^*(K_{V_2} + L)$ \cite[Proposition 0.1.3]{CD89}. 
We have $K^2_{V_{1, K_1}}=6$ and $1 \leq K^2_{V_{2, K_2}}\leq 9$ by 
\cite[Theorem 7.4]{FanoII} and \cite[Theorem 1.2(2)]{Tan24}, respectively. 
Hence (iv) holds. 
Let us show (v). 
It is enough to show 
$K_1 \subset K_2^{1/2}$. 
Pick $a \in K_1$. 
Then $a^2 \in K(V_1)^2 \subset K(V_2)$. 
Since $a^2$ is algebraic over $K(\P^1_2) =K_2$ and 
the field extension $K(V_2)/K(\P^1_2)$ is algebraically closed, 
we get $a^2 \in K(\P^1_2) = K_2$, i.e., $a \in K_2^{1/2}$. 
This completes the proof of (i)-(v). 

\medskip

We have $\Pic V_{2, K_2} = \Z H$ for some ample Cartier divisor $H$ 
on $V_{2, K_2}$. 
In particular, we can write 
$-(K_{V_{2, K_2}} +M) \sim aH$ for some $a \in \Z$. 
By (iii), we get $a>0$. 
It holds that 
\[
6[K_1:K_2] = 
[K_1:K_2] K_{V_1} \cdot_{K_1} K_{V_1} =
K_{V_1} \cdot_{K_2} K_{V_1} = (\wt{\tau}(K_{V_{2, K_2}} +M))^2 = 2 (aH)^2 = 2a^2 H^2. 
\]
By $[K_1:K_2] \in \{1, 2\}$, 
we get $a=1$ and $H^2 = 3[K_1:K_2] \in \{3, 6\}$. 
We have $-K_{V_{2, K_2}} \sim bH$ for some integer $b>0$. 
By $1 \leq K_{V_{2, K_2}}^2 \leq 9$ and 
$K_{V_{2, K_2}}^2 = (bH)^2 =b^2 H^2 = 3b^2[K_1:K_2]$, 
we get $b=1$, i.e.,  $-K_{V_{2, K_2}} \sim H$. 
It holds that 
\[
L|_{V_{2, K_2}} = M \sim aH -(-K_{V_{2, K_2}}) \sim  aH - bH =0. 
\]
Then $L$ descends via $\tau_2 : V_2 \to \P^1_2 =\P^1$, i.e., 
$L \sim \tau_2^*\MO_{\P^1}(n)$ for some $n \in \Z$. 
By the exact sequence (\ref{e1-p=2-g=9-main}) and $H^1(V_2, \MO_{V_2})=0$, we get $n>0$. 
It holds that 
\[
-\alpha^*K_{V_2} \sim -K_{V_1} + \alpha^*L \sim -K_{V_1} + 
\alpha^*\tau_2^*\MO_{\P^1}(n). 
\]
Recall that $-K_{V_1}$ is nef and big and 
inducing the flopping contraction $V_1 = Y^+ \to Z$. 
On the other hand, $\alpha^*\tau_2^*\MO_{\P^1}(1)$ is positive on the flopping curves on $V_1 = Y^+$. 
Hence $-\alpha^*K_{V_2} \sim -K_{V_1} + 
\alpha^*\tau_2^*\MO_{\P^1}(n)$ is ample, and hence so is $-K_{V_2}$. 
Therefore, $V_2$ is a smooth Fano threefold. 
On the other hand, $V_2$ has a small birational contraction 
$V_2 =\overline{Y}^+ \to \overline Z$, 
which is impossible by the classification of contractions of extremal rays \cite[(1.1.2) in Main Theorem (1.1)]{Kol91}.
\qedhere

\end{proof}




\begin{lem}\label{l-014}
Set $\overline{Z} := S_{\P^1}(0, 1, 4)$. 
Then there exist projective normal $\Q$-factorial toric threefolds $V_1$ and $V_2$ {that} satisfy the following properties. 
\begin{enumerate}
\item If $f : V \to \overline{Z}$ is a small birational morphism from 
a projective normal $\Q$-factorial threefold $V$, 
then $V \simeq V_1$ or $V \simeq V_2$. 
\item $V_1$ is not smooth and $V_2$ is smooth. 
\item There exists a morphism $\pi : V_2 \to \P^1$ with 
$\pi_*\MO_{V_2} = \MO_{\P^1}$. 
\end{enumerate}
\end{lem}


\begin{proof}
Recall that $\overline{Z} = S_{\P^1}(0, 1, 4)\subset \P^{7}$ is the projective cone over $S := S_{\P^1}(1, 4) = \P_{\P^1}(\MO_{\P^1}(1) \oplus \MO_{\P^1}(4)) \subset \P^6$ {(Remark \ref{r rat scroll})}, 
where {the closed embedding} $S \subset \P^6$ is {given}  by 
the tautological bundle $A_S := \MO_S(1)$ of the $\P^1$-bundle structure 
$\pi : S= \P_{\P^1}(\MO_{\P^1}(1) \oplus \MO_{\P^1}(4)) \to \P^1$. 
In what follows, we consider $A_S$ as a Cartier divisor {on $S$}. 
Let $\Gamma_1$ and $\Gamma_2$ be the sections of $\pi : S \to \P^1$ 
corresponding to 
the projections $\MO_{\P^1}(1) \oplus \MO_{\P^1}(4) \to \MO_{\P^1}(1)$ 
and $\MO_{\P^1}(1) \oplus \MO_{\P^1}(4) \to \MO_{\P^1}(4)$, respectively. 
Then $\Gamma_1 \cap \Gamma_2 = \emptyset$, 
$A_S|_{\Gamma_1} \simeq \MO_{\P^1}(1)$, and $A_S|_{\Gamma_2} \simeq \MO_{\P^1}(4)$. 
We have 
$A_S \sim \Gamma_2 + b F$ 
for some $b \in \Z$ and a fibre $F$ of {$\pi : S \to \P^1$}. 
{We then} get $b = {(\Gamma_2 + b F) \cdot \Gamma_1} = A_S \cdot \Gamma_1 = 
\deg (A_S|_{\Gamma_1}) = 1$, and hence  
\[
A_S = \Gamma_2 + F. 
\]
Note that $4 = \deg (A_S|_{\Gamma_2}) = A_S \cdot \Gamma_2= \Gamma_2^2 +1$, i.e., $\Gamma_2^2 = 3$ and $\Gamma_1^2 = -3$.

It is well known that 
{we have a  toric structure $S = X_{\Sigma}$ given by the following fan:}   
\[
\Sigma = \langle 
e_1, e_2,-e_1 +3e_3,  -e_2 \rangle \subset \R^2, 
\]
{where $e_1$ and $e_2$ denote the standard basis of $\R^2$}. 
We now compute the intersection numbers following notation as in  
\cite[Page 337]{CLS11}. 
Set 
\[
u_1 := e_1,\quad 
u_2 :=e_2,\quad 
u_3 := -e_1 + 3e_2,\quad  
u_4 := -e_2. 
\]
Let $D_{\rho_i}$ be the corresponding torus-invariant prime divisors on $X_{\Sigma}$. 
By the toric intersection theory \cite[Lemma 12.5.2]{CLS11}, 
we get 
\[
D_{\rho_1} \cdot D_{\rho_3} =0,\qquad  
D_{\rho_2} \cdot D_{\rho_4} =0
\]
\[
D_{\rho_1} \cdot D_{\rho_2} = 
D_{\rho_2} \cdot D_{\rho_3} = D_{\rho_3} \cdot D_{\rho_4} = D_{\rho_4} \cdot D_{\rho_1} =1. 
\]
By {the following linear equivalences \cite[Proposition 4.1.2]{CLS11}}: 
\[
0 \sim {\rm div}(\chi^{e_1^*}) = \sum_{\rho \in \Sigma(1)}
\langle e_1^*, u_{\rho}\rangle D_{\rho}
= D_{\rho_1} - D_{\rho_3}
\]
\[
0 \sim  {\rm div}(\chi^{e_2^*}) = \sum_{\rho \in \Sigma(1)}
\langle e_2^*, u_{\rho}\rangle D_{\rho}
= D_{\rho_2} +3 D_{\rho_3} -D_{\rho_4},  
\]
we get 
\[
D_{\rho_1}^2 =D_{\rho_3}^2 =0, \qquad 
D_{\rho_2}^2 = -3, \qquad D_{\rho_4}^2 = +3. 
\]
Hence 
$D_{\rho_1}$ and $D_{\rho_3}$ are fibres of $\pi$, 
$D_{\rho_2} =\Gamma_1$, and $D_{\rho_4} = \Gamma_2$. 


Note that $\overline Z$ is obtained by 
contracting a  section of $U := \P_S(\MO_S \oplus \MO_S(A_S))$ to one point. 
This is because 
$\overline{Z} \subset \P^{7}$ is the projective cone over $S \subset \P^6$ whose embedding is induced by $|A_S|$. 
We now compute the toric structure of $U$ following \cite[Page 337]{CLS11}.
Set $D_0 := 0$ and $D_1 := A_S = \Gamma_2 +F = D_{\rho_1} + D_{\rho_4}$. 
In ${\R^2_{e'_0, e'_1}} = \R e'_0 \oplus \R e'_1$, 
we set $F_0 := \R_{\geq 0} e'_1$ {and} 
$F_1 := \R_{\geq 0} e'_0$. 
{For
$(a_{0\rho_1}, a_{0\rho_2}, a_{0\rho_3}, a_{0\rho_4}):=(0, 0, 0, 0)$ and $(a_{1\rho_1}, a_{1\rho_2}, a_{1\rho_3}, a_{1\rho_4}):=(1, 0, 0, 1)$, 
we get 
$D_i = \sum_{\rho \in \Sigma(1)} a_{i\rho} D_{\rho}$.} 
For {a two-dimensional cone} $\sigma  \in \Sigma(2)$ and $j \in \{0, 1\}$, 
we consider 
\[
\Cone(u_{\rho} + a_{0\rho}e'_0+a_{1\rho}e'_1 \,|\, 
\rho \in \sigma(1)) + F_j \subset \R_{e_1, e_2}^2 \times \R^2_{e'_0, e'_1}. 
\]
{Take the canonical $\R$-linear surjection $\R_{e_1, e_2}^2 \times \R^2_{e'_0, e'_1} \to  \R_{e_1, e_2}^2 \times \R^1$}, 
where $\R^1 := \R^2_{e'_0, e'_1}/\R(e'_0+e'_1) = \R e_3$ for 
$e_3 :=e'_0 =-e'_1$ in $\R^1$. 
We have $a_{0\rho} =0$ by definition, and hence  the image cone is given by 
\[
\Cone(u_{\rho} -a_{1\rho}e'_0 \,|\, 
\rho \in \sigma(1)) \pm \R_{\geq 0} e_3 
\subset 
\R e_1 \times \R e_2 \times \R e_3
\]
{For $u_i := u_{\rho_i}$ and each $\rho \in \{ \rho_1, \rho_2, \rho_3, \rho_4\}$, 
the following hold:} 
\begin{itemize}
\item  
$\rho = \rho_1$: 
$u_{\rho_1} -a_{1\rho_1}e_3 = u_1-e_3 =e_1 - e_3$. 
\item 
$\rho = \rho_2$: 
$u_{\rho_2} -a_{1\rho_2}e_3 = u_2 = e_2$. 
\item 
$\rho = \rho_3$: 
$u_{\rho_3} -a_{1\rho_3}e_3 = u_3 = -e_1 +3e_2$.
\item 
$\rho = \rho_4$: 
$u_{\rho_4} -a_{1\rho_4}e_3 = u_4 - e_3 = -e_2 -e_3$.
\end{itemize}
Then the resulting $\P^1$-bundle {$U=X_{\widetilde{\Sigma}}$}  
is obtained by {the} $6$ rays generated by 
\[
v_1 := u_1-e_3 = (1, 0, -1), 
\quad v_2 := u_2 = (0, 1, 0),\quad  
v_3 := u_3 =(-1, 3, 0),\quad  
\]\[
v_4 := u_4-e_3 = (0, -1, -1), \quad 
v_5 := +e_3=(0, 0, 1), \quad 
v_6 := - e_3 = (0, 0, -1). 
\]
The $8$ cones of $\wt{\Sigma}$ are given by 
\[
\Cone(v_1, v_2, v_5), \quad 
\Cone(v_2, v_3, v_5), \quad 
\Cone(v_3, v_4, v_5), \quad 
\Cone(v_4, v_1, v_5), 
\]
\[
\Cone(v_1, v_2, v_6), \quad 
\Cone(v_2, v_3, v_6), \quad 
\Cone(v_3, v_4, v_6), \quad 
\Cone(v_4, v_1, v_6).  
\]
Let $D_{v_1}, ..., D_{v_6}$ be the torus-invariant prime divisors on $X_{\wt{\Sigma}}$. 
Then $D_{v_1}, ..., D_{v_4} \subset X_{\wt{\Sigma}}$ are the inverse images of $D_{u_1}, ..., D_{u_4} \subset X_{\Sigma}$. 
Hence $D_{v_5}$ and $D_{v_6}$ are the sections of $X_{\wt{\Sigma}} \to X_{\Sigma}$. 

Recall that  $f: X_{\wt{\Sigma}} =U \to X_{\Delta} = \overline{Z}$ 
is the contraction of either $D_{v_5}$ or $D_{v_6}$. 
We have $\Ex(f) = D_{v_6}$, because $D_{v_5}$ can not be contracted 
as otherwise the resulting fan $\Delta$ would have only rays contained in the half-space $\{ (x, y, z) \in \R^3\,|\, z \leq 0\}$, 
which contradicts the properness of $X_{\Delta}$. 
Then 
\[
\Delta(1) = \{ \rho_1, \rho_2, \rho_3, \rho_4, \rho_5\}. 
\]
Then {$\Delta$} has exactly the following $5$ cones: 
\[
C_1 := \Cone(v_1, v_2, v_5), \quad 
C_2 := \Cone(v_2, v_3, v_5), \quad 
\]
\[
C_3 :=\Cone(v_3, v_4, v_5), \quad 
C_4 := \Cone(v_4, v_1, v_5), \quad 
\]
\[
\Cone(v_1, v_2, v_3, v_4). 
\]
Recall that a small $\Q$-factorialisation corresponds to a  simplicial division without adding rays. 
Hence there are exactly two possibilities corresponding to 
\[
\Cone(v_1, v_2, v_3, v_4) = \Cone(v_1, v_3, v_2) \cup \Cone(v_1, v_3, v_4) = \Cone(v_2, v_4, v_1) \cup \Cone(v_2, v_4, v_3). 
\]
Set $V_1 := X_{\Delta_1}$ and $V_2 := X_{\Delta_2}$, where 
\[
\Delta_1 := \la C_1, C_2, C_3, C_4,  \Cone(v_1, v_3, v_2),  \Cone(v_1, v_3, v_4)\ra. 
\]  
\[
\Delta_2 := \la C_1, C_2, C_3, C_4,  \Cone(v_2, v_4, v_1),  \Cone(v_2, v_4, v_3)\ra 
\]
Then it is easy to check that $V_1$ is not smooth and $V_2$ is  smooth. 
Finally, there exists a morphism $\pi : V_2 \to \P^1$ 
satisfying $\pi_*\MO_{V_2} = \MO_{\P^1}$, 
because the three points $v_2 = (0, 1, 0), v_4=(0, -1, -1), v_5 = (0, 0, 1)$ lie on the plane 
$\{ (0, y, z) \in \R^3\,|\, y, z \in \R\}$ passing through the origin $(0, 0, 0)$. 
\end{proof}

\begin{lem}\label{l-023}
Set $\overline{Z} := S_{\P^1}(0, 2, 3)$. 
Then there exist projective normal $\Q$-factorial toric threefolds $V_1$ and $V_2$ {that} satisfy the following properties. 
\begin{enumerate}
\item If $f : V \to \overline{Z}$ is a small birational morphism from 
a projective normal $\Q$-factorial threefold $V$, 
then $V \simeq V_1$ or $V \simeq V_2$. 
\item $V_1$ is not smooth and $V_2$ is smooth. 
\item There exists a morphism $\pi : V_2 \to \P^1$ with 
$\pi_*\MO_{V_2} = \MO_{\P^1}$. 
\end{enumerate}
\end{lem}


\begin{proof}
Recall that $\overline{Z} = S_{\P^1}(0, 2, 3)\subset \P^{N+1}$ is the projective cone over $S := S_{\P^1}(2, 3) = \P_{\P^1}(\MO_{\P^1}(2) \oplus \MO_{\P^1}(3)) \subset \P^N$ (Remark \ref{r rat scroll}), 
where the closed embedding $S \subset \P^N$ is given by 
the tautological bundle $A_S$ of the $\P^1$-bundle structure 
$\pi : S= \P_{\P^1}(\MO_{\P^1}(2) \oplus \MO_{\P^1}(3)) \to \P^1$. 
In what follows, we consider $A_S$ as a Cartier divisor on $S$. 
Let $\Gamma_1$ and $\Gamma_2$ be the sections of $\pi : S \to \P^1$ 
corresponding to 
the projections $\MO_{\P^1}(2) \oplus \MO_{\P^1}(3) \to \MO_{\P^1}(2)$ 
and $\MO_{\P^1}(2) \oplus \MO_{\P^1}(3) \to \MO_{\P^1}(3)$, respectively. 
Then $\Gamma_1 \cap \Gamma_2 = \emptyset$, 
$A_S|_{\Gamma_1} \simeq \MO_{\P^1}(2)$, and $A_S|_{\Gamma_2} \simeq \MO_{\P^1}(3)$. 
We have 
\[
A_S \sim \Gamma_2 + b F
\]
for some $b \in \Z$ and a fibre of $\pi : S \to \P^1$. 
We then get $b = (\Gamma_2 + bF) \cdot \Gamma_1 = 
A_S \cdot \Gamma_1 = 
\deg (A_S|_{\Gamma_1}) = 2$, and hence  
\[
A_S = \Gamma_2 + 2F. 
\]
Note that $3 = \deg (A_S|_{\Gamma_2}) = A_S \cdot \Gamma_2= \Gamma_2^2 +2$, i.e., $\Gamma_2^2 = 1$ and $\Gamma_1^2 = -1$.

It is well known that 
we have a torc structure $S = X_{\Sigma}$ 
given by the following fan 
\[
\Sigma = \langle 
e_1, e_2, -e_1 +e_2, -e_2 \rangle \subset \R^2. 
\]
We now compute the intersection numbers following notation as in  \cite[Page 337]{CLS11}. 
Set 
\[
u_1 := e_1, u_2 :=e_2, u_3 := -e_1 + e_2, 
u_4 := -e_2. 
\]
Let $D_{\rho_i}$ be the corresponding torus-invariant prime divisors on $X_{\Sigma}$. 
By the toric intersection theory \cite[Lemma 12.5.2]{CLS11}, 
we get 
\[
D_{\rho_1} \cdot D_{\rho_3} =0, D_{\rho_2} \cdot 
D_{\rho_4} =0
\]
\[
D_{\rho_1} \cdot D_{\rho_2} = 
D_{\rho_2} \cdot D_{\rho_3} = D_{\rho_3} \cdot D_{\rho_4} = D_{\rho_4} \cdot D_{\rho_1} =1. 
\]
By the following linear equivalences \cite[Proposition 4.1.2]{CLS11}: 
\[
0 \sim {\rm div}(\chi^{e_1^*}) = \sum_{\rho \in \Sigma(1)}
\langle e_1^*, u_{\rho}\rangle D_{\rho}
= D_{\rho_1} - D_{\rho_3}
\]
\[
0 \sim  {\rm div}(\chi^{e_2^*}) = \sum_{\rho \in \Sigma(1)}
\langle e_2^*, u_{\rho}\rangle D_{\rho}
= D_{\rho_2} + D_{\rho_3} -D_{\rho_4},  
\]
we get 
\[
D_{\rho_1}^2 =D_{\rho_3}^2 =0, \qquad 
D_{\rho_2}^2 = -1, \qquad D_{\rho_4}^2 = +1. 
\]
Hence $D_{\rho_1}$ and $D_{\rho_3}$ are fibres of $\pi$, 
$D_{\rho_2} =\Gamma_1$ and $D_{\rho_4} = \Gamma_2$. 

Note that $\overline Z$ is obtained by 
contracting a  section of $U := \P_S(\MO_S \oplus \MO_S(A_S))$ to one point. 
This is because 
$\overline{Z} \subset \P^{N+1}$ is the projective cone over $S \subset \P^N$ whose embedding is induced by $|A_S|$. 
We now compute the toric structure of $U$ following \cite[page 337]{CLS11}. 
Set $D_0 := 0$ and $D_1 := A_S = \Gamma_2 +2F = 2D_{\rho_1} + D_{\rho_4}$. 
In $\R^2_{e'_0, e'_1} = \R e'_0 \oplus \R e'_1$, 
we set $F_0 := \R_{\geq 0} e'_1, F_1 := \R_{\geq 0} e'_0$. 
For $(a_{0\rho_1}, a_{0\rho_2}, a_{0\rho_3}, a_{0\rho_4}):=(0, 0, 0, 0)$ and $(a_{1\rho_1}, a_{1\rho_2}, a_{1\rho_3}, a_{1\rho_4}):=(2, 0, 0, 1)$, 
we get 
$D_i = \sum_{\rho \in \Sigma(1)} a_{i\rho} D_{\rho}$. 
For a two-dimensional cone $\sigma  \in \Sigma(2)$ and $j \in \{0, 1\}$, 
we consider 
\[
\Cone(u_{\rho} + a_{0\rho}e'_0+a_{1\rho}e'_1 \,|\, 
\rho \in \sigma(1)) + F_j \subset \R_{e_1, e_2}^2 \times \R^2_{e'_0, e'_1}. 
\]
Take the canonical $\R$-linear surjection $
 \R_{e_1, e_2}^2 \times \R^2_{e'_0, e'_1} \to 
\R_{e_1, e_2}^2 \times \R^1$, 
where $\R^1 := (\R^2_{e'_0, e'_1}/\R(e'_0+e'_1) = \R e_3$ for 
$e_3 :=e'_0 =-e'_1$ in $\R^1$. 
We have $a_{0\rho} =0$ by definition, and hence  the image cone is given by 
\[
\Cone(u_{\rho} -a_{1\rho}e_3 \,|\, 
\rho \in \sigma(1)) \pm \R_{\geq 0} e_3 
\subset 
\R e_1 \times \R e_2 \times \R e_3
\]
For $u_i := u_{\rho_i}$ and each $\rho \in \{\rho_1, \rho_2, \rho_3, \rho_4\}$, the following hold: 
\begin{itemize}
\item  
$\rho = \rho_1$: 
$u_{\rho_1} -a_{1\rho_1}e_3 = u_1-2e_3 =e_1 - 2e_3$. 
\item 
$\rho = \rho_2$: 
$u_{\rho_2} -a_{1\rho_2}e_3 = u_2 = e_2$. 
\item 
$\rho = \rho_3$: 
$u_{\rho_3} -a_{1\rho_3}e_3 = u_3 = -e_1 +e_2$.
\item 
$\rho = \rho_4$: 
$u_{\rho_4} -a_{1\rho_4}e_3 = u_4 - e_3 = -e_2 -e_3$.
\end{itemize}
Then the resulting $\P^1$-bundle $U = X_{\widetilde{\Sigma}}$ 
is obtained by the $6$ rays generated by 
\[
v_1 := u_1-2e_3 = (1, 0, -2), 
\quad v_2 := u_2 = (0, 1, 0),\quad  
v_3 := u_3 =(-1, 1, 0),\quad  
\]
\[
v_4 := u_4-e_3 = (0, -1, -1), \quad 
v_5 := +e_3=(0, 0, 1), \quad 
v_6 := - e_3 = (0, 0, -1). 
\]
The $8$ cones of $\wt{\Sigma}$ are given by 
\[
\Cone(v_1, v_2, v_5), \quad 
\Cone(v_2, v_3, v_5), \quad 
\Cone(v_3, v_4, v_5), \quad 
\Cone(v_4, v_1, v_5), 
\]
\[
\Cone(v_1, v_2, v_6), \quad 
\Cone(v_2, v_3, v_6), \quad 
\Cone(v_3, v_4, v_6), \quad 
\Cone(v_4, v_1, v_6).  
\]
Let $D_{v_1}, ..., D_{v_6}$ be the torus-invariant prime divisors on $X_{\wt{\Sigma}}$. 
Then $D_{v_1}, ..., D_{v_4} \subset X_{\wt{\Sigma}}$ are the inverse images of $D_{u_1}, ..., D_{u_4} \subset X_{\Sigma}$. 
Hence $D_{v_5}, D_{v_6}$ are the sections of $X_{\wt{\Sigma}} \to X_{\Sigma}$. 

Recall that  $f: X_{\wt{\Sigma}} =U \to X_{\Delta} = \overline{Z}$ 
is the contraction of either $D_{v_5}$ or $D_{v_6}$. 
We have $\Ex(f) = D_{v_6}$, because $D_{v_5}$ can not be contracted 
as otherwise the resulting fan $\Delta$ has only rays contained in the half-space $\{ (x, y, z) \in \R^3\,|\, z \leq 0\}$, 
which contradicts the properness of $X_{\Delta} = \overline Z$. 
Then 
\[
\Delta(1) = \{ \rho_1, \rho_2, \rho_3, \rho_4, \rho_5\}. 
\]
Then $\Delta$ has exactly the following $5$ cones: 
\[
C_1 := \Cone(v_1, v_2, v_5), \quad 
C_2 := \Cone(v_2, v_3, v_5), \quad 
\]
\[
C_3 :=\Cone(v_3, v_4, v_5), \quad 
C_4 := \Cone(v_4, v_1, v_5), \quad 
\]
\[
\Cone(v_1, v_2, v_3, v_4). 
\]
Recall that a small $\Q$-factorialisation corresponds to a  simplicial division without adding rays. 
Hence there are exactly two possibilities corresponding to 
\[
\Cone(v_1, v_2, v_3, v_4) = \Cone(v_1, v_3, v_2) \cup \Cone(v_1, v_3, v_4) = \Cone(v_2, v_4, v_1) \cup \Cone(v_2, v_4, v_3). 
\]
Set $V_1 := X_{\Delta_1}$ and $V_2 := X_{\Delta_2}$, where 
\[
\Delta_1 := \la C_1, C_2, C_3, C_4,  \Cone(v_1, v_3, v_2),  \Cone(v_1, v_3, v_4)\ra. 
\]  
\[
\Delta_2 := \la C_1, C_2, C_3, C_4,  \Cone(v_2, v_4, v_1),  \Cone(v_2, v_4, v_3)\ra 
\]
Then it is easy to check that $V_1$ is smooth and $V_2$ is not smooth (note that $\Cone(v_1, v_3, v_2)$ is not smooth, because the determinant of the $3 \times 3$ matrix is not in $\{ 1, -1\}$). 
Finally, there exists a morphism $\pi : V_2 \to \P^1$ 
satisfying $\pi_*\MO_{V_2} = \MO_{\P^1}$, 
because the three points $v_2 = (0, 1, 0), v_4=(0, -1, -1), v_5 = (0, 0, 1)$ lie on the plane 
$\{ (0, y, z) \in \R^3\,|\, y, z \in \R\}$ passing through the origin $(0, 0, 0)$. 
\end{proof}

\begin{cor}\label{c -K_X|_S primitive}
Let $X \subset \P^{g+1}$ be an anti-caninically embedded Fano threefold with $g \geq 7$ and $\Pic\,X = \Z K_X$. 
Let $S_\eta$ be the generic member of $|-K_X|$ and 
let $S_{\ol{\eta}}$ be the base change to the algebraic closure. 
Then the pullback $\MO_X(-K_X)|_{S_{\ol{\eta}}}$ of 
$\MO_X(-K_X)$ to $S_{\ol{\eta}}$ is not divisible by $p$ 
in $\Pic S_{\ol{\eta}}$. 
\end{cor}

\begin{proof}
Let $\mathcal S \subset X \times \P$ be the universal family parametrising 
the effective divisors that are linearly equivalent to $|-K_X|$, 
where $\P := \check{\P}^{g+1} = \P(H^0(X, -K_X))$. 
Let $\P^{\circ} \subset \P$ be the open subset parametrising all the smooth members of $|-K_X|$. 
Let $\mathcal S^{\circ}$ be its inverse image, so that we have the induced smooth projective morphism: 
\[
\rho : \mathcal S^{\circ} \to \P^{\circ}. 
 \]
 In particular, any fibre of $\rho$ over a closed point of $\P^{\circ}$ is a smooth K3 surface. 
By Proposition \ref{p-SB-lifting2}, 
Proposition \ref{p p=2 g=9 non-birat}, and 
Proposition \ref{p-p=2-g=9-main}, there exists a fibre $T$ of $\rho$ (i.e., a smooth member of $|-K_X|$) 
such that $-K_X|_T$ is not divisible by $p$ in $\Pic T$.

Note that $S_{\ol{\eta}}$ is nothing but the geometric generic fibre of $\rho$. 
Recall that we have the specialisation map 
\[
\varphi : \Pic\,(S_{\ol{\eta}}) \to \Pic\,T, 
\]
which is a group homomorphism 
satisfying $\varphi(K_X|_{S_{\ol{\eta}}}) = K_X|_T$ 
\cite[the proofs of Proposition 3.3 and Proposition 3.6]{DP12}. 
Since $K_X|_T$ is not divisible by $p$, 
$K_X|_{S_{\ol{\eta}}}$ is not devisible by $p$. 
\end{proof}

\begin{thm}\label{t-SB-lifting2}
Let $X \subset \P^{g+1}$ be an anti-caninically embedded Fano threefold with $g \geq 7$ and $\Pic\,X = \Z K_X$. 
Then $H^1(X, \Omega_X^1(nK_X))=0$ for every $n>0$ and 
$X$ lifts to $W(k)$. 
\end{thm}

\begin{proof}
Let $T$ be the generic member of $|-K_X|$ and 
let $S$ be its base change to the algebraic closure $\kappa$. 
Then $S$ is a member of $|-K_{X_\kappa}|$ for $X_{\kappa} := X \times_k \kappa$. 
Since being superspecial is a closed condition 
(Proposition \ref{p K3 ss locus}),  
$S$ is a smooth K3 surface which is not superspecial (Theorem \ref{t ex non-ss}).  
Moreover, $\MO_{X_\kappa}(-K_{X_\kappa})|_S$ is not divisible by $p$ in $\Pic S$ (Corollary \ref{c -K_X|_S primitive}). 
Then 
\[
H^1(X, \Omega_X^1(nK_X)) \otimes_k \kappa \simeq 
H^1(X_{\kappa}, \Omega^1_{X_\kappa}(nK_{X_\kappa})=0
\]
for every $n>0$ by Proposition \ref{p-SB-lifting}, which implies $H^1(X, \Omega_X^1(nK_X))=0$.
\end{proof}

\section{Fano threefolds of genus $6$}\label{s genus 6}

\subsection{Gushel-Mukai description}

We start by the following well known result. 
We include the proof, as the authors can not find an appropriate reference which works in arbitrary characteristics. 

\begin{prop}\label{p Enriques Babbage}
Let $C$ be a smooth projective curve of genus $g \geq 2$. 
Then the following hold. 
\begin{enumerate}
\item $|K_C|$ is base point free. 
\item One of the following holds. 
\begin{itemize}
\item $|K_C|$ is not very ample and the morphism $C \to D$ onto the image of $\varphi_{|K_C|}$ is a double cover with $D \simeq \P^1_k$. 
In this case, $C$ is called hyperelliptic. 
\item $|K_C|$ is very ample. 
In this case, each of $C$ and the closed immersion $\varphi_{|K_C|} : C \hookrightarrow \P^{g-1}$ is called a canonical curve. 
\end{itemize}
\item 
Assume that $|K_C|$ is very ample and we identify $C$ with the image 
of the closed immersion $\varphi_{|K_C|} : C \hookrightarrow \P^{g-1}$. 
Then one and only one of the following holds. 
\begin{enumerate}
\item $C$ is an intersection of quadric hypersurfaces in $\P^{g-1}$. 
\item $C$ is trigonal, i.e., there exists a finite surjective morphism $C \to \P^1$ of degree $3$. 
\item $C$ is a quintic plane curve, i.e., 
isomorphic to a curve on $\P^2_k$ of degree $5$. 
\end{enumerate}
\end{enumerate}
\end{prop}

\begin{proof}
The assertions (1) and (2) follow from \cite[Ch. IV, Lemma 5.1, Proposition 5.2, Proposition 5.3]{Har77}. 

Let us show (3). 
Assume that $|K_C|$ is very ample and we identify $C$ with the image 
of the closed immersion $\varphi_{|K_C|} : C \hookrightarrow \P^{g-1}$. 
Then at least one of (a)-(c)   holds by \cite[Theorem in the first page]{SD73}. 
It is enough to show (i)-(iii) below. 
\begin{enumerate}
\item[(i)] 
At least one of (a) and (b) does not hold. 
\item[(ii)] At least one of (b) and (c) does not hold. 
\item[(iii)] At least one of (c) and (a) does not hold. 
\end{enumerate}

\medskip

(i), (ii) Assume (b). 
We have a triple cover $f: C \to \P^1$. 
Fix a general fibre $D$ of $f$. 
By Lemma \ref{l trigonal curve}, 
we have $D =P_1+ P_2+ P_3$, 
where $P_1, P_2, P_3$ are distinct points and  lie on a single line $L$ on $\P^{g-1}_k$. 
Then (a) does not hold, as otherwise we would have inclusions 
\[
\{P_1, P_2, P_3\} \subset C \cap L \subset \bigcap_{C \subset Q} (Q \cap L)
\]
and $Q \cap L$ would be an effective Cartier divisor on $L=\P^1$ for some 
quadric hypersurface $Q \subset \P^{g-1}$ containing $C$. 
Thus (i) holds.

Suppose that (c) holds in addition to (b). 
In this case, we have $C \subset V \subset \P^{g-1}_k = \P^5_k$ for the Veronese surface $V$. 
Since $V$ contains no line on $\P^5_k$ and $V$ is an intersection of quadrics, 
we would get 
\[
\{P_1, P_2, P_3\} \subset C \cap L \subset 
V \cap L \subset \bigcap_{V \subset Q} (Q \cap L)
\]
and  $Q \cap L$ would be an effective Cartier divisor on $L=\P^1$ for some 
quadric hypersurface $Q \subset \P^{5}$ containing $V$. 
This is absurd. 
Thus (ii) holds.

\medskip

(iii) 
Suppose that both (a) and (c) hold. Let us derive a contradiction.  
By adjunction formula, we get $\MO_C(K_C) \simeq 
 \MO_{\P^2}(K_{\P^2}+C)|_C =
\MO_{\P^2}(-3+5)|_C = \MO_{\P^2}(2)|_C$. 
Moreover, we have 
\[
H^0(\P^2, \MO(1)) \xrightarrow{\simeq} H^0(C, \MO(1)|_C)\quad\text{and}\quad
H^0(\P^2, \MO(2)) \xrightarrow{\simeq} H^0(C, \MO(2)|_C). 
\]
Then  the closed embedding $C \subset \P^5$  induced by $|K_C|$ 
factors through the Veronese embedding $v: \P^2 \hookrightarrow \P^5$:  
\[
C \subset V \subset \P^5, \qquad V := v(\P^2).  
\]
Since $C$ is an intersection of quadrics: $C = \bigcap_{C \subset Q} Q$, we have 
\[
C \subset V \cap \bigcap_{C \subset Q} Q = \bigcap_{C \subset Q} (V \cap Q), 
\]
where $V \cap Q$ is either equal to $V$ or an effective divisor $\Gamma$ on $V (\simeq \P^2)$. 
Since each of $C$ and $V$ is an intersection of quadrics, 
we can find a quadric hypersurface $Q \subset \P^5$ such that 
$C \subset Q$ and $V \not\subset Q$. 
Set $\Gamma := V \cap Q$. 
Clearly, we have $C \subset \Gamma$. 
Via $v: \P^2 \xrightarrow{\simeq} V$, we have $v^*(\MO_{\P^5}(1)|_V) = \MO_{\P^2}(2)$. 
By $\Gamma \sim \MO_{\P^5}(2)|_V$, 
we get $\Gamma \simeq v^{-1}(\Gamma) \sim v^*(\MO_{\P^5}(2)|_V) \sim \MO_{\P^2}(4)$. 
However, this contradicts 
the scheme-theoretic inclusion 
$C \subset \Gamma$ 
and $C \sim \MO_{\P^2}(5)$. 
Thus (iii) holds. 
\qedhere

\end{proof}

\begin{lem}\label{l trigonal curve}
Let $C \subset \P^{g-1}$ be a trigonal curve of genus $g \geq 2$ and 
let $f : C \to \P^1_k$ be a finite surhective morphism of degree $3$. 
Set $D := f^*\MO_{\P^1}(1)$. 
Then the following hold. 
\begin{enumerate}
\item $h^0(C, D)=2$. 
\item $h^0(C, K_C-D) =g-2$. 
\item $f$ is separable and a general fibre of $f$ lies on a line $L$ of $\P^{g-1}_k$. 
\end{enumerate}
\end{lem}

\begin{proof}
The assertion (1) follows from Clifford's theorem: $2(h^0(C, D) -1) \leq \deg D =3$. 
By the Riemann-Roch theorem, we have 
\[
\chi(C, K_C-D) = \chi(C, \MO_C)+\deg (K_C-D) = (1-g)+(2g-2-3)= g-4. 
\]
We get $h^1(C, K_C-D) = h^0(C, D)\overset{{\rm (1)}}{=}2$ by Serre duality, which implies 
$h^0(C, K_C-D) = g-2$. 
Thus (2) holds.

Let us show (3). 
Fix a general fibre $D$ of $f$. 
By $C \not\simeq \P^1_k$, $f$ is not purely inseparable, and hence $f$ is separable. 
Hence we get $D = P_1+P_2+P_3$ for three distinct points $P_1, P_2, P_3 \in C$. 
Suppose that there exists no line $L$ on $\P^{g-1}_k$ containing $P_1, P_2, P_3$. 
Then the linear subvariety generated by $P_1, P_2, P_3$ is two-dimensional. 
Pick general points $P_4, P_5, ..., P_{g}$ on $C$. 
Then there exists no hyperplane $H$ on $\P^{g-1}$ 
passing through $P_1, ..., P_{g}$. 
By using the induced isomorphism 
\[
H^0(\P^{g-1}, \MO_{\P^{g-1}}(1)) \xrightarrow{\simeq} H^0(C, K_C), 
\]
we see that there exists no effective Cartier divisor $F$ on $C$ 
such that $K_C \sim F$ and $P_1+ \cdots + P_{g} \leq F$. 
On the other hand, we have 
\[
K_C = D + (K_C-D) = (P_1+P_2+P_3) + (K_C-D) 
\]
and $h^0(C, K_C-D) \overset{{\rm (2)}}{=}g-2$. 
Then there exists a member of $|K_C-D|$ 
passing through $P_4, P_5, ..., P_{g}$, 
because 
\[
h^0(C, K_C-D -P_4 - \cdots -P_{g}) \geq h^0(C, K_C-D) - (g-3) = 1.
\]
This is a contradiction. 
\qedhere

\end{proof}

\begin{dfn}
$\CGr(2, 5) \subset \P^{10}$ denotes the cone of $\Gr(2, 5) \subset \P^9$, 
where the closed embedding $\Gr(2, 5) \subset \P^9$ is given by the Pl\"{u}cker embedding. 
\end{dfn}

\begin{thm}\label{t genus 6 lift}
Let $X \subset \P^{7}$ be an anti-canonically embedded Fano threefold 
of genus $g=6$ such that $\Pic\,X = \Z K_X$. 
Then the following hold. 
\begin{enumerate}
    \item $X$ is isomorphic to a complete intersection 
    $\CGr(2, 5) \cap H_1 \cap H_2 \cap H_3 \cap Q$, 
    where each $H_i$ is a hyperplane of $\P^{10}$ and $Q$ is a quadric hypersurface of $\P^{10}$.
    \item $X$ lifts to $W(k)$. 
\end{enumerate}
\end{thm}
\begin{proof}
Since (1) implies (2), it suffices to show (1), 
{ i.e., $X$ is a GM variety in the sense of \cite[Definition 2.1]{DK18}.} 
   In fact, we can take 
   lifts of the embedding $\mathrm{Gr}(2, 5)\subset \P^9_k$ 
   and its cone $\CGr(2, 5) \subset \P^{10}_k$ over $W(k)$. 
   Then a lift of $X$ is given as an intersection of lifts 
   of $H_1$ $H_2$ $H_3$, and $Q$. 


In order to show (1) (i.e., $X$ is a GM variety), 
   it is enough to confirm conditions (a)--(d) in \cite[Theorem 2.3]{DK18}, 
   {because the proof of \cite[Theorem 2.3]{DK18} in this part works also in positive characteristic}. 
   The condition (a) is satisfied by definition.
   Since $-K_X$ is very ample \cite[Theorem 1.1]{FanoI} and $h^0(X,-K_X)=8$  \cite[Proposition 2.6(3)]{FanoII}, (b) holds.
   
   By \cite[Theorem 1.2]{FanoI}, $X$ is an intersection of quadrics.
   Moreover, the map
   \[
   H^0(\P^7,\sO_{\P^7}(2))\to H^0(X,\sO_{X}(2))
   \]
   is surjective by \cite[Theorem 6.2]{FanoI}, and hence
   \[
   h^0(\P^7, \mathcal{I}_X(2))=h^0(\P^7,\sO_{\P^7}(2))- h^0(X,\sO_{X}(2))=36-30=6.
   \]
   Thus (c) holds.
   Finally, we confirm (d). Let $C$ be a smooth curve that is an intersection of two general members of $|-K_X|$. 
   Then $C$ is a smooth curve of genus 6 
   such that $|K_C|$ is very ample, 
   the embedding $C \subset \P^5$ is given by $|K_C|$, and 
   $C$ is an intersection of quadrics of $\P^5$.
   Thus $C$ is not hyperelliptic, nor trigonal, nor a plane quintic (Proposition \ref{p Enriques Babbage}). 
   Then we can apply \cite[Theorem 5.2]{Muk93}, and (d) is verified as in the proof of \cite[Theorem 2.3]{DK18}.
\end{proof}

\begin{thm}\label{t g=6 to Gr(2, 5)}
Let $X \subset \P^{7}$ be an anti-canonically embedded Fano threefold 
of genus $g=6$ such that $\Pic\,X = \Z K_X$. 
Then one of the following holds. 
\begin{enumerate}
\item 
$X$ is isomorphic to a complete intersection 
    $\Gr(2, 5) \cap H_1 \cap H_2  \cap Q$, 
    where each $H_i$ is a hyperplane of $\P^9$ and $Q$ is a quadric hypersurface of $\P^9$.
\item 
There exists a finite surjective morphism $f: X \to V_5$ of degree two, 
where $V_5$ denotes the smooth Fano threefold with $r_{V_5}=2$ and $(-K_{V_5}/2)^3 = 5$. 
\end{enumerate}
\end{thm}

\begin{proof}
By Theorem \ref{t genus 6 lift}(1), we may assume that 
\[
X = \CGr(2, 5) \cap H_1 \cap H_2 \cap H_3 \cap Q \subset \P^{10}. 
\]
For 
\[
\Gr(2, 5) = \Proj\,k[x_0, ..., x_9]/(f_1, ..., f_r) \subset 
\Proj\,k[x_0, ..., x_9]=\P^9, 
\]
we have 
\[
\CGr(2, 5) = \Proj\,k[x_0, ..., x_9, x_{10}]/(f_1, ..., f_r) \subset 
\Proj\,k[x_0, ..., x_9, x_{10}]=\P^{10}.  
\]
Set 
\[
v:= [0:\cdots :0:1] \in \P^{10}, 
\]
which is nothing but the vertex of $\CGr(2, 5)$. 
Since $X$ is smooth and this is a complete intersection, 
we get $v \not\in X$. 
In particular, we have the induced morphism 
\[
\varphi : X \hookrightarrow \CGr(2, 5) \setminus \{v\} \xrightarrow{\pi}  \Gr(2, 5), 
\]
where $\CGr(2, 5) \setminus \{v\} \to  \Gr(2, 5)$ is the induced projection. 
Since every fibre of $\pi$ is one-dimensional (as it is isomorphic to $\mathbb A^1$), 
it follows from $\rho(X)=1$ that $\varphi : X \to \Gr(2, 5)$ is a finite morphism. 
There are the following two cases: 
\begin{enumerate}
\item[(1)] $v \not\in H_1 \cap H_2 \cap H_3$. 
\item[(2)] $v \in H_1 \cap H_2 \cap H_3$. 
\end{enumerate}

(1) Assume $v \not\in H_1 \cap H_2 \cap H_3$. 
Applying a suitable coordinate change,  we may assume 
\[
H_1=\{x_8=0\},\quad H_2 = \{ x_9=0\},\quad H_3 = \{ x_{10}=0\},
\]
\[
Q = \left\{ \sum_{i=0}^7 c_ix^2_i +\sum_{0 \leq i, j\leq 7} d_{i, j} x_ix_j =0\right\}. 
\]
Then $H_1, H_2, Q \subset \P^{10}$ are the cones of corresponding hypersurfaces
$H'_1, H'_2, Q' \subset \P^9$. 
Hence we get 
\[
\varphi(X) \subset \Gr(2, 5) \cap H'_1 \cap H'_2 \cap Q' =:Y. 
\]

\begin{claim*}
The following hold. 
\begin{enumerate}
\renewcommand{\labelenumi}{(\roman{enumi})}
\item $\Gr(2, 5) \cap H'_1 \cap H'_2 \cap Q'$ is a complete intersection, i.e., 
$\dim Y =3$. 
\item $Y$ is an integral scheme. 
\end{enumerate}
\end{claim*}

\begin{proof}[Proof of Claim]
Let us show (i). 
Suppose     $\dim Y \neq 3$, i.e., $\dim Y \geq 4$. 
Then its cone $\wt{Y} := \CGr(2, 5) \cap H_1 \cap H_2 \cap Q$ satisfies 
\[
\dim \wt{Y} \geq 5, 
\]
which leads to the following contradiction:  
\[
3=\dim X = \dim(\CGr(2, 5) \cap H_1 \cap H_2 \cap H_3 \cap Q) 
= \dim (\wt{Y} \cap H_3) \geq 5-1. 
\] 
Thus (i) holds.

Let us show (ii). 
We first prove that $Y$ is reduced. 
Since $X = \wt{Y} \cap H_3$ is regular, 
$\wt{Y}$ is regular around $X$. 
As $X =H_3|_{\wt Y}$ is an ample effective Cartier divisor on ${\wt Y}$, 
$X$ intersects all the irreducible components of $\wt{Y}$. 
Therefore, every irreducible component of $\wt{Y}$ is generically reduced. 
Since $\wt{Y}$ is Cohen-Macaulay, $\wt{Y}$ is reduced. 
Recall that  the induced morphism 
\[
\wt{Y} \setminus \{ v\} \to Y 
\]
is an $\mathbb A^1$-bundle. 
Locally, the corresponding ring homomorphism is given by $A \hookrightarrow A[t]$. 
Since $A[t]$ is reduced, so is $A$, i.e., $Y$ is reduced. 

Next, let us show that $Y$ is irreducible. 
Since $\wt{Y} \cap H_3 =X$ is an integral scheme, 
$X$ is contained in some irreducible component $Z$ of $\wt{Y}$. 
Suppose that there is another irreducible component $Z'$ of $\wt{Y}$. 
Since $\wt{Y}$ is connected, we may assume that $Z \cap Z' \neq \emptyset$. 
As $\wt{Y}$ is regular around $X$, we see that $X$ is disjoint from $Z \cap Z'$. 
However, this implies $X \cap Z' = \emptyset$, which is absurd. 
Hence $\wt{Y}$ is irreducible. 
Then its non-empty open subset $\wt{Y} \setminus \{v\}$ and its image $Y$ are irreducible. 
Thus (ii) holds. 
This completes the proof of Claim. 
\end{proof}

By Claim, we get $\varphi(X) = Y$, which induces a finite surjective morphism 
$\psi: X \to Y$. 
We have 
\[
\MO_Y(1)^3 = 10 = \MO_X(1)^3 =(\deg \psi) \cdot \MO_Y(1)^3. 
\]
Therefore, $\psi: X \to Y$ is a finite birational morphism. 
Since $X$ is normal, $\psi$ coincides with the normalisation of the Gorenstein variety $Y$. 
For the conductor divisor $C$ on $X$, we have 
\[
\omega_X \otimes \MO_X(C) = \psi^*\omega_Y. 
\]
Note that $C$ is an effective Cartier divisor on $X$ whose support coincides with the non-isomorphic locus of $\psi$. 
Since $-K_Y$ is ample by adjunction, we get $ C =0$, and hence $\psi$ is an isomorphism. 

\medskip

(2) Assume $v \in H_1 \cap H_2 \cap H_3$. 
In this case, the problem is reduced to the case when  
\[
H_1 = \{ x_7 =0\}, \quad 
H_2 = \{ x_8 =0\}, \quad 
H_3 = \{ x_9 =0\}. 
\]
In particular, these are the cones of the corresponding hyperplanes 
$H'_1, H'_2, H'_3 \subset \P^9$. 
By the same argument as in Claim,  
$Y := \Gr(2, 5) \cap H'_1 \cap H'_2 \cap H'_3$ is a complete intersection 
and $Y$ is an integral scheme. 
Then the induced morphism 
\[
\psi : X \to Y
\]
is a finite surjective morphism. 
We have 
\[
\MO_X(1)^3 = 10 = 2 \MO_Y(1)^3, 
\]
and hence $\psi$ is a double cover. 

It is enough to show that $Y$ is smooth. 
To this end, it suffices to prove that $\psi : X \to Y$ is flat. 
For 
\[
g_i(x_0, ..., x_6) := f_i (x_0, ..., x_6, 0, 0, 0), 
\]
we get 
\begin{itemize}
\item $Y = \Proj\,k[x_0, x_1, ...x_6]/(g_1, ..., g_r) \subset \P^6$, \item $X \subset \wt{Y} =\Proj\,k[x_0, x_1, ...x_6, x_{10}]/(g_1, ..., g_r) \subset \P^7$, and 
\item $X = \wt{Y} \cap \wt{Q} =
\Proj\,k[x_0, x_1, ...x_6, x_{10}]/(g_1, ..., g_r, h)$ for some quadric hypersurface $\wt{Q} \subset \P^7$ whose defining equation is 
$h$. 
\end{itemize}
Over $D_+(x_0) \subset \P^9$, 
the ring homomorphism corresponding to the morphism $\psi : X \to Y$ 
can be written as 
\[
\alpha : A \to A[X_{10}]/(H), 
\]
where 
we set $A :=  k[X_1, ..., X_6]/(G_1, ..., G_r)$ and $X_i, G_j, H$ are corresponding to $x_i, g, h$, respectively. 
In particular, we have 
\[
H = aX_{10}^2 + bX_{10} +c
\]
for some $a, b, c \in A$. 
In order to show that $\psi : X \to Y$ is flat, 
it suffices to prove $a \in A^{\times}$ by symmetry. 
Fix a maximal ideal $\m$ of $A$. 
It is enough to show $a \not\in \m$, 
i.e., $\dim_k \MO_{\psi^{-1}(\m)}(\psi^{-1}(\m)) =2$. 
Set $L  \subset \P^7$ to be the line passing through the vertex $\wt{v} :=[0: \cdots :0:1]$ whose image to $\P^6$ is $\m$. 
We get 
\[
\psi^{-1}(\m) =  (L \setminus \{\wt{v}\}) \cap X 
= (L \setminus \{\wt{v}\}) \cap \wt{Y} \cap \wt{Q} 
\overset{{\rm (a)}}{=} (L \setminus \{\wt{v}\})  \cap \wt{Q}  \overset{{\rm (b)}}{=} L \cap \wt{Q}, 
\]
where (a) holds by $L \setminus \{\wt{v}\}) \subset L \subset  \wt{Y}$ 
and (b) follows from $\wt{v} \not\in \wt{Q}$. 
Here $\wt{v} \not\in \wt{Q}$ is guaranteed by 
$\wt{v} \not\in X$, $\wt{v} \in \wt{Y}$, and $X = \wt{Y} \cap \wt{Q}$. 
Moreover, $\wt{v} \not\in \wt{Q}$ and $\wt{v} \in L$ imply 
that $L \cap \wt{Q} = \wt{Q}|_{L}$ is an effective Cartier divisor on a line $L = \P^1$ of degree two. 
Hence $L \cap \wt{Q}$ is of length two. 
\end{proof}

\subsection{Separability of the double cover}

\begin{thm}\label{t g=6 separable}
Let $X \subset \P^{7}$ be an anti-canonically embedded Fano threefold 
of genus $g=6$ such that $\Pic\,X = \Z K_X$. 
Let 
\[
f : X \to V_5
\]
be a finite surjective morphism of degree two, 
where $V_5$ denotes the smooth Fano threefold with $r_{V_5}=2$ and $(-K_{V_5}/2)^3 = 5$. 
Then $f$ is separable, i.e., 
the induced field extension $K(X)/K(V_5)$ is  separable. 
\end{thm}

\begin{proof}
Set $V := V_5$. 
Recall that there exists a Cartier divisor $H$ on $V$ such that 
$K_X \sim f^*(K_V+H)$ and we have an exact sequence 
\[
0 \to \MO_V \to f_*\MO_X \to \MO_V(-H) \to 0. 
\]
In particular, $H$ is an ample Cartier divisor on $V=V_5$, and hence $-K_V \sim 2H$. 
In order to show that $f$ is separable, it suffices to show 
\[
c_3(\Omega^1_V(2H)) \neq 0. 
\]
Set $G := \Gr(2, 5)$. 
In what follows, we will use standard notation on Schubert calculus 
(cf. \cite[Section 4.1]{3264}).  
For the closed immersion $\iota : V \hookrightarrow G$, 
we set $\overline{(-)} := \iota^*(-)$ (e.g., $\overline{\sigma}_1 := \iota^*\sigma_1$). 
Recall that the following hold for the universal subbundle $S$ and quotient bundle $Q$ on $G=\Gr(2, 5)$.  
\begin{enumerate}
\item $S$ and $Q$ are locally free sheaves with ${\rm rank}\,S =2$ and 
${\rm rank}\,Q =3$. 
\item We have $\Omega_G^1 \simeq S \otimes Q^*$, 
where $Q^*$ denotes the dual of $Q$. 
\item $c_t(Q^*) = 1- \sigma_1t + \sigma_2t^2 - \sigma_3t^3$ \cite[5.6.2 in page 178]{3264}. 
\item $c_t(S) = 1 -\sigma_1t + \sigma_{1, 1}t^2$ \cite[5.6.2 in page 178]{3264}. 
\item $\sigma_1 = H =\MO_V(1)$. 
\end{enumerate}
Set  $c_i := c_i(\Omega^1_G(2))$ and $\overline{c}_i := \iota^*c_i= c_i(\Omega^1_G(2)|_V)$

\setcounter{step}{0}
\begin{step}\label{s1 g=6 separable}
It holds that 
\[
c_3(\Omega^1_V(2)) = -10 \overline{\sigma}^3_1 +6 \overline{c}_1\overline{\sigma}^2_1 -3 \overline{c}_2\overline{\sigma}_1+ \overline{c}_3. 
\]
\end{step}

\begin{proof}[Proof of Step \ref{s1 g=6 separable}]
For $X_6 := G, X_5:=G \cap (1), X_4 := G \cap (1)^2, X_3 := G \cap (1)^3$ ($\dim X_n =n$), 
we have a conormal exact sequence 
\[
0 \to \MO_{X_{n-1}}(-1) \to \Omega^1_{X_n}|_{X_{n-1}} \to \Omega^1_{X_{n-1}} \to 0. 
\]
Taking the tensor product with $\MO(2)$, we get 
\[
c_t(\Omega^1_{X_{n-1}}(2)) c_t(\MO(1)) = c_t(\Omega^1_{X_n}(2)|_{X_{n-1}}). 
\]
Applying this procedure three times, we get 
\[
c_t(\Omega^1_V(2)) c_t(\MO(1))^3 = c_t(\Omega^1_G(2)|_V). 
\]
For $c_i := c_i(\Omega^1_G(2))$ and $\overline{c}_i := i^*c_i= c_i(\Omega^1_G(2)|_V)$, 
the following holds in the formal power series ring: 
\begin{eqnarray*}
\sum_{i=0}^3c_i(\Omega^1_V(2))t^i &=& 
c_t(\Omega^1_V(2)) = 
c_t(\Omega^1_G(2)|_V)c_t(\MO(1))^{-3} \\
&=& 
(1 + \overline{c}_1t+  \overline{c}_2t^2+  \overline{c}_3t^3 + \cdots +  \overline{c}_6t^6)
(1-\overline{\sigma}_1 t+\overline{\sigma}^2_1t^2 -\overline{\sigma}_1^3t^3)^3\\
&=& 
(1 +  \overline{c}_1t+  \overline{c}_2t^2+  \overline{c}_3t^3 + \cdots +  \overline{c}_6t^6)
(1-3\overline{\sigma}_1 t+6\overline{\sigma}^2_1t^2 -10\overline{\sigma}_1^3t^3)\\
&=& \cdots +(-10 \overline{\sigma}^3_1 +6 \overline{c}_1\overline{\sigma}^2_1 -3 \overline{c}_2\overline{\sigma}_1+ \overline{c}_3) t^3+ \cdots. 
\end{eqnarray*}
This completes the proof of Step \ref{s1 g=6 separable}. 
\end{proof}

In what follows, we use the following formulas 
Recall that the Chern character $\ch E$ of a locally free sheaf $E$ 
of rank $r$ is given as follows \cite[Appendix A, \S 3]{Har77}: 
\[
\ch(E) = r+c_1 + \frac{1}{2}(c_1^2- 2c_2) + \frac{1}{6}(c_1^3 -3c_1c_2 +3c_3)+ \cdots. 
\]
For another locally free sheaf $F$ of finite rank,  we have $\ch(E \otimes F) = \ch(E)\ch(F)$ \cite[Appendix A, \S 3, \S 4]{Har77}.

\begin{step}\label{s2 g=6 separable}
For $c'_i := c_i(\Omega^1_G)$, it holds that 
\begin{align*}
c'_1 &=  - 5\sigma_1,\\
\frac{1}{2}(c'^2_1- 2c'_2) &=   \frac{7}{2}\sigma_1^2 - 3\sigma_{1,1} - 2\sigma_2,   \\
   \frac{1}{6}(c'^3_1 -3c'_1c'_2 +3c'_3)
   &=  -\frac{11}{6} \sigma_1^3 + \frac{5}{2}\sigma_1\sigma_{1,1} + 2\sigma_1\sigma_2 - \sigma_3. 
\end{align*}
\end{step}

\begin{proof}[Proof of Step \ref{s2 g=6 separable}]
We get  
\begin{eqnarray*}
&& 6+c'_1 + \frac{1}{2}(c'^2_1- 2c'_2) + \frac{1}{6}(c'^3_1 -3c'_1c'_2 +3c'_3)+ \cdots\\
&=&\ch(\Omega_G^1) = \ch(S \otimes Q^*) = 
\ch(S) \ch(Q^*) \\
&=& (2 -\sigma_1 + \frac{1}{2} (\sigma_1^2 -2 \sigma_{1, 1})+ \frac{1}{6}(-\sigma_1^3 +3 \sigma_1\sigma_{1, 1})+\cdots)\\
&\times&
(3 -\sigma_1 +\frac{1}{2}(\sigma_1^2 -2 \sigma_2) 
+\frac{1}{6}(-\sigma_1^3+3\sigma_1\sigma_2-3\sigma_3)+ \cdots)\\
&=&6 - 5\sigma_1  + \left( \frac{7}{2}\sigma_1^2 - 3\sigma_{1,1} - 2\sigma_2 \right) 
+ \left( -\frac{11}{6} \sigma_1^3 + \frac{5}{2}\sigma_1\sigma_{1,1} + 2\sigma_1\sigma_2 - \sigma_3 \right) +\cdots
\end{eqnarray*}
This completes the proof of Step \ref{s2 g=6 separable}. 
\end{proof}

\begin{step}\label{s3 g=6 separable}
It holds that 
\begin{align*}
c_1 &= 7\sigma_1,\\
c_2 &= 19 \sigma_1^2 +3\sigma_{1, 1} +2 \sigma_2,\\
c_3 &= 145 \sigma_1^3 +14\sigma_1\sigma_{1, 1} 
+10\sigma_1\sigma_2-2\sigma_3. 
\end{align*}
\end{step}

\begin{proof}[Proof of Step \ref{s3 g=6 separable}]
It holds that 
\begin{eqnarray*}
&& 6+c_1 + \frac{1}{2}(c^2_1- 2c_2) + \frac{1}{6}(c^3_1 -3c_1c_2 +3c_3)+ \cdots\\
&=&\ch(\Omega_G^1(2)) = \ch(\MO_G(2) \otimes \Omega_G^1) = \ch(\MO_G(2))
\ch(\Omega_G^1)  \\
&=& (1+(2\sigma_1) + \frac{1}{2!} (2\sigma_1)^2 + \frac{1}{3!} (2\sigma_1)^3) 
\\
&\times&
 \left(6+c'_1 + \frac{1}{2}(c'^2_1- 2c'_2) + 
 \frac{1}{6}(c'^3_1 -3c'_1c'_2 +3c'_3)+ \cdots\right).
\end{eqnarray*}
We then get 
\begin{align*}
c_1 &= c'_1 + 12 \sigma_1\\
\frac{1}{2}(c_1^2 -2c_2) &= 
\frac{1}{2}(c_1'^2 -2c'_2) + 2\sigma_1 c'_1 + 12\sigma_1^2\\
\frac{1}{6}(c^3_1 -3c_1c_2 +3c_3)& = 
\frac{1}{6}(c'^3_1 -3c'_1c'_2 +3c'_3) + 
\sigma_1(c'^2_1- 2c'_2)+ 
2\sigma_1^2c'_1 + 8\sigma_1^3. 
\end{align*}
In what follows, we compute $c_1, c_2, c_3$ by using these equations together with  Step \ref{s2 g=6 separable}. 
We have 
$c_1 = c'_1 + 12 \sigma_1  = -5 \sigma_1 +12 \sigma_1 = 7\sigma_1$. 
It holds that 
\begin{align*}
&\frac{1}{2}(c_1^2 -2c_2) = 
\frac{1}{2}(c_1'^2 -2c'_2) + 2\sigma_1 c'_1 + 12\sigma^2_1\\
&= \left(\frac{7}{2}\sigma_1^2 - 3\sigma_{1,1} - 2\sigma_2\right) + 2\sigma_1 (-5\sigma_1) +12 \sigma_1^2 \\
\Leftrightarrow&\,\, \frac{49}{2} \sigma_1^2  -c_2 
= \frac{11}{2}\sigma_1^2 - 3\sigma_{1,1} - 2\sigma_2\\
\Leftrightarrow&\,\, c_2 = 19 \sigma_1^2 +3\sigma_{1, 1} +2 \sigma_2.   
\end{align*}
We get 
\begin{eqnarray*}
&&\frac{1}{6}(c^3_1 -3c_1c_2 +3c_3) \\
&=& 
\frac{1}{6}(c'^3_1 -3c'_1c'_2 +3c'_3) + 
\sigma_1(c'^2_1- 2c'_2)+ 
2\sigma_1^2c'_1 + 8\sigma^3_1\\
&=& 
\left(-\frac{11}{6} \sigma_1^3 + \frac{5}{2}\sigma_1\sigma_{1,1} + 2\sigma_1\sigma_2 - \sigma_3\right) + 2\sigma_1 
\left(\frac{7}{2}\sigma_1^2 - 3\sigma_{1,1} - 2\sigma_2\right) 
+50 \sigma_1^3 + 8\sigma_1^3\\
&=& \left(65 - \frac{11}{6}\right)\sigma_1^3 + \frac{-7}{2}\sigma_1\sigma_{1, 1} 
-2\sigma_1\sigma_2-\sigma_3.
\end{eqnarray*}
Therefore, 
\begin{align*}
&\,\,\frac{1}{6}(343 \sigma_1^3 -21\sigma_1 (19 \sigma_1^2 +3\sigma_{1, 1} +2 \sigma_2) +3c_3) = \left(65 - \frac{11}{6}\right)\sigma_1^3 + \frac{-7}{2}\sigma_1\sigma_{1, 1} 
-2\sigma_1\sigma_2-\sigma_3\\
\Leftrightarrow &\,\,
\frac{1}{6}(-56 \sigma_1^3 -63\sigma_1\sigma_{1, 1} -42\sigma_1 \sigma_2 +3c_3) = \left(65 - \frac{11}{6}\right)\sigma_1^3 + \frac{-7}{2}\sigma_1\sigma_{1, 1} 
-2\sigma_1\sigma_2-\sigma_3\\
\Leftrightarrow&\,\, 
\frac{1}{6} \cdot 3c_3 = \frac{145}{2}\sigma_1^3 + 7\sigma_1\sigma_{1, 1} 
+5\sigma_1\sigma_2-\sigma_3\\
\Leftrightarrow&\,\, c_3= 145 \sigma_1^3 +14\sigma_1\sigma_{1, 1} 
+10\sigma_1\sigma_2-2\sigma_3. 
\end{align*}
This completes the proof of Step \ref{s3 g=6 separable}. 
\end{proof}

\begin{step}\label{s4 g=6 separable}
$c_3(\Omega^1_V(2)) \neq 0$.     
\end{step}

\begin{proof}[Proof of Step \ref{s4 g=6 separable}]
By  Step \ref{s1 g=6 separable} and  Step \ref{s3 g=6 separable}, we get 
\begin{align*}
c_3(\Omega^1_V(2)) 
&=-10\overline{\sigma}^3_1 +6\overline{\sigma}^2_1\overline{c}_1 -3\overline{\sigma} \overline{c}_2 +\overline{c}_3\\
&=-10\overline{\sigma}_1^3 +6\overline{\sigma}_1^2(7\overline{\sigma}_1)  -3\overline{\sigma}_1(19 \overline{\sigma}_1^2 +3\overline{\sigma}_{1, 1} +2 \overline{\sigma}_2)\\
&+ (145 \overline{\sigma}_1^3 +14\overline{\sigma}_1\overline{\sigma}_{1, 1} 
+10\overline{\sigma}_1\overline{\sigma}_2-2\overline{\sigma}_3)\\
&= (-10 +42 -57 +145)\overline{\sigma}_1^3 + (14 -9) \overline{\sigma}_1\overline{\sigma}_{1, 1}
+(-6+10)\overline{\sigma}_1\overline{\sigma}_2  -2\overline{\sigma}_3\\
&= 120 \overline{\sigma}_1^3 + 5 \overline{\sigma}_1\overline{\sigma}_{1, 1}
+4\overline{\sigma}_1\overline{\sigma}_2  -2\overline{\sigma}_3. 
\end{align*}
These  intersection numbers on $V=V_5$ can be calculated by 
those on $G = \Gr(2, 5)$ as follows: 
\[
\overline{\sigma}_1^3 = (\sigma_1|_V)^3 = \sigma_1^6, \qquad 
\overline{\sigma}_1 \overline{\sigma}_{1, 1} = \sigma_1^4\sigma_{1, 1}, \qquad 
\overline{\sigma}_1\overline{\sigma}_2 = \sigma_1^4\sigma_2,\qquad 
\overline{\sigma}_3 = \sigma_1^3 \sigma_3. 
\]
Recall that our abbreviated notation means 
$\sigma_i := \sigma_{i, 0}.$ 
By Schubert calculus on $\Gr(2, 5)$ (especially by Pieri's forumla \cite[4.2.4]{3264}), we have 
\[
\sigma_3\sigma_1^3 = \sigma_{3, 0}\sigma_1^3 = \sigma_{3, 3}, \qquad 
\sigma_{2, 1}\sigma_1^3 = (\sigma_{3, 1} + \sigma_{2, 2})\sigma_1^2 = 2\sigma_{3, 3},
\]
\[
\sigma_2 \sigma_1^4 = (\sigma_{3, 0} + \sigma_{2, 1})\sigma_1^3 = 3\sigma_{3, 3}, \quad 
\sigma_{1, 1} \sigma_1^4 = \sigma_{2, 1}\sigma_1^3 = 2 \sigma_{3, 3}, 
\quad 
\sigma_1^6 = (\sigma_{2, 0} + \sigma_{1, 1}) \sigma^4_{1}
=5 \sigma_{3, 3}.  
\]
We then get 
\begin{align*}
c_3(\Omega^1_V(2)) 
 &= 120 \overline{\sigma}_1^3 + 5 \overline{\sigma}_1\overline{\sigma}_{1, 1}
+4\overline{\sigma}_1\overline{\sigma}_2  -2\overline{\sigma}_3\\
&=120 \sigma_1^6 + 5 \sigma^4_1\sigma_{1, 1}
+4\sigma^4_1\sigma_2  -2\sigma_1^3\sigma_3\\
&=(120 \cdot 5 +5 \cdot 2 +4 \cdot 3 - 2) \sigma_{3, 3} \neq 0, 
\end{align*}
where we used the fact that $\sigma_{3, 3}$ is 
rationally equivalent to the class $[P]$ 
for some 
closed point $P$ on $\Gr(2, 5)$. 
This completes the proof of Step \ref{s4 g=6 separable}. 
\end{proof}
Step \ref{s4 g=6 separable} completes the proof of Theorem \ref{t g=6 separable}. 
\end{proof}

\subsection{$h^1(X, \Omega_X^1)=1$ ($g=6$)}

\begin{prop}\label{p g=6 ANV}
Let $X \subset \P^{7}$ be an anti-caninically embedded Fano threefold  
of genus $6$ with $\Pic\,X = \Z K_X$. 
Then the following hold. 
\begin{enumerate}
\item $h^1(X, \Omega_X^1)=1$. 
\item $H^1(X, \Omega_X^1(-A))=0$ for every ample Cartier divisor $A$ on $X$. 
\end{enumerate}
\end{prop}

\begin{proof}
By Theorem \ref{t g=6 to Gr(2, 5)}, one of the following holds. 
\begin{enumerate}
\item[(A)] 
$X \simeq \Gr(2, 5) \cap H_1 \cap H_2  \cap Q$, 
    where each $H_i$ is a hyperplane of $\P^9$ and $Q$ is a quadric hypersurface of $\P^9$.
\item[(B)] 
There exists a finite surjective morphism $f: X \to V$ of degree two, 
where $V$ denotes the smooth Fano threefold with $r_{V}=2$ and $(-K_{V}/2)^3 = 5$. 
\end{enumerate}
{Recall that  both $\Gr(2,5)$ and $V$ are quasi-$F$-split (\cite[Theorem 2]{Mehta-Ramanathan} and \cite[Theorem A]{Kawakami-Tanaka(dPvar)}) and 
every 
smooth quasi-$F$-split variety satisfies Akizuki-Nakano vanishing \cite{Petrov}.}
If (A) holds, then the assertions (1) and (2) hold by Lemma \ref{lem:ANV for CI}. 

In what follows, we treat the case when (B) holds. 
Let us show (1). 
The proof consists of 4 steps.

\setcounter{step}{0}

\begin{step}\label{s1 g=6 ANV}
$h^1(X, \Omega_X^1) \geq 1$. 
\end{step}

\begin{proof}[Proof of Step \ref{s1 g=6 ANV}]
Recall that there exists a lift $\mathcal X$ of $X$ over $W(k)$ (Theorem \ref{t genus 6 lift}). 
Then the assertion follows from the upper semicontinuity of cohomologeis \cite[Ch. III, Theorem 12.8]{Har77}: 
\[
h^{1, 1}(X)  \geq h^{1, 1}(X_{\overline{K}}) =1. 
\]
This completes the proof of Step \ref{s1 g=6 ANV}. 
\end{proof}

\begin{step}\label{s2 g=6 ANV}
$h^1(X, f^*\Omega_V^1) = 1$. 
\end{step}

\begin{proof}[Proof of Step \ref{s2 g=6 ANV}]
We have an exact sequence 
\[
0 \to \MO_V \to f_* \MO_X \to \MO_V(-L) \to 0. 
\]
Taking the tensor product with $\Omega_V^1$, 
this exact sequence induces the following one: 
\[
0 \to \Omega^1_V \to f_*f^*\Omega_V^1 \to \Omega^1_V(-L) \to 0. 
\]
By Akizuki-Nakano vanishing on $V$, 
we have $H^0(V, \Omega^1_V(-L))= H^1(V, \Omega^1_V(-L))=0$, 
which implies  
\[
1 = h^1(V, \Omega^1_V) = h^1(V, f_*f^*\Omega^1_V) = h^1(X, f^*\Omega^1_V). 
\]
This completes the proof of Step \ref{s2 g=6 ANV}. 
\end{proof}

\begin{step}\label{s3 g=6 ANV}
The following hold. 
\begin{enumerate}
\item[(i)] The standard sequence 
\[
0 \to f^*\Omega^1_V \xrightarrow{\alpha} \Omega_X^1 \to \Omega^1_{X/V} \to 0
\]
is exact. 
\item[(ii)] If $H^1(X, \Omega^1_{X/V}) =0$, then $h^1(X, \Omega^1_X) \leq 1$. 
\end{enumerate}
\end{step}

\begin{proof}[Proof of Step \ref{s3 g=6 ANV}]
By Step \ref{s2 g=6 ANV}, (i) implies (ii). 
Hence it suffices to show (i). 
It is well known that the sequence 
\[
f^*\Omega^1_V \xrightarrow{\alpha} \Omega_X^1 \to \Omega^1_{X/V} \to 0. 
\]
is exact  \cite[Ch. II, Proposition 8.11]{Har77}. 
Since $f: X \to V$ is separable, i.e., generically \'etale (Theorem \ref{t g=6 separable}), 
$\alpha : f^* \Omega^1_V \to \Omega_X$ is an isomorphism at the generic point of $X$. 
In particular, $\alpha$ is injective at the generic point. 
This, together with torsion freeness of $f^*\Omega^1_V$, 
implies that $\alpha$ is injective.     
This completes the proof of Step \ref{s3 g=6 ANV}. 
\end{proof}

Let $L$ be an ample Cartier divisor on $V$ satisfying $-K_V \sim 2L$.

\begin{step}\label{s4 g=6 ANV}
$X$ is isomorphic to a closed subscheme of $P:=\P_V(\MO_V \oplus \MO_V(-L))$. 
\end{step}

\begin{proof}[Proof of Step \ref{s4 g=6 ANV}]
By \cite[Definition 2.1, Remark 2.2, Lemma 2.3]{FanoIII}, 
there exists a Cartier divisor $M$ on $V$ such that 
$\MO_X(K_X) \simeq f^*\MO_V(K_V+M)$ and we have an exact sequence 
\[
0 \to \MO_V \to f_*\MO_X \to \MO_V(-M) \to 0.
\]
We then get $L \sim M$, and hence we may assume that $M=L$. 
It holds that $\Ext^1(\MO_V(-L), \MO_V) \simeq H^1(V, \MO_V(L))=0$, 
and hence the above exact sequence splits, i.e., $\MO_V \to f_*\MO_X$ splits. 
By \cite[the proof of Proposition 0.1.3]{CD89} (note that $\MO_V \to f_*\MO_X$ splits if and only if 
$f$ corresponds to a splittable admissible triple), 
there is a closed immersion $j^{\circ} : X \hookrightarrow P^{\circ}$ 
to the $\A^1$-bundle 
\[
P^{\circ} := \Spec_V(\bigoplus_{d=0}^{\infty} \MO_V(-dL)). 
\]
Since $P^{\circ}$ is an open subscheme of $P$, we obtain a closed immersion $j : X \to P$ over $V$. 
This completes the proof of Step \ref{s4 g=6 ANV}. 
\end{proof}

In what follows, we consider $X$ as a closed subscheme on $P = \P_V(\MO_V \oplus \MO_V(-L))$ (Step \ref{s4 g=6 ANV}).

\begin{step}\label{s5 g=6 ANV}
The following hold. 
\begin{enumerate}
\item[(i)] For the coherent ideal sheaf $I := I_{P, X}$ on $P$ 
that defines $X$, 
the standard sequence 
\[
0 \to I/I^2 \xrightarrow{\beta} \Omega^1_{P/V}|_X \to \Omega^1_{X/V}
\to 0 
\]
is exact. 
\item[(ii)] $H^1(X, \Omega^1_{X/V})=0$. 
\end{enumerate}
\end{step}

\begin{proof}[Proof of Step \ref{s5 g=6 ANV}]
We first prove that the implication (i) $\Rightarrow$ (ii). 
Recall that both $I/I^2$ and $\Omega^1_{P/V}|_X$ are invertible sheaves on $X$. 
By $\Pic X = \Z K_X$ and the Kodaira vanishing on $X$, 
we have that $H^i(X, M)=0$ for every $i \in \{1, 2\}$ and invertible sheaf $M$ on $X$. 
In particular, we get $H^1(X, \Omega^1_{P/V}|_X) = H^2(X, I/I^2)=0$, 
which implies $H^1(X, \Omega^1_{X/V})=0$. 
This completes the proof of the implication (i) $\Rightarrow$ (ii).

It suffices to show (i). 
It is well known that the sequence 
\[
I/I^2 \xrightarrow{\beta} \Omega^1_{P/V}|_X \to \Omega^1_{X/V}
\]
is exact \cite[Ch. II, Proposition 8.12]{Har77}. 
By the  same argument as in Step \ref{s3 g=6 ANV}(i), it is enough to prove that $\beta$ is injective at the generic point of $X$. 
Applying the base change by $(-) \times_V \Spec K(V)$ (i.e., taking the generic fibres over $V$), 
it suffices to show the injectivity of $\beta \times_V \Spec K(V)$. 
For the algebraic closure $\overline{K(V)}$ of $K(V)$, 
a field extension $K(V) \to \overline{K(V)}$ is faithfully flat. 
Then it is enough to prove that 
$\beta \times_V \Spec \overline{K(V)}$ is injective, 
which follows from the fact that $P \times_V \Spec \overline{K(V)} \simeq \P^1_{\overline{K(V)}}$ and 
$X\times_V \Spec\,\overline{K(V)}$ is a zero-dimensional reduced closed subscheme on 
$P \times_V \Spec\,\overline{K(V)}$ (Theorem \ref{t g=6 separable}). 
This completes the proof of Step \ref{s5 g=6 ANV}. 
\end{proof}
Step \ref{s1 g=6 ANV}, Step \ref{s3 g=6 ANV}, and  Step \ref{s5 g=6 ANV} complete the proof of (1).

\medskip

Let us show (2). 
Fix an ample Cartier divisor $A$ on $X$. 
By $f^*L \sim -K_X$ and $\Pic X = \Z K_X$, we can write $A \sim f^*A_V$ for some ample Cartier divisor $A_V$ on $V$. 
By the same argument as in Step \ref{s2 g=6 ANV} and Step \ref{s5 g=6 ANV}, 
we get $H^1(X, f^*\Omega_V^1(-A)) = 0$ and $H^1(X, \Omega^1_{X/V}(-A))=0$, respectively. 
Then Step \ref{s3 g=6 ANV}(i) implies $H^1(X, \Omega_X^1(-A))=0$. 
Thus (2) holds. 
\end{proof}

\subsection{$H^0(X, \Omega_X^2)=0$ ($g=6$)}

The purpose of this subsection is to prove 
$H^0(X, \Omega_X^2)=0$ 
for an anti-canonically embedded Fano threefold $X \subset \P^7$ 
of genus $6$ with $\Pic X= \Z K_X$. 
The key part is to prove $H^0(V_5, \Omega_{V_5}^2(1))=0$. 
Although this result is known in characteristic zero, 
we cannot apply the proof as far as the authors know. 
Our proof is totally different from the one in characteristic zero. 
The main idea is to apply two-ray game  in order to  
reduce the computation of $H^0(V_5, \Omega_{V_5}^2(1))=0$ to a similar cohomology on $\P^3$. 

\begin{nota}\label{n V5 2ray}
Let $V$ be a smooth Fano threefold with $r_V =2$ and $(-K_V/2)^3 =5$. 
Fix an ample Cartier divisor $D_V$ satisfying $2D_V \sim -K_V$. 
Take a conic $\Gamma$ on $V$ (i.e., $\Gamma$ is a curve on $V$ with $D_V \cdot \Gamma =2$ and let $g: W \to V$ be the blowup along $\Gamma$. 
It is easy to see that $W$ is a smooth Fano threefold. 
By \cite[Table 5 in Section 9]{FanoIII}, 
$W$ is of No.~2-22. 
Then we have the blowup $f: W \to \P^3$ along a smooth rational curve $B$ on $\P^3$ of degree $4$. 
Set $E :=\Ex(f)$, $D_W := g^*D_V$, and $D_{\P^3} := f_*D_W = f_*g^*D_V$. 
\end{nota}

\begin{lem}\label{l V5 2ray easy}
We use Notation \ref{n V5 2ray}. 
Then the following hold. 
\begin{enumerate}
\item $\MO_{\P^3}(D_{\P^3}) \simeq \MO_{\P^3}(3)$. 
\item $D_W =f^*D_{\P^3} -E$. 
\end{enumerate}
\end{lem}

\begin{proof}
Let us show (1). 
For a hyperplane $H$ on $\P^3$, we have 
\[
-K_W \sim f^*H + g^*D_V = f^*H+D_W. 
\]
We then get 
\[
4H \sim -K_{\P^3} \sim -f_*K_W \sim f_*f^*H + f_*D_W 
= H + D_{\P^3}. 
\]
Hence (1) holds.

Let us show (2). 
We can write 
\[
D_W = f^*D_{\P^3} + n E
\]
for some $n \in \Z$. 
Fix a curve $\ell_f$ on $W$ contracted by $f$. 
Then $E \cdot \ell_f =-1$ and $f^*({\rm divisor}) \cdot \ell_f=0$. 
We have 
\[
1 = -K_W  \cdot \ell_f = (f^*H + D_W) \cdot \ell_f= D_W \cdot \ell_f 
=(f^*D_{\P^3}+nE) \cdot \ell_f = -n. 
\]
Hence $n=-1$. Thus (2) holds. 
\end{proof}

\begin{lem}\label{l V5 to P3}
We use Notation \ref{n V5 2ray}. 
Then the following hold. 
\begin{enumerate}
\item 
It holds that $\wedge^2 (I_B/I^2_B) \simeq \MO_B(-14)$, 
where $\MO_B(-14)$ denotes the invertible sheaf on $B\simeq \P^1$ of degree $-14$. 
\item There exists the following exact sequence: 
\[
0 \to I_B \otimes \Omega^2_{\P^3} \to f_*(\Omega_W^2(-E)) \to 
\wedge^2 (I_B/I^2_B) \to 0.
\]
\item $H^0(V, \Omega_V^2(D_V)) \simeq 
H^0(W, \Omega_W^2(D_W)) \simeq 
H^0(\P^3, I_B \otimes \Omega^2_{\P^3}(D_{\P^3}))$. 
\end{enumerate}
\end{lem}

\begin{proof}
The assertion (1) follows from $\deg N_{B/\P^3} = 2g(B) -2 +(-K_{\P^3}) \cdot B = 0 -2 +4 \cdot 4= 14$ (cf. \cite[Lemma 3.21(2)]{FanoII}). 

Let us show the implication (2) $\Rightarrow$ (3). 
It follows from $g_*\Omega_W^2 = \Omega_V^2$ and $g^*D_V = D_W$ 
that  $H^0(V, \Omega_V^2(D_V)) \simeq 
H^0(W, \Omega_W^2(D_W))$.  
By (1) and (2), we get an exact sequence 
\[
0 \to I_B \otimes \Omega^2_{\P^3} \to f_*(\Omega_W^2(-E)) \to 
\MO_B(-14) \to 0.
\]
Recall that $f^*D_{\P^3} -E = D_W$ and $\MO_{\P^3}(D_{\P^3}) = \MO_{\P^3}(3)$ (Lemma \ref{l V5 2ray easy}). 
Taking the tensor product with $\MO_{\P^3}(D_{\P^3})$, we get the following exact sequence:  
\[
0 \to I_B \otimes \Omega^2_{\P^3}(D_{\P^3}) \to f_*(\Omega_W^2(D_W)) \to 
\MO_B(-2) \to 0, 
\]
which induces $H^0(\P^3, I_B \otimes \Omega^2_{\P^3}(D_{\P^3})) \xrightarrow{\simeq} H^0(W, \Omega_W^2(D_W))$. 
This completes the proof of the implication (2) $\Rightarrow$ (3). 

It is enough to prove (2). 
On $W$ and $E$, we have the following exact sequences 
\begin{equation}\label{e1 V5 to P3}
0 \to \Omega_W^2(-E) \to \Omega_W^2 \to \Omega_W^2|_E \to 0, 
\end{equation}
\begin{equation}\label{e2 V5 to P3}
0 \to \Omega^1_E(-E) \to \Omega_W^2|_E \to \Omega_E^2  \to 0, 
\end{equation}
where the latter one is obtained by applying $\wedge^2$ to 
the conormal exact sequence 
$0 \to \MO_W(-E)|_E \to \Omega_W^1|_E \to \Omega_E^1  \to 0$. 
Since $f|_E : E \to B$ is a $\P^1$-bundle, we get $f_*\Omega_E^2 = f_*\omega_E = 0$. This, together with (\ref{e2 V5 to P3}), induces 
\begin{equation}\label{e3 V5 to P3}
\varphi : f_*(\Omega^1_E(-E)) \xrightarrow{\simeq} f_*(\Omega^2_W|_E).
\end{equation}
By the exact sequence 
\[
0 \to f_*( f^*\Omega_B^1(-E)) \to f_*(\Omega_E^1(-E)) \to f_* (\Omega_{E/B}^1(-E)) =0, 
\]
we obtain 
\[
f_*(\Omega_E^1(-E)) \simeq f_*( f^*\Omega_B^1(-E)) \simeq \Omega_B^1 \otimes_{\MO_{\P^3}} f_*\MO_W(-E) = 
\Omega_B^1 \otimes_{\MO_{\P^3}} I_B \simeq
\Omega_B^1 \otimes_{\MO_B} I_B/I_B^2.  
\]
Therefore, we get 
\begin{equation}\label{e4 V5 to P3}
f_*(\Omega^2_W|_E) \xrightarrow{\simeq, \varphi^{-1}} f_*(\Omega^1_E(-E)) \simeq \Omega_B^1 \otimes I_B/I_B^2.
\end{equation}

\begin{claim*}
    $R^1f_*(\Omega^2_W(-nE))=0$ for every integer $n>0$. 
\end{claim*}

\begin{proof}[Proof of Claim]
We prove the assertion by descending induction on $n$. 
If $n \gg 0$, then the assertion by the Serre vanishing theorem. 
Fix an integer $n>0$. Assume  
$R^1f_*(\Omega^2_W(-(n+1)E))=0$ and let us show $R^1f_*(\Omega^2_W(-nE))=0$. 
By (\ref{e1 V5 to P3}), 
it suffices to show $R^1f_*(\Omega^2_W(-nE)|_E)=0$. 
By (\ref{e2 V5 to P3}),  it is enough to prove 
$R^1f_*(\Omega^1_E(-(n+1)E))=0$ and $R^1f_*\Omega^2_E(-nE)=0$. 
Both of them hold by the fact that $-E$ is $f$-ample and 
$f|_E :E \to B$ is a $\P^1$-bundle. 
This completes the proof of Claim. 
\end{proof}

Applying $f_*$ to (\ref{e1 V5 to P3}), Claim induces the following exact sequence: 
\[
0 \to f_*(\Omega_W^2(-E)) \to f_*\Omega_W^2 \to f_*(\Omega_W^2|_E) \to 0. 
\]
On $\P^3$ and $B$, we have the following exact sequences
\[
0 \to I_B \otimes \Omega^2_{\P^3} \to \Omega^2_{\P^3} 
\to \Omega^2_{\P^3}|_B 
\to 0, 
\]
\[
0 \to \wedge^2 I_B/I_B^2 \to \Omega^2_{\P^3}|_B \xrightarrow{\alpha} \Omega^1_B 
\otimes I_B/I_B^2 \to 0, 
\]
where the latter one is obtained by applying $\wedge^2$ to 
the conormal exact sequence $0 \to I_B/I_B^2 \to \Omega^1_{\P^3}|_B \to \Omega^1_B \to 0$. 
Then we have the following commutative diagram: 
\[
\begin{tikzcd}
0 \arrow[r] & f_*\Omega_W^2(-E) \arrow[r] &  f_*\Omega_W^2 
\arrow[r] & f_*(\Omega_W^2|_E) \arrow[r] &  0\\
0 \arrow[r] & I_B\otimes \Omega_{\P^3}^2 \arrow[r] 
\arrow[u] &  \Omega_{\P^3}^2 
\arrow[r]
\arrow[u, equal]  & \Omega_{\P^3}^2|_B \arrow[r] 
\arrow[u, "\beta"] &  0.
\end{tikzcd}
\]
In order to show (2), it is enough to prove ($\star$) below. 
\begin{enumerate}
\item[($\star$)] $\Ker(\beta) = \wedge^2 (I_B/I_B^2)$. 
\end{enumerate}
Indeed,  ($\star$) and the snake lemma induce the following short exact sequence: 
\[
0 \to I_B\otimes \Omega^2_{\P^3} \to f_*\Omega_W^2(-E) 
\to \wedge^2 (I_B/I_B^2) \to 0. 
\]

\medskip

Therefore, 
it is enough to show ($\star$). 
We have two natural homomorphisms
\[
\alpha : \Omega^2_{\P^3}|_B \to \Omega^1_B \otimes I_B/I_B^2, \qquad 
\beta : \Omega^2_{\P^3}|_B\to  f_{*}(\Omega^2_{W}|_E), 
\]
where $\alpha$ is induced by the conormal sequence for $B \subset \P^3$, 
whilst $\beta$ is the pullback homomorphism. 

Let us show that $\beta(\Ker \alpha)=0$. 
We have $\Ker \alpha = \wedge^2 I_B/I_B^2= (I_B/I_B^2) \wedge  (I_B/I_B^2)$. 
Pick $\zeta \in \Ker \alpha$. 
Then $\zeta  = \sum z dx \wedge dy$ for $z \in \MO_B$ and $x, y \in I_B$, 
where $\sum$ is a finite sum.  
Since our aim is to show $\beta(\zeta)=0$, 
we may assume that $\zeta  =dx \wedge dy$. 
Via the (inverse of) natural isomorphism 
$ f_*(\Omega^2_W|_E)\xrightarrow{\simeq, \varphi^{-1}} f_*(\Omega_E^1 \otimes I_E/I_E^2)$ given in (\ref{e3 V5 to P3}), we get 
\[
\varphi^{-1}(\beta(\zeta)) = \varphi^{-1}((d(f^*x) \wedge d(f^*y))|_E) = d(f^*x|_E) \otimes f^*y =  d0 \otimes f^*y =0,  
\]
where $f^*x|_E=0$ follows from $x \in I_B$. 

 Therefore, we get $\Ker \alpha \subset \Ker \beta$. 
Hence we obtain the induced surjection $\psi : \Im(\alpha) \twoheadrightarrow \Im(\alpha)$. 
Since $\alpha$ and $\beta$ are surjective, 
we get the following isomorphism: 
\[
\Im(\alpha) = \Omega^1_B \otimes (I_B/I_B^2) \overset{(\ref{e4 V5 to P3})}{\simeq} f_*(\Omega_W^2|_E) = \Im(\beta). 
\]
Therefore, $\psi : \Im(\alpha) \twoheadrightarrow \Im(\alpha)$ is a surjective endomorphism (up to isomorphisms), and hence $\psi$ is an isomorphism. 
We then get 
$\Ker \beta = \Ker \alpha = \wedge^2 I_B/I_B^2$, i.e., ($\star$) holds. 
\qedhere

\end{proof}

\begin{prop}\label{p V5 Omega^2(1)}
Let $B$ be a curve on $\P^3$. 
Then 
\[
H^0(\P^3, \Omega^2_{\P^3} \otimes \MO_{\P^3}(3) \otimes I_B)=0, 
\]
where $I_B$ denotes the coherent ideal sheaf on $\P^3$ that defines $B$. 
\end{prop}
\begin{proof}
Applying $\wedge^2$ to the Euler sequence 
\[
0 \to \Omega^1_{\P^3} \xrightarrow{f} \bigoplus_{i=0}^3 
\MO_{\P^3}(-1)e_i \xrightarrow{g} \MO_{\P^3} \to 0, 
\]
we get an exact sequence
\[
0 \to \Omega^2_{\P^3} \xrightarrow{f \wedge f} \bigoplus_{0\leq i < j \leq 3} \MO_{\P^3}(-2) e_i \wedge e_j \xrightarrow{h} \Omega_{\P^3}^1 \to 0. 
\]
Then  the composite homomorphism
\[
\wt{h} := f \circ h : \bigoplus_{0\leq i < j \leq 3} \MO_{\P^3}(-2) e_i \wedge e_j \xrightarrow{h} \Omega_{\P^3}^1 \overset{f}{\hookrightarrow} \bigoplus_{i=0}^3 \MO_{\P^3}(-1)e_i 
\]
can be written as follows: 
\[
\wt{h}(e_i \wedge e_j) =e_i  \otimes g(e_j) - e_j \otimes g(e_i) = e_i \otimes x_j - e_j \otimes x_i 
= x_je_i -x_ie_j, 
\]
where $\P^3 = \Proj\,k[x_0, x_1, x_2, x_3]$ 
and hence $H^0(\P^3, \MO_{\P^3}(1)) = kx_0 \oplus kx_1 \oplus kx_2 \oplus k x_3$. 
Since $f$ is injective, we get an exact sequence  
\[
0 \to H^0(\P^3, \Omega^2_{\P^3}(3)) \to \bigoplus_{0\leq i < j \leq 3} H^0(\P^3, \MO_{\P^3}(1)) e_i \wedge e_j \xrightarrow{\wt{h}} \bigoplus_{0\leq i \leq 3} H^0(\P^3, \MO_{\P^3}(2)) e_i.
\]
It holds that 
\[
h^0(\P^3, \Omega^2_{\P^3}(3)) = 6 \times 
h^0(\P^3, \MO_{\P^3}(1)) 
- h^0(\P^3, \Omega^1_{\P^3}(3)) 
\]
\[
= 6 \cdot 4 - (4 \cdot h^0(\P^3, \MO_{\P^3}(2)) - h^0(\P^3, \MO_{\P^3}(3))) = 24 - ( 4 \cdot 10 -20) = 4. 
\]
For $0 \leq \alpha \leq 3$ and $0 \leq i < j \leq 3$, we have 
\[
\wt{h}(x_{\alpha} e_i \wedge e_j)
= x_{\alpha}x_j e_i - x_{\alpha}x_i e_j. 
\]
By 
\begin{itemize}
\item $\wt{h}(x_0 e_1 \wedge e_2) = x_0x_2e_1 - x_0x_1 e_2$, 
\item $\wt{h}(x_1 e_2 \wedge e_0) = x_1x_0e_2 - x_1x_2 e_0$,
\item $\wt{h}(x_2 e_0 \wedge e_1) = x_2x_1e_0 - x_2x_0 e_1$,
\end{itemize}
we get 
\[
\zeta_{012} := x_0e_1 \wedge e_2 +x_1e_2 \wedge e_0 + x_2 e_0 \wedge e_1 \in \Ker (\wt{h}). 
\]
We define $\zeta_{013},\zeta_{023},\zeta_{123} \in \Ker (\wt{h})$ in the same way. 
By symmetry, it holds that 
\[
\zeta_{012}, \zeta_{013}, \zeta_{023}, \zeta_{123} \in \Ker(\wt{h}).  
\]
Then these elements form a $k$-linear basis of $H^0(\P^3, \Omega^2_{\P^3}(3))$, as they are $k$-linearly independent: 
\[
H^0(\P^3, \Omega^2_{\P^3}(3)) = 
k \zeta_{012} \oplus k \zeta_{013}\oplus k  \zeta_{023} \oplus k  \zeta_{123}.
\]
By construction, we get the following commutative diagram 
in which each horizontal sequence is exact: 
\[
\begin{tikzcd}
0 \arrow[r] & 
H^0(\Omega^2_{\P^3}(3)) \arrow[r] & 
\displaystyle{\bigoplus_{0\leq i < j \leq 3} H^0(\MO_{\P^3}(1)) e_i \wedge e_j} \arrow[r, "\wt{h}"] & 
\displaystyle{\bigoplus_{0\leq i \leq 3} H^0(\MO_{\P^3}(2)) e_i}\\
0 \arrow[r] & 
H^0(I_B\Omega^2_{\P^3}(3)) \arrow[r]\arrow[u, hook] & 
\displaystyle{\bigoplus_{0\leq i < j \leq 3} H^0(I_B\MO_{\P^3}(1)) e_i \wedge e_j} \arrow[r, "\wt{h}"]\arrow[u, hook] & 
\displaystyle{\bigoplus_{0\leq i \leq 3} H^0(I_B\MO_{\P^3}(2)) e_i}. \arrow[u, hook]
\end{tikzcd}
\]

Take an element  
$\xi \in 
H^0(\P^3, I_B \otimes \Omega^2_{\P^3}(3))$. 
It suffices to show $\xi =0$. 
We can write 
\[
\xi = 
a \zeta_{012}+b\zeta_{013}+c \zeta_{023}+d \zeta_{123} \in 
k \zeta_{012} \oplus k \zeta_{013}\oplus k  \zeta_{023} \oplus k  \zeta_{123}
=
H^0(\P^3, \Omega^2_{\P^3}(3)) 
\]
for some $a, b, c, d \in k$. 
We have 
\[
\xi = a \zeta_{012}+b\zeta_{013}+c \zeta_{023}+d \zeta_{123} = 
(ax_0 + dx_3)e_1 \wedge e_2+ \cdots 
\in \displaystyle{\bigoplus_{0\leq i < j \leq 3} H^0(\MO_{\P^3}(1)) e_i \wedge e_j}, 
\]
i.e., $ax_0 + dx_3$ is the coefficient of $e_1 \wedge e_2$. 
By the above commutative diagram, we get 
\[
ax_0 + dx_3 \in H^0(I_B\MO_{\P^3}(1)) =H^0(\P^3, I_B \otimes \MO_{\P^3}(1)). 
\]
For the coordinate hyperplane $H_i := \{ x_i =0\}$, 
we have $H_0 \simeq \P^2$, $H_3 \simeq \P^2$, and $H_0 \cap H_3 \simeq \P^1$. 
Note that we may assume that $H_0, H_1, H_2, H_3$ are general hyperplanes by applying a suitable coordinate change in advance. 
In particular, we get  $B \cap H_i \neq \emptyset$ and $B \cap H_i \cap H_j= \emptyset$ for $i \neq j$. 
Fix a closed point $P_i \in B \cap H_i$ for each $i$. 
Since $\xi$ is vanishing at $P_3$ and $P_3 \in H_3 \setminus H_0$, 
we get $a=0$, 
because the coefficient of $e_1 \wedge e_2$ is $ax_0 + dx_3$. 
Similarly, it holds that $a=b=c=d=0$, i.e., $\xi =0$. 
Therefore, we get 
$H^0(I_B \otimes \Omega^2_{\P^3}(3)) = H^0(I_B \Omega^2_{\P^3}(3)) 
= 0$, as required. 
\end{proof}

\begin{thm}\label{t V5 Omega^2(1)}
Let $V$ be a  smooth Fano threefold with $r_V =2$ and $(-K_V/2)^3 = 5$. 
Let $D_V$ be an ample Cartier divisor $D_V$ satisfying $-K_V \sim 2D_V$. 
Then 
\[
H^0(V, \Omega^2_V \otimes \MO_V(D_V))=0.
\]
\end{thm}

\begin{proof}
We use Notation \ref{n V5 2ray}. 
By Lemma \ref{l V5 2ray easy}(1) and Lemma \ref{l V5 to P3}(3), 
it is enough to show $H^0(\P^3, I_B \otimes \Omega^2_{\P^3}(3)) =0$, 
which follows from Proposition \ref{p V5 Omega^2(1)}. 
\end{proof}

\begin{thm}\label{t g=6 2form}
Let $X \subset \P^7$ be an anti-canonically embedded Fano threefold 
of genus $6$ with $\Pic X= \Z K_X$. 
    Then $H^0(X, \Omega_X^2)=0$. 
\end{thm}

\begin{proof}
By Theorem \ref{t g=6 to Gr(2, 5)}, one of the following holds. 
\begin{enumerate}
\item[(A)] 
$X \simeq \Gr(2, 5) \cap H_1 \cap H_2  \cap Q$, 
    where each $H_i$ is a hyperplane of $\P^9$ and $Q$ is a quadric hypersurface of $\P^9$.
\item[(B)] 
There exists a finite surjective morphism $f: X \to V$ of degree two, 
where $V$ denotes the smooth Fano threefold with $r_{V}=2$ and $(-K_{V}/2)^3 = 5$. 
\end{enumerate}
If (A) holds, then the assertion  holds by Lemma \ref{lem:ANV for CI}. 
Hence we may assume (B). 
For the ample Cartier divisor $L$ on $V$ satisfying $-K_V \sim 2L$, 
we have an exact sequence 
$0 \to \MO_V \to f_*\MO_X \to \MO_V(-L) \to 0$. 
Taking the tensor product with $\Omega_V^1$, 
we get 
\[
0 \to \Omega_V^1 \to f_*f^*\Omega_V^1 \to \Omega_V^1(-L) \to 0. 
\]
By $H^3(V, \Omega_V^1) = H^3(V, \Omega_V^1(-L))=0$ (Theorem \ref{t V5 Omega^2(1)}), we get $H^3(X, f^*\Omega_V^1)=0$. 
Recall that the sequence 
\[
0 \to f^*\Omega_V^1 \to \Omega_X^1 \to \Omega_{X/V}^1 \to 0. 
\]
is exact (Step \ref{s3 g=6 ANV} in the proof of Proposition \ref{p g=6 ANV}). 
Since $f$ is generically \'etale (Theorem \ref{t g=6 separable}), we have $\dim \Supp \Omega_{X/V}^1 \leq 2$. 
Then the induced map 
\[
H^3(X,  f^*\Omega_V^1) \to H^3(X, \Omega_X^1)
\]
is surjective. 
Therefore, $H^3(X, f^*\Omega_V^1)=0$ implies $H^3(X, \Omega_X^1)=0$. 
We then get $H^0(X, \Omega_X^2)=0$ by Serre duality. 
\end{proof}

\section{Proofs of main theorems}\label{s main proofs}

In this section, we prove the main theorems in this paper, i.e., 
Theorem \ref{Introthm:W(k)-lift}, 
\ref{Introthm:ANV}, \ref{Introthm:E_1-degeneration}, 
and \ref{Introthm:hoge number} 
for the case when $\Pic X = \Z K_X$ (i.e., $\rho(X)=r_X=1$) 
and $|-K_X|$ is very ample. 
In this case, $X \subset \P^{g+1}$ 
is an anti-canonically embedded Fano threefold 
for  the closed embedding  $X \subset \P^{g+1}$  induced by $|-K_X|$. 

\begin{thm}[\textup{a special case of Theorem \ref{Introthm:W(k)-lift}}]\label{mainthm:lift}
Let $X$ be a smooth Fano threefold. 
Assume that $\Pic X = \Z K_X$ and $|-K_X|$ is very ample.  
    Then $X$ lifts to $W(k)$.
\end{thm}
\begin{proof}
Let $X \subset \P^{g+1}$ be the closed embedding induced by $|-K_X|$. 
Then the assertion holds by the following results: 
\begin{itemize}
\item $g \geq 7$: Theorem \ref{t-SB-lifting2}. 
\item $g=6$: Theorems \ref{t genus 6 lift}. 
\item $g \leq 5$: \cite[Proposition 2.8 and  the first paragraph of the proof of 
 Lemma 9.29]{FanoII}. 
\end{itemize}    
\end{proof}




\begin{thm}[\textup{a special case of Theorem \ref{Introthm:hoge number}}]\label{mainthm:hoge number}
Let $X$ be a smooth Fano threefold. 
Assume that $\Pic X = \Z K_X$ and $|-K_X|$ is very ample.  
Take a lift $\mathcal{X}$ of $X$ over $W(k)$, whose existence is ensured by Theorem \ref{Introthm:W(k)-lift}. 
Let $X_{\overline{K}}$ be the geometric generic fibre of $\mathcal{X}$ over $W(k)$. 
Then 
\[
h^j(X,\Omega^i_X)=h^j(X_{\overline{K}},\Omega^i_{X_{\overline{K}}})
\]
for all $i, j \geq 0$.
\end{thm}

\begin{proof}
Let $X \subset \P^{g+1}$ be the closed embedding induced by $|-K_X|$. 
    By \cite[Theorem 3.5 and Corollary 3.7]{Kaw2} (cf. \cite[Theorem 2.4]{FanoI}),
    we have 
    \[
    h^{i,j}(X) := h^j(X, \Omega_X^i) =0\quad\text{ for\quad all}\quad(i,j)\in\{(0,1),(0,2),(0,3),(1,0) 
    \}. 
    \]
    By Serre duality: $h^{i, j}(X) = h^{3-i, 3-j}(X)$, 
    it suffices to determine $h^{2,0}(X) ,h^{1,1}(X)$, and $h^{1,2}(X)$.
    Note that we have $\chi(X,\Omega^1_X)=\chi(X,\Omega^1_{X_{\overline{K}}})$ by the flatness of $\mathcal{X}$ over $W(k)$, where $X_{\overline{K}}$ denotes the geometric generic fibre of $\mathcal X \to \Spec W(k)$. 
    Therefore, if $h^{2,0}(X)=h^{2,0}(X_{\overline{K}})$ and $h^{1,1}(X)=h^{1,1}(X_{\overline{K}})$, then
    we obtain 
    \begin{align*}
        h^{1,2}(X)&=\chi(X,\Omega^1_X)-h^{1,0}(X)+h^{1,1}(X)+h^{1,3}(X)\\
        &=\chi(X_{\overline{K}},\Omega^1_{X_{\overline{K}}})-h^{1,0}(X_{\overline{K}})+h^{1,1}(X_{\overline{K}})+h^{1,3}(X_{\overline{K}})\\
        &=h^{1,2}(X_{\overline{K}}).
    \end{align*}
    It is enough to show $h^{2,0}(X)=h^{2,0}(X_{\overline{K}})$ and $h^{1,1}(X)=h^{1,1}(X_{\overline{K}})$.
    If $g\in\{3,4,5\}$ (resp.~$g=6$), then these equalities hold 
    by Theorem \ref{thm:CI} (resp. Proposition \ref{p g=6 ANV} and Theorem \ref{t g=6 2form}).
    If $g \geq 7$, then $X$ is rational or birational to a smooth cubic threefold (Proposition \ref{p rat or birat to cubic}(2)), which implies $h^{2,0}(X)=0=h^{2,0}(X_{\overline{K}})$ and $h^{1,1}(X)=1=h^{1,1}(X_{\overline{K}})$ by Lemma \ref{l-sm-RET} and Theorem \ref{t-rat3-cohs}, respectively.
\end{proof}

\begin{thm}[\textup{a special case of Theorem \ref{Introthm:ANV}}]\label{mainthm:Akizuki-Nakano}
Let $X$ be a smooth Fano threefold. 
Assume that $\Pic X = \Z K_X$ and $|-K_X|$ is very ample.  
Then $X$ satisfies Akizuki-Nakano vanishing, i.e., 
$H^j(X, \Omega_X^i \otimes \MO_X(-A))=0$ 
if $i+j <3$ and $A$ is an ample Cartier divisor $A$ on $X$.
\end{thm}
\begin{proof}
Let $X \subset \P^{g+1}$ be the closed embedding induced by $|-K_X|$. Fix an ample Cartier divisor $A$ on $X$. 
{Since $\Pic X = \Z K_X$ and $|-K_X|$ is very ample, we get $H^0(X,\sO_X(A))>0$.} 
Thus $H^0(X,\Omega^i_X(-A)) \hookrightarrow H^0(X,\Omega^i_X)=0$ for all $i> 0$ by Theorem \ref{mainthm:hoge number}.
By \cite[Corollary 4.5]{FanoI}, we have $H^j(X,\sO_X(-A))=0$ for all $j>0$.
Hence it suffices to show that $H^1(X,\Omega^1_X(-A))=0$. 

If $g=6$ (resp.~$g \geq 7$), 
then we get the vanishing 
$H^1(X,\Omega^1_X(-A))=0$ 
by Proposition \ref{p g=6 ANV} (resp.~Theorem \ref{t-SB-lifting2}). 
Then we may assume that $g \in \{3, 4, 5\}$ \cite[Theorem 1.1]{FanoII}. 
If $g=3$ (resp.~$4$, resp.~$5$), then $X$ is a hypersurface of $\P^4$ (resp.~complete intersection of $\P^5$, resp.~complete intersection of $\P^6$) by \cite[Proposition 2.8]{FanoII}, and hence 
the vanishing 
$H^1(X,\Omega^1_X(-A))=0$ 
holds by Lemma \ref{lem:ANV for CI}. 
\end{proof}

\begin{thm}[\textup{a special case of Theorem \ref{Introthm:E_1-degeneration}}]\label{mainthm:E_1-deg}
Let $X$ be a smooth Fano threefold. 
Assume that $\Pic X = \Z K_X$ and $|-K_X|$ is very ample.  
Then the Hodge to de Rham spectral sequence degenerates at $E_1$. 
Moreover, $H^i_{\cris}(X/W)$ is torsion-free for all $i\geq 0$, where $W :=W(k)$.
\end{thm}
\begin{proof}
    The assertion holds by Theorem \ref{mainthm:lift}, Theorem \ref{mainthm:hoge number}, and Proposition \ref{prop:hodge number to degeneration and torsion-freeness}.
\end{proof}


\begin{prop}\label{prop:hodge number to degeneration and torsion-freeness}
    Let $X$ be a smooth projective variety.
    Assume that there exists a lift $\mathcal X$ of $X$ over $W:=W(k)$ and 
    let $X_{\overline{K}}$ be its geometric generic fibre.
    Suppose that $h^{i,j}(X)=h^{i,j}(X_{\overline{K}})$ for all $i,j\geq 0$.
    Then the following hold.
    \begin{enumerate}
        \item The Hodge to de Rham spectral sequence degenerates at $E_1$.
        \item $H^{i}_{\cris}(X/W)$ is torsion-free for all $i\geq 0$.
    \end{enumerate}
\end{prop}
\begin{proof}
    For every integer $i \geq 0$, the following holds: 
    \begin{align*}
        h^{i,0}(X)+h^{i-1,1}(X)+\cdots+h^{0,i}(X)&\geq \dim_k
        H^i(X, \Omega^{\bullet}_X)\\
        &\geq 
    \mathrm{rank}_{W}H^{i}_{\cris}(X/W)\\
    &= b_i(X)\\ 
    &\overset{(\star)}{=}b_i(X_{\overline{K}})\\
    &=h^{i,0}(X_{\overline{K}})+h^{i-1,1}(X_{\overline{K}})+\cdots+h^{0,i}(X_{\overline{K}})
    \end{align*}
     Here $(\star)$ follows from \cite{KM74} when $k$ is an algebraic closure of a finite field, and the general case is reduced to this case by the smooth proper base change theorem for \'etale and Crystalinne cohomologies \cite[Theorem II]{PTZ23}. 
    Note that the first equality is equivalent to the degeneration at $E_1$ of the Hodge to de Rham spectral sequence,
    and the second equality is equivalent to the torsion-freeness of $H^{i}_{\cris}(X/W)$ and $H^{i+1}_{\cris}(X/W)$
    by the exact sequence
    \[
0 \to H^i_{\Cris}(X/W) \otimes_W k 
\to H^i(X, \Omega_X^{\bullet}) \to {\rm Tor}_1^W(H^{i+1}_{\Cris}(X/W), k) \to 0. 
\]
 By assumption, the above inequalities have to be all equalities, and hence the assertions hold.
\end{proof}

\begin{thm}\label{t line}
Let $X$ be a smooth Fano threefold. 
Assume that $\Pic X = \Z K_X$ and $|-K_X|$ is very ample. 
Then there exists a line $L$ on $X$, i.e., 
a curve $L$ on $X$ satisfying $-K_X \cdot L = 1$.
\end{thm} 

\begin{proof}
By Theorem \ref{mainthm:lift}, 
there exists a lift $\mathcal X$ of $X$ over $W(k)$. 
Let $X_{\overline K}$ be the geometric generic fibre of $\mathcal X \to \Spec W(k)$. 
It is enough to show ($\star$) below 
(cf. \cite[the proof of Lemma 9.29]{FanoII}). 
\begin{enumerate}
\item[($\star$)] There exists a line on $X_{\overline K}$. 
\end{enumerate}
If $g \neq 7$, then ($\star$) follows from known results as follows: 
\begin{itemize}
\item $g \leq 5$: \cite[Corollary 1.4]{Can21}. 
\item $g \geq 8$: \cite[Theorem 0.2]{Tak89}. 
\item $g = 6$:  \cite[Lemma 2.6(a), Theorem 4.7(d)]{DK19}. 
\end{itemize}
Assume $g=7$. 
In this case, there exists a conic on $X_{\overline K}$ (cf. Theorem \ref{t conic g=7}). 
Then ($\star$)  holds by the same argument as in \cite[Corollary 4.4.10]{IP99}. 
\end{proof}

\begin{cor}\label{c rationality}
Let $X$ be a smooth Fano threefold. 
Assume that $\Pic X = \Z K_X$, $|-K_X|$ is very ample, and $g \in \{ 7, 9, 10, 12\}$, where $g$ is the integer defined by $2g-2 =(-K_X)^3$. 
Then $X$ is rational. 
\end{cor}

\begin{proof}
If $g \neq 7$, then the assertion follows from 
Proposition \ref{p rat or birat to cubic}(1). 
Assume $g=7$. 
Take a line $L$ on $X$, whose existence is guaranteed by Theorem \ref{t line}. 
By \cite[Proposition 5.2, Theorem 8.1, Corollary 8.3]{FanoII}, 
there exist a birational map $X \dashrightarrow Y^+$ and a del Pezzo fibration 
$\pi : Y^+ \to \P^1$ with $(-K_{Y^+})^2 \cdot \pi^{-1}(b) =5$, 
where $b$ is a closed point of $\P^1$. 
Then $Y^+$ is rational \cite[Theorem 1.4]{BT24}, and hence so is $X$. 
\end{proof}

\section{Appendix: Veronese surfaces}\label{s Veronese}

\begin{nota}\label{n Veronese}
We work over an algebraically closed field $k$ of arbitrary characteristic. 
Let $S \subset \P^5_k = \Proj\,k[x, y, z, s, t, u]$ be the Veronese surface, i.e., 
$S$ is the image of the Veronese embedding 
\[
\nu : \P^2_k \hookrightarrow \P^5_k, \qquad [X:Y:Z] \mapsto [X^2:Y^2:Z^2:YZ:ZX:XY]. 
\]
Set $V$ to be the secant variety of $S \subset \P^5_k$, which is given by 
\[
V := \overline{V^{\circ}}
\qquad\text{and}\qquad 
V^{\circ} := \bigcup_{\substack{P, Q \in S\\ P \neq Q}} \overline{PQ}, 
\]
where 
$\overline{V^{\circ}}$ is the closure of $V^{\circ}$ in $\P^5_k$ and 
$\overline{PQ}$ denotes the line passing through $P$ and $Q$. 
We equip $V$ with the reduced scheme structure.  
\end{nota}

\begin{rem}\label{r Veronese sym matrix}
Let $\Sigma_3(k)$  be the set consisting of all the $3 \times 3$ symmetric matrices. 
Set $\overline{\Sigma}_3(k) := (\Sigma_3(k) \setminus \{O\})/k^{\times}$, 
where $O$ denotes the zero matrix. 
We then have the following bijective corresponding:  
\begin{align*}
\P^5_k(k) & \xleftrightarrow{\text{bije.}} \overline{\Sigma}_3(k)\\ 
[a:b:c:d:e:f] &\leftrightarrow
\begin{bmatrix}
a & f & e\\
f & b & d\\
e & d & c
\end{bmatrix}\\
S &\xleftrightarrow{(\star)} \{ [A] \in \overline{\Sigma}_3(k) \,|\, A \in \Sigma_3(k) \setminus \{O\}, \rank A =1\}\\
V &\xleftrightarrow{(\star\star)} \{ [A] \in \overline{\Sigma}_3(k) \,|\, A \in \Sigma_3(k) \setminus \{O\}, \rank A \leq 2 \}
\end{align*}  
where 
$[A] \in \overline{\Sigma}_3(k)$ denotes the image of $A \in \Sigma_3(k) \setminus \{O\}$, and 
\[
\begin{bmatrix}
a & f & e\\
f & b & d\\
e & d & c
\end{bmatrix}
:= 
\left[ 
\begin{pmatrix}
a & f & e\\
f & b & d\\
e & d & c
\end{pmatrix}\right]. 
\]
The correspondence $(\star)$ holds, because  the defining homogeneous ideal of $S$ is given by 
\[
(xy - u^2, yz - s^2, zx - t^2, xs-tu, yt-us, zu-st), 
\]
whose generators are nothing but the determinants of all the $2 \times 2$ minor matrices. 
The correspondence $(\star\star)$ follows from  
\[
\det \begin{pmatrix}
x & u & t\\
u & y & s\\
t & s & z
\end{pmatrix} = xyz +2stu -xs^2 -yt^2 -zu^2
\]
and Lemma \ref{l Veronese secant} below. 
\end{rem}


\begin{lem}\label{l Veronese secant}
We use Notation \ref{n Veronese}. 
Then the following hold. 
\begin{enumerate}
\item 
$V = \{ xyz +2stu -xs^2 -yt^2 -zu^2=0\}$. 
Moreover, $V$ is a $4$-dimensional variety. 
\item $V = \bigcup_{P \in S} T_PS$, 
where $T_PS$ denotes the embedded tangent plane of $S \subset \P^5$. 
\end{enumerate}
\end{lem}


\begin{proof}
Recall that $V$ is a variety, because $V$ is the closure of the image of 
an irreducible set $ ((X \times X) \setminus \Delta_X) \times \P^1_k$. 
By \cite[Page 15 and Theorem 1.4]{Zak93}, 
(1) implies (2) and we get $\dim V \geq 4$. 
Hence it suffices to show (1) under assuming $\dim V \geq 4$ and $V$ is a variety. 
Then the assertion (1) follows from the same argument as in \cite[page 145]{Har95}. 
\qedhere


\end{proof}




\begin{prop}\label{p Veronese cases}
We use Notation \ref{n Veronese}. 
Fix a closed point $R \in \P^5$ and 
let  $g : S \dashrightarrow T$ be the projection from $R \in \P^5$ 
to its image $T \subset \P^4_k$. 
Then the following hold. 
\begin{enumerate}
\item $R \not\in S$ if and only if $g$ is a morphism 
(i.e., the largest open subset $S'$ of $S$ on which $g$ is defined 
 is  equal to $S$). 
\item 
If $g : S \to T$ is a morphism, then one of the following holds. 
\begin{enumerate}
\item $R \not\in V$ and 
$g : S \to T$ is an isomorphism. In particular, $T \simeq \P^2_k$. 
\item 
$R \in V \setminus S$, $g: S \to T$ is birational, and 
$\dim {\rm Sing}\,T = 1$, where ${\rm Sing}\,T$ denotes the non-smooth locus of $T$. 
\end{enumerate}
\end{enumerate}
    \end{prop}

\begin{proof}
Let us show (1). 
It is clear that $g$ is a morphism if $R \not\in S$. 
Assume $R \in S$. 
Recall that the projection $\P^5 \dashrightarrow \P^4$ from $R$ 
is nothing but the composition of the blowup $\mu : \Bl_R \P^5 \to \P^5$ at $R$ and the induced $\P^1$-bundle 
$\Bl_R \P^5 \to \P^4$ (whose fibres are corresponding to the lines on $\P^5$ passing through $R$). 
In particular, $\pi|_{\Ex(\mu)} : \Ex(\mu) \xrightarrow{\simeq} \P^4$. 
If $R \in S$, then the induced morphism $\Bl_R S \to T$ 
would  not contract the $(-1)$-curve $\Bl_R S \cap \Ex(\mu)$, 
which implies that $g: S \dashrightarrow T$ can not be a morphism. 
Thus (1) holds.

Let us show (2). Assume that $g : S \to T$ is a morphism. 
By definition, $T$ is a projective surface.

\begin{claim*}
$g$ is birational. 
\end{claim*}

\begin{proof}[Proof of Claim]
Since $H^0(\P^5, \MO_{\P^5}(1)) \xrightarrow{\simeq} H^0(S, \MO_{\P^5}(1)|_S)$ 
is an isomorphism, 
the induced morphism $g: S \to T$ is defined by a base point free linear system corresponding to a $5$-dimensional $k$-vector subspace of $H^0(S, \MO_{\P^5}(1)|_S)$. 
We then get $h^0(T, \MO_T(1)) = 5$, where $\MO_T(1) := \MO_{\P^4}(1)|_T$. 
For the $\Delta$-genus 
$\Delta(T, \MO_T(1))$, it holds that 
\begin{align*}
0 
&\leq \Delta(T, \MO_T(1))\\
&= \dim T + \MO_T(1)^2 -h^0(T, \MO_T)\\
&= 2 + \frac{(\MO_{\P^5}(1)|_S)^2}{\deg g} -5\\
&= \frac{4}{\deg g} -3, 
\end{align*}
which implies $\deg g=1$. This complete the proof of Claim. 
\end{proof}

Assume $R \not\in V$. Then $g$ is an isomorphism, i.e., (a) holds. 
Assume $R \in V \setminus S$. 
Let us show that (b) holds. 
By Claim, it suffices to show that $\dim {\rm Sing}\,T=1$. 
We have $R \in V = \bigcup_{P \in S} T_PS$ (Lemma \ref{l Veronese secant}), i.e., 
$R \in T_PS$ for some closed point $P \in S$. 
By $P \neq R$, we have the line  $M := \overline{PR}$ passing through $P$ and 
$R$. 
Then the zero-dimensional closed subscheme $M \cap S$ 
satisfies is of degree two
(indeed, $\deg (M \cap S) \geq 2$ holds by the fact that $M = \overline{PR} \subset T_PS$. The other inequality $\deg (M \cap S) \leq 2$ follows from the fact that $S$ is an intersection of quadrics (Remark \ref{r Veronese sym matrix})). 
Via the Veronese embedding $\nu : \P^2 \xrightarrow{\simeq} S$, 
there exists a unique line $L$ on $\P^2$ satisfying $\nu^{-1}(M \cap S) \subset L$. 
Then 
we get $M \cap S = \nu(\nu^{-1}(M \cap S)) \subset 
\nu(L) =:C$. Here $C\subset \P^5$ is a conic. 
Since the line $M = \overline{PR}$ is the smallest linear variety in $\P^5_k$ containing $M \cap S$, 
we have $M \subset \Pi_C$ for the plane $\Pi_C$ containing the conic $C$.  
Note that there exist infinitely many lines $\{ L_{\lambda}\}_{\lambda \in \Lambda}$ on the plane $\Pi_C (\simeq \P^2)$ 
passing through $R$. 
Since each $L_{\lambda}$ corresponds to a fibre of the $\P^1$-bundle 
$\pi : \Bl_R\,\P^5 \to \P^4$, 
$T$ is singular at $\pi(L_{\lambda})$, because 
$\deg(L_{\lambda} \cap S) \geq \deg(L_{\lambda} \cap C) \geq 2$. 
Hence $T$ has infinitely many singular points, i.e., $\dim ({\rm Sing} T) =1$. 
This completes the proof of (2). 
\end{proof}

The following result is claimed without a proof in \cite[Exercise B.5]{Pro}. 
We give a proof for the sake of completeness.

\begin{lem}\label{l Veronese no quad0}
We use Notation \ref{n Veronese}. 
Fix a closed point $R \in \P^5_k$ and 
let  $g : S \dashrightarrow T$ be the projection from $R \in \P^5_k$ 
to its image $T \subset \P^4_k$. 
Assume that  $k$ is of characteristic $\neq 2$ and $g$ is a morphism. 
Then there exists no smooth quadric hypersurface containing $T$. 
\end{lem}

\begin{proof}   

\setcounter{step}{0}
\begin{step}\label{s1 Veronese no quad0}
In order to prove Lemma \ref{l Veronese no quad0}, 
we may assume that $R = [1:1:0:0:0:0] \in \P^5$. 
\end{step}

\begin{proof}[Proof of Step \ref{s1 Veronese no quad0}]
In what follows, we use the identification 
\[
\P^5_k = \P^5_k(k) = \overline{\Sigma}_3(k) := (\Sigma_3(k) \setminus \{ O \})/k^{\times}
\]
given in Remark \ref{r Veronese sym matrix}, 
where $\Sigma_3(k)$ is the set consisting of all the $3 \times 3$ symmetrix matrices. 
We have the left $\GL_3(k)$-action: 
\[
\GL_3(k) \times \overline{\Sigma}_3(k) \to \overline{\Sigma}_3(k), \qquad 
(P, [A]) \mapsto [PAP^t], 
\]
where $[A]$ denotes the equivalent class to which $A \in \Sigma_3(k) \setminus \{ O \}$ belongs. 
By linear algebra, 
there are exactly three $\GL_3(k)$-orbits of $\overline{\Sigma}_3(k)$ 
corresponding to $\rank =1, 2, 3$. 
By Remark \ref{r Veronese sym matrix}, 
$S$ (resp. $V \setminus S$) coincides with the orbit consisting of the matrices of rank $1$ 
(resp. $2$). 
Since the $\GL_3(k)$-action on $\P^5_k$ is given by linear transforms, 
we may replace $R$ by an arbitrary closed point of $\P^5_k$ corresponding to $[A] \in \overline{\Sigma}_3(k)$ satisfying $\rank A = 2$; 
$[1:1:0:0:0:0] \in \P^5_k$ is such a point, because the corresponding matrix is given by 
$\begin{pmatrix}
1 & 0 & 0\\
0 & 1 & 0\\
0 & 0 & 0
\end{pmatrix}$ (Remark \ref{r Veronese sym matrix}). 
This completes the proof of Step \ref{s1 Veronese no quad0}. 
\qedhere

\end{proof}

\begin{step}\label{s2 Veronese no quad0}
The assertion of Lemma \ref{l Veronese no quad0} 
holds if $R = [1:1:0:0:0:0] \in \P^5$. 
\end{step}

\begin{proof}[Proof of Step \ref{s2 Veronese no quad0}]
Assume $R = [1:1:0:0:0:0] \in \P^5_k$. 
Take an automorphism: 
\[
\sigma: \P^5 \xrightarrow{\simeq} \P^5, \qquad [x:y:z:s:t:u] \mapsto 
[x: x-y :z:s:t:u], 
\]
which satisfies $\sigma(R) = R' :=[1:0:0:0:0:0]$. 
Recall that we have 
\[
S = \Proj\,\frac{k[x, y, z, s, t, u]}{(xy - u^2, yz -s^2, zx -t^2, xs-tu, yt -us, zu-st)}. 
\]
Then 
\[
S' = \sigma(S) = \Proj\,\frac{k[x, y, z, s, t, u]}{(x(x-y) - u^2, (x-y)z -s^2, zx -t^2, xs-tu, (x-y)t -us, zu-st)}. 
\]
Let us determine the defining equations of 
the projection $T$ of $S'$ from $R' = [1:0:0:0:0:0]$ 
to $H := \{ x =0\}$. 
This projection satisfies   
\begin{align*}
f: \P^5 \setminus [1:0:0:0:0:0] &\to \P^4\\
[1 : b:c:\alpha:\beta:\gamma]&\mapsto [b:c:\alpha:\beta:\gamma]. 
\end{align*}
The $k$-algebra homomorphism corresponding to $\wt{f}: D_+(x) \cap D_+(y) \to D_+(y)$ can be written as 
\begin{align*}
\varphi : k[Z, S, T, U] &\to k[y, y^{-1}, z, s, t, u]\\
Z &\mapsto z/y\\
S &\mapsto s/y\\
T &\mapsto t/y\\
U &\mapsto u/y. 
\end{align*}

The scheme-theoretic image of $\wt{f}$ is corresponding to the kernel of 
\[
\psi : k[Z, S, T, U] \xrightarrow{\varphi} k[y, z, s, t, u] 
\]
\[
\to 
\frac{k[y, z, s, t, u]}{(1-y - u^2, (1-y)z -s^2, z -t^2, s-tu, (1-y)t -us, zu-st)}. 
\]
After substituting $y = -u^2 +1, z =t^2, s = tu$, 
$\psi$ can be written as 
\begin{align*}
\psi : k[Z, S, T, U] &\to k[t, u, (-u^2+1)^{-1}]\\
Z &\mapsto t^2/(-u^2+1)\\
S &\mapsto tu/(-u^2+1)\\
T &\mapsto t/(-u^2+1)\\
U &\mapsto u/(-u^2+1). 
\end{align*}
Then 
\[
S^2 -TU, ST -UZ \in \Ker(\psi). 
\]
Since $(S^2-TU, ST-UZ)$ is a prime ideal (Lemma \ref{l Veronese projection prime}), 
we get $\Ker(\psi) = (S^2-TU, ST-UZ)$ by counting dimension.

Therefore, the image $T$ of the projection of $S$ is the closure in $\P^4$ 
of 
\[
k[Z, S, T, U]/(S^2-TU, ST-UZ) \subset D_+(y) = \A^4, 
\]
and hence 
\[
T  = \Proj\,\frac{k[y, z, s, t, u]}{(s^2 -tu, st-uz)} \subset 
\Proj\,k[y, z, s, t, u] =\P^4. 
\]
Then it suffices to show that  
\[
\{\alpha (s^2 -tu) + \beta ( st -uz) =0\} \subset \P^4 
\]
is not smooth for every $(\alpha, \beta) \in k^2 \setminus \{(0, 0)\}$, 
which follows from the fact that $[1:0:0:0:0]$ is always a singular point. 
This completes the proof of Step \ref{s2 Veronese no quad0}. 
\end{proof}
Step \ref{s1 Veronese no quad0} and 
Step \ref{s2 Veronese no quad0} complete the proof of Lemma \ref{l Veronese no quad0}. 
\end{proof}

\begin{lem}\label{l Veronese projection prime}
Let $k$ be an algebraically closed field. 
Then the ideal $(S^2 -TU, ST-UZ)$ is a prime ideal of $k[Z, S, T, U]$. 
\end{lem}

\begin{proof}
Set $I := (S^2 -TU, ST-UZ)$. 
For $R := k[Z, U]$, we have an $R$-module isomorphism: 
\begin{equation}\label{e1 Veronese projection prime}
k[Z, S, T, U]/I = R S \oplus R[T] \simeq R \oplus R[T]. 
\end{equation}
In particular, the induced $k$-algebra homomorphism $R \to k[Z, S, T, U]/I$ is flat. 
Then, for  the function field $K :=K(R)$ of $R$, 
we get an injective $k$-algebra homomorphism 
\[
k[Z, S, T, U]/I \hookrightarrow k[Z, S, T, U]/I \otimes_R K. 
\]
Thus it is enough to show that 
$k[Z, S, T, U]/I \otimes_R K$ is an integral domain. 

Note that 
\[
k[Z, S, T, U]/I \otimes_R R[U^{-1}] \simeq k[S, U, U^{-1}]. 
\]
Then 
\[
k[Z, S, T, U]/I \otimes_R K \simeq W^{-1}(k[S, U, U^{-1}]) 
\]
for some multiplicatively closed subset $W$ of $k[S, U, U^{-1}]$. 
By $k[Z, S, T, U]/I \otimes_R K \neq 0$ (\ref{e1 Veronese projection prime}), we get $W^{-1}(k[S, U, U^{-1}]) \neq 0$, and hence 
$W^{-1}(k[S, U, U^{-1}])$ is an integral domain. 
\end{proof}

\begin{thm}\label{t Veronese no quad}
We use Notation \ref{n Veronese}. 
Fix a closed point $R \in \P^5_k$ and 
let  $g : S \dashrightarrow T$ be the projection from $R \in \P^5_k$ 
to its image $T \subset \P^4_k$. 
Assume that $g$ is a morphism, i.e., $R \not\in S$ (Proposition \ref{p Veronese cases}). 
Then $T$ is not contained any smooth quadric hypersurface in $\P^4_k$. 
\end{thm}

\begin{proof}

\setcounter{step}{0}

\begin{step}\label{s1 Veronese no quad}
The assertion of Theorem \ref{t Veronese no quad} holds 
if $R \not\in V$. 
\end{step}

\begin{proof}[Proof of Step \ref{s1 Veronese no quad}]
Assume $R \not\in V$. 
Suppose that $T \subset U$ for some  smooth quadric hypersurface $U \subset \P^4$. 
We have that $\MO_U(T) \simeq \MO_U(d)$ for some $d>0$, where 
$\MO_U(d) := \MO_{\P^4}(d)|_U$. 
By $T \simeq \P^2$ and the adjunction formula, we get the following contradiction: 
\[
9 = K_{\P^2}^2 =
(K_U + T)^2 \cdot T = 
\]
\[
\MO_U(-3+d)^2 \cdot \MO_U(d) = 
\MO_{\P^4}(-3+d)^2 \cdot \MO_{\P^4}(d) \cdot \MO_{\P^4}(2) = 
2d(-3+d)^2 \in 2\Z. 
\]     
This completes the proof of Step \ref{s1 Veronese no quad}. 
\end{proof}

\begin{step}\label{s2 Veronese no quad}
The assertion of Theorem \ref{t Veronese no quad} holds 
if $R \in V$. 
\end{step}


\begin{proof}[Proof of Step \ref{s2 Veronese no quad}]
Assume $R \in V \setminus S$. 
The assertion of Theorem \ref{t Veronese no quad} has been already established for the case when $k$ is of characteristic zero 
(Lemma \ref{l Veronese no quad0}). 
In what follows, we assume that $k$ is of characteristic $p>0$ 
and reduce the problem to the case of characteristic zero by constructing suitable $W(k)$-lifts. 
Set $Q := \Frac\,W(k)$ and let $\overline{Q}$ be its algebraic closure. 

We have  $S \subset V \subset \P^5_k = \Proj\,k[x,y,z,s,t,u]$. 
Recall that the Veronese surface $S$ and its secant variety $V$ are given 
as follows (Remark \ref{r Veronese sym matrix}): 
\[
S = \Proj\,\frac{k[x,y,z,s,t,u]}{(xy-u^2, yz-s^2, zx -t^2, xs-tu, yt-us, zu-st)} 
\subset  \P^5_k, 
\]
\[
V = \Proj\,\frac{k[x,y,z,s,t,u]}{(xyz+2stu-xs^2-yt^2-zu^2)} 
\subset  \P^5_k. 
\]
Since these descriptions work in arbitrary characteristic, 
we get a lift $\wt{S} \subset \wt{V} \subset \P^5_{W(k)}$ 
of $S \subset V \subset \P^5_{k}$: 
\[
\wt{S} = \Proj\,\frac{W(k)[x,y,z,s,t,u]}{(xy-u^2, yz-s^2, zx -t^2, xs-tu, yt-us, zu-st)} 
\subset  \P^5_{W(k)}, 
\]
\[
\wt{V} = \Proj\,\frac{W(k)[x,y,z,s,t,u]}{(xyz+2stu-xs^2-yt^2-zu^2)} 
\subset  \P^5_{W(k)}. 
\]
The closed point $R \in V \setminus S$ admits a lift $\wt{R} \subset \wt{V}$ 
\cite[Th\'eor\`eme 18.5.17]{EGAIV4}, 
which is a section of the induced morphism 
$\wt{V} \to \Spec\,W(k)$ \cite[Proposition 2.8.14]{Fu15}. 
Then the base change $R_{\ol{Q}} :=R \times_{W(k)} \overline{Q}$ 
is a closed point on $V_{\ol{Q}}$. 

\begin{claim*}
There exist a closed subscheme $\wt{T} \subset \P^4_{W(k)}$ 
and a $W(k)$-morphism $\wt{g} : \wt{S} \to \wt{T}$ such that 
\begin{enumerate}
\item[(i)] $\wt{T} \subset \P^4_{W(k)}$  is a $W(k)$-lift of $T \subset \P^4_k$, 
\item[(ii)] $\wt{g} : \wt{S} \to \wt{T}$ is a $W(k)$-lift of 
$g : S \to T$, and 
\item[(iii)] 
the base change 
$g_{\ol Q} : S_{\overline Q} \to T_{\overline Q}$ of $\wt{g}$ by $(-) \times_{W(k)} \ol{Q}$ is the projection 
from the point $R_{\ol{Q}}$ on the secant variety $V_{\ol Q}$ of $S_{\overline Q} 
\subset \P^5_{\ol{Q}}$. 
\end{enumerate}    
\end{claim*}

\begin{proof}[Proof of Claim]
For the blowup $\wt{\mu} : \Bl_{\wt{R}} \P^5_{W(k)} \to \P^5_{W(k)}$ 
of $\P^5_{W(k)}$ along $\wt{R}$, 
the $\P^1$-bundle structure $\wt{\pi} : \Bl_{\wt{R}} \P^5_{W(k)} \to \P^4_{W(k)}$
is induced by the complete linear system $|\wt{\mu}^*\MO_{\P^5_{W(k)}}(1) - \wt{E}|$ for $\wt{E} := \Ex(\wt{\mu})$.  
We then define $\wt{T}$ as the image of the composite morphism 
\[
(\wt{S} \simeq) \wt{\mu}^{-1}(\wt{S}) \hookrightarrow 
\Bl_{\wt{R}} \P^5_{W(k)} \xrightarrow{\wt{\pi}} \P^4_{W(k)}, 
\]
which is given by the complete linear system $|\wt{\mu}^*\MO_{\P^5_{W(k)}}(1) - \wt{E}|$. 
In order to show (i)-(iii), it is enough to prove that 
the restriction map 
\[
\rho : H^0(\Bl_{\wt{R}} \P^5_{W(k)}, 
\wt{\mu}^*\MO_{\P^5_{W(k)}}(1) - \wt{E}) \to 
H^0(\Bl_R \P^5_{k}, 
\mu^*\MO_{\P^5_{k}}(1) - E)
\]
is surjective. 
By the exact sequence 
\[
H^0(\Bl_{\wt{R}} \P^5_{W(k)}, 
\wt{\mu}^*\MO_{\P^5_{W(k)}}(1) - \wt{E}) \xrightarrow{\rho}
H^0(\Bl_R \P^5_{k}, 
\mu^*\MO_{\P^5_{k}}(1) - E)
\]
\[
\to 
H^1(\Bl_{\wt{R}} \P^5_{W(k)}, 
\wt{\mu}^*\MO_{\P^5_{W(k)}}(1) - \wt{E})
\xrightarrow{\times p}
H^1(\Bl_{\wt{R}} \P^5_{W(k)}, 
\wt{\mu}^*\MO_{\P^5_{W(k)}}(1) - \wt{E})
\]
\[
\to H^1(\Bl_R \P^5_{k}, 
\mu^*\MO_{\P^5_{k}}(1) - E)
\]
and Nakayama's lemma, the problem is reduced to showing 
\begin{equation}\label{e1 Veronese no quad}
H^1(\Bl_R \P^5_{k}, 
\mu^*\MO_{\P^5_{k}}(1) - E) =0. 
\end{equation}
Since $\Bl_R \P^5_{k}$ is toric, 
(\ref{e1 Veronese no quad}) follows from the fact that $\mu^*\MO_{\P^5_{k}}(1) - E$ is nef 
and \cite[Corollary 1.7]{Fuj07}. 
This completes the proof of Claim. 
\end{proof}


Suppose that there exists a smooth quadric hypersurface $U \subset \P^4_k$ 
containing $T$. 
Then Lemma \ref{l containing quad lift} enables us to find  
a $W(k)$-lift $\wt{T} \subset \wt{U}$ of $T \subset U$. 
We then get $T_{\ol{Q}} \subset U_{\ol{Q}}$, and hence $T_{\ol{Q}}$ is contained 
in a smooth quadric hypersurface $U_{\ol{Q}} \subset \P^4_{\ol{Q}}$. 
This contradicts Claim and Lemma \ref{l Veronese no quad0}. 
This completes the proof of Step \ref{s2 Veronese no quad}. 
\qedhere

\end{proof}
Step \ref{s1 Veronese no quad} and 
Step \ref{s2 Veronese no quad} complete the proof of Theorem \ref{t Veronese no quad}. 
\end{proof}

\begin{lem}\label{l containing quad lift}
Let $T = U \cap U' \subset \P^4_k$ be a complete intersection of 
hypersurfaces of degree $a$ and $b$. 
Take an arbitrary $W(k)$-lift $\wt{T} \subset \P^4_{W(k)}$ of $T \subset \P^4_k$. 
Then the restriction map 
\[
H^0(\P^4_{W(k)}, I_{\wt{T}} \otimes \MO_{\P^4_{W(k)}}(n)) \to 
H^0(\P^4_{k}, I_{T} \otimes \MO_{\P^4_{k}}(n))
\]
is surjective for every $n \in \Z$. 
\end{lem}

\begin{proof}
We have the following commutative diagram: 
\[
\begin{tikzcd}
& 0 \arrow[d] & 0 \arrow[d] & 0 \arrow[d] & \\
0 \arrow[r] & I_{\wt T} \arrow[r, "\times p"] \arrow[d] & I_{\wt T}\arrow[r] \arrow[d] & I_{T} \arrow[r] \arrow[d] & 0 \\
0 \arrow[r] & \MO_{\P^4_{W(k)}} \arrow[r, "\times p"] \arrow[d] & \MO_{\P^4_{W(k)}} \arrow[r] \arrow[d] & \MO_{\P^4_k} \arrow[r] \arrow[d] & 0 \\
0 \arrow[r] & \MO_{\wt{T}} \arrow[r, "\times p"] \arrow[d] & \MO_{\wt{T}} \arrow[r] \arrow[d] & \MO_{T} \arrow[r] \arrow[d] & 0. \\
& 0 & 0 & 0 &
\end{tikzcd}
\]
The second and third horizontal sequences are exact, because 
$\P^4_{W(k)}$ and $\wt{T}$ are $W(k)$-lifts of $\P^4_k$ and $T$, respectively. 
Applying the snake lemma, we see that 
all the horizontal and vertical sequences are exact. 

Fix $n \in \Z$. 
In order to show 
$H^1(\P^4_{W(k)}, I_{\wt T} \otimes \MO_{\P^4_{W(k)}}(n))=0$, 
it is enough, by Nakayama's lemma, to prove $H^1(\P^4_k, I_T \otimes \MO_{\P^4_{k}}(n))=0$. 
The problem is reduced to the surjectivity of 
the restriction map 
\[
\rho : H^0(\P^4_k, \MO_{\P^4_{k}}(n)) \to H^0(T, \MO_{\P^4_{k}}(n)|_T). 
\]
By $T = U \cap U'$, we have the decomposition 
\[
\rho :  H^0(\P^4_k, \MO_{\P^4_{k}}(n)) \xrightarrow{\rho_1} 
H^0(U, \MO_{\P^4_{k}}(n)|_U) \xrightarrow{\rho_2}
H^0(T, \MO_{\P^4_{k}}(n)|_T). 
\]
For $\MO_U(c) := \MO_{\P^4_k}(c)|_U$, we have short exact sequences 
\[
0 \to \MO_{\P^4}(-a) \to \MO_{\P^4} \to \MO_U \to 0 
\]
\[
0 \to \MO_U(-b) \to \MO_{U} \to \MO_{T} \to 0. 
\]
In particular, $\rho_1$ is surjective. 
By the upper exact sequence, 
$H^1(U, \MO_U(m))=0$ for every $m \in \Z$. 
This, together with the lower exact sequence, implies that $\rho_2$ is surjective. 
\end{proof}

\section{Appendix: 
algebraicity of $H^1(X, \Omega_X^1)$}\label{s-rat3-Dol}


The purpose of this section is to prove that the $k$-linear map  
\[
\dlog : \Pic\,X \otimes_{\Z} k \to H^1(X, \Omega_X^1)
\]
is an isomorphism for a smooth projective threefold $X$ 
which is rational or birational to a smooth cubic threefold. 
The proof is divided into the following two steps (I) and (II). 
\begin{enumerate}
    \item[(I)] $h^1(X, \Omega_X^1)=\rho(X)$ (Subsection \ref{ss1-rat3-Dol}). 
    \item[(II)] $\Pic\,X \otimes_{\Z} k \to H^1(X, \Omega_X^1)$ is injective (Subsection \ref{ss2-rat3-Dol}). 
\end{enumerate}
The main tool for (I) is the Crystalline cohomology. 
The injectivity in (II) will be settled by an argument due to van der Geer--Katsura \cite{vdGK03}. 


\subsection{$h^{1, 1}(X) = \rho(X)$}\label{ss1-rat3-Dol}



\subsubsection{Dolbeault cohomologies vs de Rham cohomologies}

\begin{lem}\label{l-sm-RET1}
Let $f: X \to Y$ be a projective birational morphism of smooth varieties. 
Then $f_*\Omega_X^i = \Omega_Y^i$.
\end{lem}

\begin{proof}
The assertion hold because the composition of the pullback map $\Omega_Y^i\to f_*\Omega_X^i$ and the restriction $f_*\Omega_X^i \hookrightarrow (f_*\Omega_X^i)^{**} =\Omega_Y^i$ is an isomorphism, 
where $(-)^{**}$ denotes the double dual.
\end{proof}

\begin{lem}\label{l-sm-RET}
Let $X$ be a smooth projective threefold which is rational or birational to a smooth cubic threefold. 
Then the following hold. 
\begin{enumerate}
    \item $H^j(X, \MO_X)=0$ for every $j>0$. 
    \item $H^0(X, \Omega_X^i) =0$ for every $i>0$. 
\end{enumerate}
\end{lem}

\begin{proof}
By \cite[Theorem 1]{CR11} and Lemma \ref{l-sm-RET1}, 
we may assume that $X = \P^3$ or $X \subset \P^4$ is a smooth cubic threefold. 
If $X=\P^3$, then the assertion is well known. 
If $X \subset \P^4$ is a smooth cubic threefold, then (1) holds by 
\[
0 \to \MO_{\P^4}(-3) \to \MO_{\P^4} \to \MO_X \to 0. 
\]
The assertion (2) follows from Theorem \ref{thm:CI}. 
\end{proof}



\begin{lem}\label{l-CP}
Let $X$ and $Y$ be smooth projective threefolds such that 
$X$ and $Y$ are birational. 
Then there exists  
a birational morphism $f: Z \to X$ 
from a smooth projective threefold $Z$ which is obtained by a sequence of blowups with smooth centres 
\[
Z = Y_{\ell} \xrightarrow{g_{\ell}} Y_{\ell -1}
\xrightarrow{g_{\ell-1}}  \cdots \xrightarrow{g_1} Y_0 := Y, 
\]
i.e., each $g_i : Y_i \to Y_{i-1}$ is a blowup along a point or a smooth curve. 
\end{lem}

\begin{proof}
The assertion is a direct consequence of \cite[Proposition 4.2]{CP08}. 
\end{proof}

\begin{prop}\label{p-rat3-Dol}
Let $X$ be a smooth projective  threefold 
 which is rational or birational to a smooth cubic threefold. 
Then 
\[
H^1(X, \Omega_X^1) \simeq H^2(X, \Omega_X^{\bullet}). 
\]
\end{prop}

\begin{proof}
By the spectral sequence 
\[
E_1^{i, j} :=H^j(X, \Omega_X^i) \Rightarrow H^{i+j}(X, \Omega_X^{\bullet}) =: E^{i+j}
\]
and $E^{i, 0}_1 = E^{0, j}_1=0$  for all $i>0$ and all $j>0$ (Lemma \ref{l-sm-RET}), 
it is enough to show (1) and (2) below. 
\begin{enumerate}
    \item $E_2^{1, 1} = E_3^{1, 1} = E_4^{1, 1} = \cdots $. 
    \item $d:  E_1^{1, 1} =H^1(X, \Omega_X^1)  \to 
    H^1(X, \Omega_X^2) = E_1^{2, 1}$ is zero. 
\end{enumerate}
The assertion (1) follows from the definition of spectral sequences and $E_2^{3, 0} = 0$. 

Let us show (2). 
If $X=\P^3$ or $X \subset \P^4$ is a smooth cubic threefold, 
then (2) follows from the fact that 
the spectral sequence on $X$ generates at $E_1$ 
(Theorem \ref{thm:CI}). 
If $X =:X_0 \to X_1 \to \cdots \to X_n$ is a sequence of blowups with smooth centres, and 
$X_n=\P^3$ or $X_n \subset \P^4$ is a smooth cubic threefold, 
then the assertion follows from \cite[Corollary 4.9(1)]{AZ17}. 

Let us reduce the general case to this case. 
We can find a birational morphism $f: Y \to X$, where $Y$ is obtained by a sequence of blowups with smooth centres starting from $Z$, 
where $Z =\P^3$ or $Z$ is a smooth cubic threefold (Lemma \ref{l-CP}). 
Consider the following commutative diagram: 
\[
\begin{tikzcd}
 H^1(Y, \Omega_Y^1)  \arrow[r, "d_Y=0"] &    H^1(Y, \Omega_Y^2)\\
 H^1(X, \Omega_X^1)  \arrow[r, "d"] \arrow[u, "f^*"] & 
    H^1(X, \Omega_X^2), \arrow[u, "f^*"]
\end{tikzcd}
\]
where each vertical arrow is the natural map. 
Here each $f^*$ is injective by $f_*\Omega_Y^r = \Omega_X^r$ (Lemma \ref{l-sm-RET1}) and the Leray spectral sequence: 
\[
H^i(X, R^jf_*\Omega_Y^r) \Rightarrow H^{i+j}(Y, \Omega_Y^r), 
\]
where $r \in \{1, 2\}$. Thus (2) holds. 
    \qedhere

\end{proof}

\subsubsection{Crystalline cohomology vs Picard number}

For foundations on crystalline cohomologies, we refer to \cite{Ber74} and \cite{BO78}.

\begin{nasi}[Basic properties on Crystalline cohomology]\label{n-cris-basic}
Let $X$ be a smooth projective variety. Set $W := W(k)$. 
We here summarise some known results. 
\begin{enumerate}
    \item $H^i_{\cris}(X/W)$ is a finitely generated $W$-module for every $i \in \Z$. 
    \item $H^i_{\cris}(X/W)=0$ when $i <0$ or $i > 2\dim X$. 
    \item $H^0_{\cris}(X/W) \simeq W$ and $H^{2\dim X}_{\cris}(X/W) \simeq W$. 
    \item $H^1_{\cris}(X/W)$ is a free $W$-module of rank $\leq 2g$ for $g := \dim H^1(X, \MO_X)$ 
    \cite[Proposition 3.11 and Remarque 3.11.2 in page 618, Proposition 7.3.6(b)(ii) in page 656]{Ill79}. 
\end{enumerate}
\end{nasi}

\begin{prop}\label{p-rat3-pic-cris}
Let $X$ be a smooth projective rationally chain connected variety (e.g., a smooth Fano variety). 
Then 
\[
\rank_W H^2_{\cris}(X/W) = \rho(X)
\]
for $W:=W(k)$. 
\end{prop}


\begin{proof}
See \cite[Corollary 3.3 and Proposition 4.2]{GJ18}. 
\end{proof}

\subsubsection{Crystalline cohomology vs de Rham Cohomology}



\begin{prop}\label{p-cris-birat}
Let $f : X \to Y$ be a birational morphism of smooth projective threefolds. 
Fix $i \in \Z$. 
Then the following hold. 
\begin{enumerate}
    \item $f^*:  H^i_{\cris}(Y/W) \to H^i_{\cris}(X/W)$ is injective. 
    \item $H^i_{\cris}(X/W)$ is torsion free if and only if $H^i_{\cris}(Y/W)$ is torsion free. 
\end{enumerate}
\end{prop}

\begin{proof}
We first treat the case when $f$ is a blowup along a smooth subvariety $Z$. 
Set $c := {\rm codim}_Y\,Z \in \{2, 3\}$. 
By \cite[Corollary 4.8]{AZ17}, we obtain 
\[
H^i_{\Cris}(X/W) \simeq H^i_{\Cris}(Y/W) \oplus 
\bigoplus_{s=1}^{c-1} H^{i-2s}_{\cris}(Z/W). 
\]
Then (1) holds. 
The assertion (2) follows from the fact that $H^j_{\cris}(Z/W)$ is torsion free for every $j$ (\ref{n-cris-basic}). 
This completes the proof for the case when  $f$ is a blowup along a smooth subvariety. 
In particular, the assertion holds if $f$ is a sequence of blowups along smooth subvarieties.

Let us reduce the general case to the case treated in the previous paragraph. 
There exist birational morphisms of smooth projective threefolds  
\[
h : Z \xrightarrow{g} X \xrightarrow{f} Y, 
\]
where $h$ is a sequence of blowups along smooth subvarieties (Lemma \ref{l-CP}). 
By the functoriality: 
\[
h^* : H^i_{\cris}(Y/W) \xrightarrow{f^*} H^i_{\cris}(X/W) \xrightarrow{g^*} H^i_{\cris}(Z/W), 
\]
the injectivity of $h$ implies the injectivity of $f^*$. 
Thus (1) holds. 

Let us show (2). 
If $H^i_{\cris}(X/W)$ is torsion free, 
then (1) implies that 
$H^i_{\cris}(Y/W)$ is torsion free. 
Conversely, assume that 
$H^i_{\cris}(Y/W)$ is torsion free. 
Then $H^i_{\cris}(Z/W)$ is torsion free by the previous paragragh. 
Since $g^*$ is injective by (1), also $H^i_{\cris}(X/W)$ is torsion free. 
Thus (2) holds. 
\qedhere



\end{proof}

\begin{thm}\label{t-cris-birat}
Let $X$ and $Y$ be smooth projective threefolds such that 
$X$ and $Y$ are birational. 
Fix $i \in \Z$. 
Then $H^i_{\cris}(X/W)$ is torsion free if and only if $H^i_{\cris}(Y/W)$ is torstion free. 
\end{thm}

\begin{proof}
The assertion follows from Proposition \ref{p-cris-birat}(2) 
by taking birational morphisms $Z \to X$ and $Z \to Y$ from a smooth projective threefold $Z$. 
\end{proof}

\begin{thm}\label{t-rat-cris-birat}
Let $X$ be a smooth projective threefold which is rational or birational to a smooth cubic threefold. 
Then the following hold. 
\begin{enumerate}
    \item $H^i_{\cris}(X/W)$ is torsion free for every $i \in \Z$. 
    \item $\dim_k H^i(X, \Omega_X^{\bullet}) = \rank_W H^i_{\cris}(X/W)$ for every $i \in \Z$. 
\end{enumerate}
\end{thm}

\begin{proof}
Let us show (1). 
By Theorem \ref{t-cris-birat}, we may assume that $X = \P^3$ or $X \subset \P^4$ is a smooth cubic threefold. 
In this case, the assertion follows from Theorem \ref{thm:CI}  and 
Proposition \ref{prop:hodge number to degeneration and torsion-freeness}. 
Thus (1) holds. 
The assertion (2) follows from (1) and 
the univesal coefficient exact sequence  
(\cite[Ch. VII. Remarque 1.1.11 b)]{Ber74} or \cite[Summary 7.26.2]{BO78}): 
\[
0 \to H^i_{\Cris}(X/W) \otimes_W k 
\to H^i(X, \Omega_X^{\bullet}) \to {\rm Tor}_1^W(H^{i+1}_{\Cris}(X/W), k) \to 0. 
\]
\end{proof}

We are ready to show the main theorem of this subsection. 

\begin{thm}\label{t-rat3-cohs}
Let $X$ be a smooth projective  threefold 
which is rational or birational to a smooth cubic threefold. 
Then 
\[
\dim_k H^1(X, \Omega_X^1) = 
\dim_k H^2(X, \Omega_X^{\bullet}) = \rank_W H^2_{\cris}(X/W)
= \rho(X). 
\]
\end{thm}

\begin{proof}
The first (resp. second, resp. third) equality follows from 
Proposition \ref{p-rat3-Dol} (resp. Theorem \ref{t-rat-cris-birat}, resp. Proposition \ref{p-rat3-pic-cris}). 

\end{proof}

\subsection{Injectivity of $\Pic X \otimes k \to H^1(X, \Omega_X^1)$}\label{ss2-rat3-Dol}



\begin{prop}\label{p-vdGK}
Let $X$ be a smooth projective variety. 
Assume that
\begin{enumerate}
\item[(A)] $H^0(X, \Omega_X^1) =0$, and  
\item[(B)] $H^0(X, B\Omega_X^2)=0$. 
\end{enumerate}
Then the induced group homomorphism
    \[
    \dlog : \Pic\,X \otimes_{\Z} \F_p \to H^1(X, \Omega_X^1)
    \]
    is injective. 
\end{prop}



\begin{proof}
By the following exact sequence in \'etale topology: 
\[
1 \to \G_m \xrightarrow{p} \G_m\to \G_m/\G_m^p \to 1, 
\]
we get the following induced injection: 
\[
\Pic\,X \otimes_{\Z} \F_p \hookrightarrow H_{\et}^1(X, \G_m/\G_m^p). 
\] 
By \cite[Ch.~0, Collolaire 2.1.18]{Ill79},  we have the following exact sequence in \'etale topology: 
\[
0 \to \G_m/\G_m^p
\xrightarrow{\dlog} Z\Omega_X^1 \xrightarrow{1-C} \Omega^1_X \to 0. 
\]
Hence 
the map $\dlog : \Pic\,X \otimes_{\Z} \F_p \to H^1(X, \Omega_X^1)$ 
is decomposed as follows: 
\[
 \Pic\,X \otimes_{\Z} \F_p \hookrightarrow 
H_{\et}^1(X, \G_m/\G_m^p) 
\overset{{\rm (A)}}{\hookrightarrow} H^1(X, Z\Omega_X^1) \xrightarrow{H^1(i)} 
H^1(X, \Omega_X^1), 
\]
where $i : Z\Omega_X^1 \hookrightarrow F_{*}\Omega_X^1$ denotes the natural inclusion. 
Here $H^1(i)$ is injective by (B) and  the following exact sequence: 
\[
0 \to Z\Omega_X^1 \xrightarrow{i} F_{*}\Omega_X^1 \xrightarrow{F_{*}d} B\Omega_X^2 \to 0. 
\]
Thus the assertion holds. 
\end{proof}

\begin{thm}\label{t-vdGK}
Let $X$ be a smooth projective variety. 
Assume that
\begin{enumerate}
\item[(A)] $H^0(X, \Omega_X^1)=0$, 
\item[(B)] $H^0(X, B\Omega_X^2)=0$, and 
\item[(C)] $H^1(X, B\Omega_X^1)=0$. 
\end{enumerate}
Then the induced $k$-linear map 
\begin{equation}\label{e1-vdGK}
    \dlog : \Pic\,X \otimes_{\Z} k \to H^1(X, \Omega_X^1) 
\end{equation}
    is injective. 
\end{thm}


Although the following proof is identical to that of \cite[Theorem 4.9]{vdGK03}, we include it for the sake of completeness. 

\begin{proof}
Suppose that the map (\ref{e1-vdGK}) not injective. 
Let us derive a contradiction. 
Recall that, via $\Pic\,X \simeq H^1(X, \MO_X^{\times})$, 
the map (\ref{e1-vdGK}) satisfies the following commutative diagram: 
\[
\begin{CD}
\displaystyle{ \varinjlim_{\mathfrak{U}} \left(\check{H}^1(\mathfrak{U}, \MO_X^{\times}) 
\otimes_{\Z} k\right)} @>>> \displaystyle{\varinjlim_{\mathfrak{U}} 
\check{H}^1(\mathfrak{U}, \Omega_X^1)}\\
@VV\simeq V @VV\simeq V\\
H^1(X, \MO_X^{\times}) \otimes_{\Z} k @>(\ref{e1-vdGK})>> H^1(X, \Omega_X^1), 
\end{CD}
\]
where $\mathfrak U$ runs over all the affine open covers of $X$ and we have used the following isomorphism: 
\[
\displaystyle{\left( \varinjlim_{\mathfrak{U}} \check{H}^1(\mathfrak{U}, \MO_X^{\times})\right) 
\otimes_{\Z} k}
\simeq 
\displaystyle{ \varinjlim_{\mathfrak{U}} \left(\check{H}^1(\mathfrak{U}, \MO_X^{\times}) 
\otimes_{\Z} k\right)}. 
\]
Since $\varinjlim_{\mathfrak U}$ preserves injectivity, 
we can find an affine open cover $\mathfrak{U} = (U_i)_{i \in I}$ of $X$ 
such that 
\begin{equation}\label{e2-vdGK}
    \dlog : \check{H}^1(\mathfrak{U}, \MO_X^{\times}) \otimes_{\Z} k \to \check{H}^1(\mathfrak{U}, \Omega_X^1) 
\end{equation}
is not injective. 
Then there exist $a_1, ..., a_{\ell} \in k$ and 
$\zeta_1, ..., \zeta_{\ell} \in \check{H}^1(\mathfrak{U}, \MO_X^{\times})$ 
such that 
\begin{enumerate}
\item $\sum_{\nu=1}^{\ell} a_{\nu} \zeta_{\nu} \neq 0$ in $\check{H}^1(\mathfrak{U}, \MO_X^{\times}) \otimes_{\Z} k$ (where $a_{\nu} \zeta_{\nu} := \zeta_{\nu} \otimes a_{\nu}$), and 
    \item $\sum_{\nu=1}^{\ell} a_{\nu} \dlog(\zeta_{\nu}) = 0$ in $\check{H}^1(\mathfrak{U}, \Omega_X^1)$.
\end{enumerate}
We choose $\ell$ to be the smallest integer satisfying the above conditions. 
We see that $\zeta_1, ..., \zeta_{\ell}$ 
are $\F_p$-linearly independent in $H^1(X, \MO_X^{\times}) \otimes_{\Z} \F_p$, as otherwise we could make $\ell$ 
smaller. 
By (1) and Proposition \ref{p-vdGK}, we may assume that $\ell \geq 2$, $a_1 =1$, and $a_2 \not\in \F_p$. 
For each $\nu \in \{1, ..., \ell\}$, 
the element $\zeta_{\nu} \in \check{H}^1(\mathfrak{U}, \MO_X^{\times})$ is represented by 
elements $\{f_{ij}^{(\nu)} \in \Gamma(U_{ij}, \MO_X^{\times})\}_{i, j \in I}$ 
satisfying the cocycle condition 
$(f_{ij}^{(\nu)}|_{U_{ijk}}) (f_{jk}^{(\nu)}|_{U_{ijk}})(f_{ki}^{(\nu)}|_{U_{ijk}}) =1$ ($i, j, k \in I$). 
By (2), we can find elements
$\omega_{i} \in \Gamma(U_{i}, \Omega_{X}^{1})$ such that
\begin{equation}\label{e3-vdGK}
\sum_{\nu=1}^{\ell} a_{\nu}\frac{df_{ij}^{(\nu)}}{f_{ij}^{(\nu)}} 
= \omega_{i}|_{U_{ij}} -\omega_{j}|_{U_{ij}} 
\qquad \text{in}\qquad \Gamma(U_{ij}, \Omega_X^1). 
\end{equation}
Recall that we have the following exact sequence: 
\begin{equation}\label{e4-vdGK}
0 \to B\Omega_X^1 \to Z\Omega_X^1 \xrightarrow{C} \Omega_X^1 \to 0. 
\end{equation}
Since the Cartier operator $C: Z\Omega_X^1 \to \Omega_X^1$ is surjective 
(\ref{e4-vdGK}), 
there exists an element
$\widetilde{\omega}_i \in \Gamma(U_i, Z\Omega_{X}^{1})$ such that
$C({\widetilde{\omega}}_{i}) = \omega_{i}$ for every $i \in I$. 
Recall that we have 
\begin{equation}\label{e5-vdGK}
C\left( a^{1/p}_{\nu}\frac{df_{ij}^{(\nu)}}{f_{ij}^{(\nu)}} \right) = 
a_{\nu}\frac{df_{ij}^{(\nu)}}{f_{ij}^{(\nu)}} 
\end{equation}
by \cite[(2.1.23) in page 518]{Ill79}. 
It follows from 
(\ref{e3-vdGK}), (\ref{e4-vdGK}), and  (\ref{e5-vdGK})
that there exist elements $\xi_{ij} \in \Gamma(U_{ij}, B\Omega_X^1)$ such that 
\[
\xi_{ij} + 
\sum_{\nu=1}^{\ell} a^{1/p}_{\nu}\frac{df_{ij}^{(\nu)}}{f_{ij}^{(\nu)}} 
= \widetilde{\omega}_{i}|_{U_{ij}} -\widetilde{\omega}_{j}|_{U_{ij}} 
\qquad \text{in}\qquad \Gamma(U_{ij}, \Omega_X^1). 
\]
Since $\left\{ \frac{df_{ij}^{(\nu)}}{f_{ij}^{(\nu)}}\right\}$ satisfies the cocycle condition, 
so does $\{ \xi_{ij} \}$, and hence we obtain $\{ \xi_{ij} \} \in \check{H}^{1}(\mathfrak U, B\Omega^1_{X})$. 
By $\check{H}^{1}(\mathfrak U, B\Omega^1_{X}) \simeq 
H^1(X, B\Omega_X^1) \overset{{\rm (C)}}{=} 0$, there exist 
elements $\xi_{i} \in \Gamma(U_{i}, B\Omega_X^1)$ such 
that $\xi_{ij}|_{U_{ij}} = \xi_{i} - \xi_{j}$. 
Therefore, we have
\begin{equation}\label{e6-vdGK}
\sum_{\nu=1}^{\ell} {a}^{1/p}_{\nu} \frac{df_{ij}^{(\nu)}}{f_{ij}^{(\nu)}} =
({\widetilde{\omega}}_{i} -\xi_i)|_{U_{ij}} - ({\widetilde{\omega}}_{j} -\xi_j)|_{U_{ij}} 
\qquad \text{in}\qquad \Gamma(U_{ij}, \Omega_X^1). 
\end{equation}
Subtracting (\ref{e6-vdGK}) from (\ref{e3-vdGK}), we get 
\begin{enumerate}
    \item[(2)'] $\sum_{\nu=1}^{\ell} b_{\nu} \dlog(\zeta_{\nu}) = 0$ in $\check{H}^1(\mathfrak{U}, \Omega_X^1)$. 
\end{enumerate}
for $b_{\nu} := a_{\nu} - a^{1/p}_{\nu} \in k$. 
In order to conclude a contradiction, 
it is sufficient to show 
\begin{enumerate}
\item[(1)'] $\sum_{\nu=1}^{\ell} b_{\nu} \zeta_{\nu} \neq 0$ in $\check{H}^1(\mathfrak{U}, \MO_X^{\times}) \otimes_{\Z} k$.
\end{enumerate}
Indeed, (1)' and (2)', together with $b_1 = 0$ and $b_2 \neq 0$, 
contradicts the minimality of $\ell$. 
The condition (1)' follows from $b_2 \neq 0$ and 
the fact that $\zeta_1, ..., \zeta_{\ell} \in \check{H}^1(\mathfrak{U}, \MO_X^{\times}) \otimes_{\Z} \F_p$ are $\F_p$-linearly independent. 
\end{proof}



\begin{thm}\label{t-rat3-Pic-Dol}
Let $X$ be a smooth projective  threefold which is rational or birational to a smooth cubic threefold. 
Then the induced $k$-linear map 
\[
\dlog : \Pic\,X \otimes_{\Z} k \to H^1(X, \Omega_X^1)
\]
is bijective. 
\end{thm}

\begin{proof}
By Lemma \ref{l-sm-RET}, (A)-(C) of  Theorem \ref{t-vdGK} hold, 
and hence $\dlog :  \Pic\,X \otimes_{\Z} k \to H^1(X, \Omega_X^1)$ is injective. 
Then this is bijective by Theorem \ref{t-rat3-cohs}. 
\qedhere


\end{proof}

\begin{rem}
Let $X$ be a smooth projective  threefold which is rational or birational to a smooth cubic threefold. 
Then we get $\Br(X)$ is $p$-torsion free by using the argument as in \cite[Subsection 2.4]{AV}. 
Then we see that the assumptions in \cite[Subsection 2.1]{AV} hold for $X$. 
This gives an alternative proof of Theorem \ref{t-rat3-Pic-Dol}. 
The authors learned this method by Emiliano Ambrosi. 
\end{rem}

\bibliographystyle{skalpha}
\bibliography{bibliography.bib}

\end{document}